\documentclass[11pt]{article}
\usepackage{graphicx, amsmath, amsfonts, amssymb, amsthm, fullpage, authblk, bbm, comment, authblk, enumitem, etoc, booktabs}
\usepackage{natbib}
\usepackage[colorlinks=true,citecolor=blue,linkcolor=blue,urlcolor=blue]{hyperref}
\usepackage[dvipsnames]{xcolor}
\usepackage{setspace}
\usepackage{algorithm}
\usepackage{algpseudocode}
\algrenewcommand\algorithmicrequire{\textbf{Input:}}
\algrenewcommand\algorithmicensure{\textbf{Output:}}

\usepackage{cleveref}
\crefname{assumption}{assumption}{assumptions}
\Crefname{assumption}{Assumption}{Assumptions}

\newcommand{\R}{\mathbb R}
\newcommand{\E}{\mathbb E}
\newcommand{\Prob}{\mathbb P}
\newcommand{\TV}{\operatorname{D}_{\mathrm{TV}}}
\newcommand{\KL}{\operatorname{D}_{\mathrm{KL}}}
\newcommand{\Ham}{\operatorname{D}_{\mathrm{Ham}}}
\newcommand{\ind}{\mathbbm{1}}
\newcommand{\tr}{\operatorname{tr}}

\newcommand{\dint}{\, \mathrm{d}}
\DeclareMathOperator{\var}{Var}
\DeclareMathOperator*{\argmin}{arg\,min}

\DeclareMathOperator*{\intr}{int}
\DeclareMathOperator*{\cov}{Cov}

\newtheorem{theorem}{Theorem}
\newtheorem{lemma}[theorem]{Lemma}
\newtheorem{proposition}[theorem]{Proposition}

\newtheorem{remark}{Remark}

\newtheorem{assumption}{Assumption}
\newtheorem{definition}{Definition}

\title{Minimax signal detection in regression under differential privacy}
\author{Elly K.~H.~Hung}
\affil{Department of Statistics, University of Warwick}

\begin{document}
\maketitle

\begin{abstract}
We establish minimax separation rates across sparsity levels for regression signal detection under $(\epsilon, \delta)$-central differential privacy, considering linear regression with an unknown error distribution and generalised linear models (GLMs) with an unknown intercept. Without sparsity restrictions, the separation rates for linear regression and GLMs with binary responses match the private Gaussian mean testing rate up to polylogarithmic factors. For GLMs with more than two possible response outcomes, additional information in the marginal response distribution permits faster detection of sufficiently small signals when the dimension exceeds the sample size $n$. Our upper bounds are based on a reduction to private mean testing, tests of dependence between the response and sparse covariate projections, and a response-only test that is adjusted for the unknown intercept. The resulting procedure adapts to unknown sparsity and matches our lower bounds up to polylogarithmic factors. We also construct polynomial-time tests that attain these rates when $\epsilon\gtrsim n^{-1/2}$, and compare their numerical performance with existing private procedures.
\end{abstract}

\addtocontents{toc}{\protect\setcounter{tocdepth}{0}}

\section{Introduction}
Determining whether a collection of covariates is associated with a response is a classical problem in statistical inference. In high-dimensional regression, this question can be formulated as the global hypothesis testing problem
\begin{equation*}
    H_0:\beta=0
    \qquad\text{against}\qquad
    H_1:\|\beta\|_2\ge r,
\end{equation*}
where $\beta\in\mathbb R^d$ is the regression coefficient and $r$ measures the signal strength. It is well known that detecting whether a signal is present can be substantially easier than estimating the regression coefficient itself. For example, in the dense isotropic Gaussian linear model with $d = o(n)$ and unit error variance, the minimax $\ell_2$ estimation error is of order $\sqrt{d/n}$, whereas the minimax separation rate for signal detection is of order $d^{1/4}/\sqrt n$ \citep{ingster2010detection}. 

The gap between testing and estimation makes signal detection particularly relevant in high-dimensional applications, many of which involve sensitive data such as genomic records or consumer transactions. This motivates studying how hypothesis testing can be carried out while protecting individual-level data. A leading framework for privacy is differential privacy \citep[DP,][]{dwork2006differential}, which has been widely used in both industry \citep[e.g.][]{ding2017telemetry} and government applications \citep[e.g.][]{abowd2022confidentiality}. At a high level, differential privacy requires the distribution of a released output to be stable under replacement of a single observation. 

Intuitively, stronger privacy guarantees require larger signals to reliably distinguish between the null and alternative hypotheses.  We study how the minimax separation rate depends on the sample size, dimension, sparsity level, and privacy parameters, as well as on the regression model and its nuisance parameters.

\subsection{Contributions}
We characterise, up to polylogarithmic factors, the minimax separation rates for private signal detection across sparsity levels in linear regression with an unknown error distribution and generalised linear models (GLMs) with an unknown intercept.
\begin{itemize}
    \item In the dense regime, the rates for linear regression and GLMs with two response support points match the private Gaussian mean testing rate established by \cite{narayanan2022private}. We attain this rate through a transformation that reduces regression testing to mean testing and accommodates the unknown nuisance parameters. In particular, for linear regression, no moment, symmetry or smoothness assumptions on the errors are needed.
    \item Our lower bounds identify when nuisance adjustment can absorb changes in the response marginal caused by a regression signal. This affects both the non-private and private detection rates. For GLMs with more than two response support points, the response marginal retains second-order information about the signal despite the unknown intercept, allowing faster detection of sufficiently small signals when $d\gtrsim n$. We therefore construct a response-only test that targets this information.
    \item To exploit sparsity, we test for non-zero covariances between bounded transformations of the response and sparse covariate projections. Combining the sparse covariance test with the mean-reduction and response-only tests yields a procedure that adapts to the unknown sparsity level. As a byproduct, we characterise the minimax separation rates for private sparse mean testing of sub-Gaussian vectors up to polylogarithmic factors. To our knowledge, this characterisation under central DP is new and we present it in \Cref{app:sparse-mean} as a result of independent interest.
\end{itemize}

We also construct polynomial-time tests that match the minimax rates up to polylogarithmic factors when $\epsilon\gtrsim n^{-1/2}$, and provide numerical experiments comparing their performance with existing private procedures.

\subsection{Related work}
\noindent \textbf{Non-private regression testing.} Without privacy constraints, global hypothesis testing in high-dimensional regression has been studied extensively, so we only list a few examples. \cite{ingster2010detection} established the asymptotic detection boundaries for Gaussian linear regression across sparsity levels, distinguishing known and unknown error variance, while \cite{carpentier2022estimation} obtained non-asymptotic separation rates under sub-Gaussian errors with unknown variance. \cite{mukherjee2015hypothesis} derived analogous detection boundaries in the context of binary regression. For GLMs, \cite{goeman2011testing} developed asymptotic theory for a low dimensional null against a high dimensional alternative under canonical link functions, allowing for nuisance parameters; this framework was extended in \cite{guo2016tests} to accommodate diverging dimension and a broader class of link functions. 

\medskip
\noindent \textbf{Private regression inference.}
 Early work developed private procedures for ordinary least squares (OLS), including confidence intervals and approximations to classical $t$-statistics \citep{sheffet2017differentially, barrientos2019significance}. \cite{Alabi2023Ftest} constructed differentially private $F$-tests for testing linear relationships and mixtures, while \citet{pena2024uncertainty} studied likelihood-based hypothesis testing and model comparison for Gaussian linear regression. Bias-corrected coefficient estimates and asymptotic confidence intervals were considered in \cite{lin2026differentially}. For M-estimators more generally, \cite{avellamedina2021privacy} developed asymptotic theory for private analogues of the Wald, score and likelihood ratio tests. Bayesian approaches for private inference have been developed for linear regression \citep{bernstein2019bayesian} and for logistic and Poisson regression \citep{kulkarni2021glm}. \cite{hung2026differentially} consider private hypothesis testing with censored survival data in the Cox regression model. For high-dimensional sparse linear regression, \cite{cai2026highdimensional} developed differentially private methods for debiased inference on regression coefficients and multiple testing with false discovery rate control. These contributions primarily study private analogues of classical inferential procedures and their validity. Recent work also studies minimax support recovery in sparse linear regression with bounded covariates \citep{gu2026learning}.

\medskip
\noindent \textbf{Optimality in private hypothesis testing.}
Another line of work develops finite-sample private tests together with optimality results. \cite{canonne2019structure} characterised the sample complexity of private testing between simple hypotheses, while \cite{awan2018binomial} derived uniformly most powerful private tests for binomial data. For distribution testing, \cite{acharya2018discretedists} established optimal rates for identity testing of discrete distributions and closeness testing in the sparse regime. \cite{canonne2020private} developed identity tests for high-dimensional product Bernoulli distributions and  Gaussians with known covariance, and the  Gaussian sample complexity bounds were sharpened by \cite{narayanan2022private}. More recently, private nonparametric testing has been studied through minimax analysis of permutation tests \citep{kim2026permutation} and smooth two-sample testing \citep{kim2026sharp}. The central DP rate for the white-noise-with-drift model follows as a special case of the federated DP framework of \cite{cai2024federated}. In the local model of DP, minimax testing has also been studied for simple hypotheses \citep[e.g.][]{pensia2023simple} and for goodness-of-fit testing in discrete and continuous settings \citep[e.g.][]{lamweil2022minimax}. 

Hence, to the best of our knowledge, minimax separation rates for private signal detection in regression remain largely unresolved.

\subsection{Organisation and notation}
Section 2 introduces the models and states the main results. Sections 3 and 4 develop the lower bounds and corresponding tests. Section 5 presents polynomial-time implementations and numerical experiments, and Section 6 discusses future directions. Additional technical details are deferred to the appendices.

We use the following notation. For a positive integer $n$, we write $[n]:=\{1, 2, \dots, n\}$. For $x \in \mathbb{R}$ and $a\leq b$, we define the clipping operator $[x]_{a}^b := \min\{b, \max\{a, x\}\}$. For a vector $v \in \mathbb{R}^d$, we write $\|v\|_2$ for its Euclidean norm, and write $\|v\|_0 := | \{j \in [d]: v_j \neq 0\}|$ for its number of non-zero coordinates. We say a random variable $X \in \mathbb{R}$ is sub-Gaussian if $\|X\|_{\psi_2} := \inf\{c>0: \mathbb{E}[\exp(X^2/c^2)] \leq 2\} < \infty$, and we say a random vector $X \in \mathbb{R}^d$ is sub-Gaussian if $\|X\|_{\psi_2} := \sup_{u: \|u\|_2 \leq 1} \|u^\top X\|_{\psi_2} < \infty$. For a distribution $P$, we write $P^{\otimes n}$ to denote its $n$-fold product. We use $C,c,c_1,\ldots$ to denote positive constants that may depend only on the fixed constants in our assumptions but not on $n,d,s^*,\epsilon$ or $\delta$. Their values may change from line to line unless otherwise specified. For two non-negative sequences $\{a_n\}$ and $\{b_n\}$, we write $a_n \lesssim  b_n$ if there exists a constant $C>0$ such that $a_n \leq  Cb_n$ for all $n$, and we write $a_n \asymp b_n$ if both $a_n \lesssim b_n$ and $b_n \lesssim a_n$. We use $\lesssim_{\log}$ and $\asymp_{\log}$ for the corresponding comparisons up to polylogarithmic factors. We also write $a_n=O(b_n)$ if $|a_n|\lesssim b_n$, and $a_n=o(b_n)$ if $a_n/b_n \to 0$ as $n\to\infty$, where $b_n>0$. We denote by $\operatorname{Lap}(b)$ the centred Laplace distribution with scale $b$, and use $\TV$, $\KL$ and $\Ham$ for total variation distance, Kullback--Leibler divergence, and Hamming distance, respectively.

\section{Problem setup and main results}
We observe $n$ independent and identically distributed (i.i.d.) pairs $\{(X_i, Y_i)\}_{i=1}^n$ taking values in $\mathbb{R}^d \times \mathcal{Y}$. The conditional distribution of $Y$ given $X=x$ is specified by the linear predictor $x^\top \beta$ and an unknown nuisance parameter $\nu \in \mathcal{N}$. Here, $\beta\in\mathbb R^d$ is the regression coefficient of interest, and $\nu$ may differ under the null and alternative hypotheses. We consider testing 
\begin{equation}
    \mathrm{H_0}: \beta=0 \qquad \mathrm{against} \qquad \mathrm{H_1}(r, R): r \leq \|\beta\|_2 \leq R
    \label{eq:hypotheses-def}
\end{equation}
for $0<r\leq R \leq \infty$. Taking $R=r_0$, where $r_0>0$ is a sufficiently small model-dependent constant, restricts the alternatives to a neighbourhood of zero. 

\begin{assumption}[Covariate distribution]
    Fix constants $0<\underline{\lambda}\leq\overline{\lambda}$ and $C_X\geq\sqrt{8\underline{\lambda}/3}$. We assume the covariate distribution satisfies $\mathbb{E}[X]=0$,  $ \underline{\lambda} I_d\preceq\E[XX^\top]\preceq\overline \lambda I_d$, and $\|X\|_{\psi_2} \leq C_X$. 
\label{assp:covariate-dist}
\end{assumption}

The restriction on $C_X$ ensures that the covariate class contains $N(0,\underline{\lambda}I_d)$, which we use in the lower-bound construction. 

\subsection{Model classes}
\label{subsec:models}
We focus on two classes of regression models, which have different nuisance parameters. 

\medskip 
\noindent \textbf{Linear regression.} We consider the model 
\begin{equation}
    Y = X^\top \beta + \xi, \qquad \xi \perp X
    \label{eq:linear_regression_model}
\end{equation}
where the distribution of $\xi$ is unknown.

\begin{assumption}[Linear regression error]
\label{assp:linear-error}
There exist fixed constants $M_\xi>0$ and $p_\xi\in(0,1)$
such that $\Prob(|\xi|\leq M_\xi)\geq p_\xi$.
\end{assumption}
\Cref{assp:linear-error} only requires that errors have bounded magnitude with a fixed positive probability. In particular, it imposes no smoothness, symmetry, or tail assumptions on the error distribution.

\medskip
\noindent \textbf{Generalised linear models.}
We consider GLMs based on a fixed non-degenerate, one-parameter natural exponential family, with open natural parameter space $\Theta$. Conditional on its covariates $X=x$, the response distribution is given by 
\begin{equation}
p_{\beta,\gamma}(y\mid x)=h(y)\exp\left\{y\vartheta(\gamma+x^\top\beta)-b \left(\vartheta(\gamma+x^\top\beta)\right)\right\}.
\label{eq:GLM-model-def}
\end{equation}
Here, $h:\mathcal{Y} \rightarrow \mathbb{R}_{\geq 0}$ is the base density, $\vartheta: \mathbb{R} \rightarrow \Theta$ maps the linear predictor to the natural parameter, $b$ is the log-partition function, and $\gamma \in \Gamma$ is an unknown nuisance intercept. 

\begin{assumption}[GLM regularity] We assume that $\vartheta$ is either non-decreasing or non-increasing on $\mathbb{R}$. The intercept parameter space $\Gamma$ is a compact interval with non-empty interior, and there exists an open interval $U\supset\Gamma$ such that $\vartheta\in C^3(U)$ and $\inf_{a\in U}|\vartheta'(a)|>0$. 
\label{assp:glm-regularity}
\end{assumption}

The conditions in \Cref{assp:glm-regularity} hold for Bernoulli regression with logit or probit links, and Poisson and exponential regression with log links, among others. The global monotonicity condition ensures that the effect of the linear predictor does not reverse outside the intercept region. The remaining conditions provide the smoothness needed for local likelihood expansions and ensure that the response distribution remains locally sensitive to changes in the linear predictor.

\subsection{Preliminaries on differential privacy}

\begin{definition}[$(\epsilon, \delta)$-differential privacy] \label{def-cdp} 
For $\epsilon > 0$ and $\delta \geq 0$, a privacy mechanism $M$ is $(\epsilon, \delta)$-differentially private if it is a conditional distribution on a space $\mathcal{R}$ such that
\begin{equation*}
M(R \in A \mid D) \leq e^\epsilon M(R \in A \mid D') + \delta
\end{equation*}
for any measurable set $A$ in the $\sigma$-algebra of $\mathcal{R}$ and any pair of datasets $(D, D')$ differing in at most one record.
\end{definition}  

A test $\phi$ is a (randomised) decision rule taking values in $\{0, 1\}$, where $\phi=1$ denotes rejection of the null hypothesis. Our decision rules can be viewed as post-processing a private test statistic, which preserves differential privacy.

\begin{lemma}[Post-processing]
Let $M: \mathcal{D} \rightarrow \mathcal{R}$ be $(\epsilon,\delta)$-DP, and let $f:\mathcal{R} \to \mathcal S$ be any, possibly randomised, mapping independent of the data conditional on $M(D)$. Then $f\circ M$ is $(\epsilon,\delta)$-DP. 
\label{lemma:post-processing}
\end{lemma}

Many of our tests will be based on privatising test statistics using the Laplace mechanism.

\begin{lemma}[Laplace mechanism]
Let $\Delta_T
:= \sup_{D, D': \Ham(D, D')=1}|T(D)-T(D')|$ denote the global sensitivity of $T: \mathcal{D} \rightarrow \mathbb{R}$. Then 
\begin{equation*}
    M(D)=T(D)+W,
    \qquad  W\sim\operatorname{Lap} \left(\frac{\Delta}{\epsilon}\right),
\end{equation*}
 is $(\epsilon, 0)$-DP provided that $\Delta \geq \Delta_T$. By \Cref{lemma:post-processing}, for every $\tau\in\mathbb R$, $\phi(D):= \mathbbm 1\{M(D)>\tau\}$ is an $(\epsilon, 0)$-DP test. 
\label{lemma:laplace-test}
\end{lemma}

\subsection{Minimax separation rates}
For a class of joint distributions $\mathcal P=\{P_{\beta,\nu}:\beta\in\mathbb R^d,\nu\in\mathcal N\}$, let $\nu$ collect all nuisance quantities, including the covariate distribution. For $1\leq s^*\leq d$ and $0<r<R\leq\infty$, define
\begin{equation*}
 \mathcal P_0:=\{P_{0,\nu}:\nu\in\mathcal N\},  \qquad
 \mathcal P_{1}(s^*, r,R):=
 \{P_{\beta,\nu}:r\leq\|\beta\|_2\leq R,\ \|\beta\|_0\leq s^*,\ \nu\in\mathcal N\}.
\end{equation*}

We define the minimax separation radius as 
\begin{align*}
r^*(n,d,s^*, \epsilon,\delta,\mathcal{P}, R):=\inf\Big\{r \in (0, R) \;:\;
&\exists\,(\epsilon,\delta)\text{-DP test } \phi \text{ such that} \\
&\sup_{P\in\mathcal{P}_0} P^{\otimes n}(\phi=1) +
\sup_{P\in\mathcal{P}_{1}(s^*, r, R)} P^{\otimes n}(\phi=0) \leq \frac{1}{3} \Big\},
\end{align*}
where the probabilities are taken over the randomness of the samples and the private test. The threshold $1/3$ can be replaced with any fixed constant in $(0, 1)$, which would lead to the same rates up to constant factors.

Let $\mathcal P_{\mathrm{Lin}}$ denote the class of joint distributions of $(X,Y)$ under the linear regression model. Here, $\beta$, the covariate distribution, and the error distribution vary over those satisfying \Cref{assp:covariate-dist} and \Cref{assp:linear-error}, respectively. For a fixed GLM response family $\mathcal F$ satisfying \Cref{assp:glm-regularity}, let $\mathcal{P}_{\mathrm{GLM}}(\mathcal{F})$ denote the corresponding class of joint distributions. We likewise allow $\beta$, the intercept $\gamma$, and the covariate distribution to vary, with the latter satisfying \Cref{assp:covariate-dist}. To state the rates, let 
\begin{equation}
    \varrho_{s^*, D} := \min \left \{ \frac{D^{1/4}}{n^{1/2}} +\frac{D^{1/4}}{n^{3/4}\epsilon^{1/2}}, \sqrt{\frac{s^*}{n}} \right\}  +\frac{1}{n\epsilon}
    \quad \mathrm{and} \quad
    \ell_{n,d,\epsilon} :=1+\log\{(n+1)(d+1)(1+\epsilon^{-1})\}.
    \label{eq:separation-notation}
\end{equation}
The following theorem gives the minimax separation rate for linear regression and GLMs with two-point response support, and the local minimax separation rate for GLMs with larger response support.

\begin{theorem}
\label{thm:main-result}
Fix $\epsilon>0$ and $0\leq\delta<1$ with
$\delta\leq c_\delta \epsilon$, where $c_\delta>0$ is an absolute constant. For the GLM statements, fix a response family $\mathcal{F}$ satisfying \Cref{assp:glm-regularity}.

For $\mathcal P=\mathcal P_{\mathrm{Lin}}$, or
$\mathcal P=\mathcal P_{\mathrm{GLM}}(\mathcal F)$ with two response support points, there exists a constant $c_P>0$ such that if $\ell_{n,d,\epsilon}^3 \varrho_{s^*, d} \leq c_P$, then
\begin{equation}
    \varrho_{s^*, d} \lesssim r^*(n,d, s^*,\epsilon,\delta,\mathcal P,\infty) \lesssim \ell_{n, d, \epsilon}^3 \,\varrho_{s^*, d}.
    \label{eq:lower-linear}
\end{equation}

For $\mathcal P_{\mathrm{GLM}}(\mathcal F)$ with response support containing more than two points, there exist constants $r_0, c_P>0$ such that, if
$\ell_{n,d,\epsilon}^{3}\varrho_{s^*, d\wedge n}\leq c_P$, then
\begin{equation}
    \varrho_{s^*, d \wedge n} \lesssim r^*(n,d, s^*,\epsilon,\delta,\mathcal {P}_{\mathrm{GLM}}(\mathcal{F}), r_0) \lesssim \ell_{n, d, \epsilon}^3\, \varrho_{s^*, d \wedge n}.
    \label{eq:lower-GLM}
\end{equation}
\end{theorem}

The upper bounds are attained by pure-DP $(\delta=0)$ tests, while the lower bounds also apply for approximate DP with $\delta\leq c_\delta\epsilon$. The minimum selects the better of two bounds, one depending on dimension and the other on sparsity. The first combines the non-private dense rate with a dimension-dependent privacy contribution reflecting uncertainty about the signal direction. The $(n\epsilon)^{-1}$ term is a dimension-free privacy cost already present in two-point hypothesis testing. In the sparse regime, the privacy contribution is independent of both dimension and sparsity, up to logarithmic factors. This contrasts with both the $d^{1/4}$ dependence in the dense directional privacy term and the $s^*/(n\epsilon)$ privacy cost for pure-DP $\ell_2$ estimation of sparse Gaussian means \citep[e.g.][]{georgiev2022privacy}.

For GLMs with more than two points in their response distribution, the effective dimension is $d\wedge n$ rather than $d$. This can improve the detection rate when $d\gtrsim n$, and ignoring logarithmic factors, sparsity does not further improve this bound when $s^*\gtrsim\max\{\sqrt n,\epsilon^{-1}\}$. The next section explains this distinction through the information remaining in the response marginal after accounting for the nuisance parameter.

\section{Lower bounds}
\label{sec:lower-bounds}
In this section, we establish the lower bounds in \Cref{thm:main-result}. We begin with a general bound for private testing based on divergences between distributions. 

\begin{proposition}
    Fix $\epsilon, \delta \geq 0$. Let $\mathcal{X}$ be a Polish space with its Borel $\sigma$-field, and let $\mathcal{P}_0$ and $\mathcal{P}_1$ denote two classes of distributions on $\mathcal{X}$. Let $P_0\in\mathcal{P}_0$ and $P_1\in\mathcal{P}_1$. Let $\pi$ be a probability measure on a parameter space such that $P_\beta\in\mathcal{P}_1$ for every $\beta\in\operatorname{supp}(\pi)$, and define $\bar{P}^{(n)}:=\int P_\beta^{\otimes n}\,\pi(\mathrm d\beta)$. Suppose there exist $\alpha, \tau, \omega \geq 0$ such that the following conditions hold:
    \begin{enumerate}
    \item [(a)] we have $\TV(P_0, P_1) \leq \alpha$. 
    \item [(b)]  for each $\beta \in \mathrm{supp}(\pi)$, there exists a distribution $Q_\beta$ such that $\TV(P_\beta, Q_\beta) \leq \tau$.
    \item[(c)] we can bound $\KL(\bar{Q}^{(n)} \parallel P_0^{\otimes n}) \leq n\omega^2$, where $\bar{Q}^{(n)} := \int Q_\beta^{\otimes n} \pi (\mathrm{d}\beta)$. 
\end{enumerate}
Then for any $(\epsilon, \delta)$-DP test $\phi$, we can lower bound the testing error by 
\begin{align}
    &\sup_{P\in\mathcal{P}_0} P^{\otimes n}(\phi=1) + \sup_{P\in\mathcal{P}_1} P^{\otimes n}(\phi=0) \notag \\
    &\qquad \geq \frac{1}{2} \max \Big\{1 - \sqrt{n}\omega - n \tau,  0.9\exp \left(-10 \epsilon n \{\alpha \wedge (\omega + \tau)\} \right) - 10\delta n (\alpha \wedge (\omega + \tau)) \Big\}.
    \label{eq:general_lower}
\end{align}
\label{prop:general_lower_bound}
\end{proposition}
The two terms in \eqref{eq:general_lower} give the non-private and privacy-constrained minimax lower bounds for testing, respectively. A related transport inequality for $(\epsilon,0)$-DP tests was recently established by \citet{kim2026sharp}. Our formulation differs in two ways. First, to accommodate approximate ($\delta > 0)$ DP, we follow \cite{acharya2021differentially} in applying group privacy constraints to datasets generated from a coupling of the null and alternative mixture distributions. This allows us to control the difference in rejection probabilities through the Hamming distance between these datasets. We then use the transport inequality of \cite{marton1986simple} to bound the expected Hamming distance in terms of the KL divergence. Second, when a direct KL comparison between the original mixture and the null is difficult to bound, we introduce proxy distributions for the mixture components. The proof is given in \Cref{app:general_lower_bound}. Applying \Cref{prop:general_lower_bound} to our regression setting then requires suitable alternatives and proxy distributions.

\subsection{Role of the nuisance intercept}
The different rates in \Cref{thm:main-result} reflect how much of the change in the response marginal distribution can be explained by perturbing the nuisance parameter. For GLMs, averaging over the centred covariate distribution removes the first-order marginal perturbation, leaving a leading change of order $\|\beta\|_2^2$. In the GLM \eqref{eq:GLM-model-def}, a perturbation of the intercept on the same scale can absorb part of this remaining signal.

For $a\in\Gamma$, let $P_a$ denote the response distribution at linear predictor $a$, with density or probability mass function $p_a$, and write $\E_a$ for the expectation under $P_a$. Let $S_a := (\partial_a p_a)/p_a$ and $H_a := (\partial_a^2 p_a)/p_a$ denote the first- and second-order likelihood directions for perturbations of the linear predictor. Let $R_a$ be the residual from the orthogonal projection of $H_a$ onto the span of $S_a$ in $L_2(P_a)$, so that
\begin{equation}
    \E_a[R_a^2] = \inf_{c \in \mathbb{R}} \E_a \left[(H_a + 2 c S_a)^2 \right].
    \label{eq:projection-resdiual-def}
\end{equation}
We provide an explicit form for $R_a$ in \Cref{lemma:projection}. 

\begin{lemma} Under \Cref{assp:glm-regularity}, for every $a\in \Gamma$, we have that $\E_a[R_a^2]=0$ if and only if $|\mathrm{supp}(Y)|=2$. For families with larger response support, $\E_a[R_a^2]$ is bounded away from zero uniformly on $\Gamma$. 
\label{lemma:two-point-iff}
\end{lemma}

\Cref{lemma:two-point-iff} explains why the separation rates in \Cref{thm:main-result} differ across response families. When the response takes only two possible values, its distribution is determined by the probability of either value. The effect of a small signal on this marginal probability can be cancelled by adjusting the intercept. Similarly, in the Gaussian linear model, the error variance can be adjusted to match the response marginal. These least favourable constructions leave the dependence between $X$ and $Y$ as the only source of information about the signal. For GLMs with larger response supports, the response marginal retains information that adjusting the intercept cannot remove. For example, in exponential regression, variation in $X^\top\beta$ induces a mixture of response distributions with different rates. This mixture has greater variance than an exponential distribution with the same mean. Under the reference law $Y\sim\mathrm{Exp}(\exp(a))$, the corresponding residual is $R_a(Y)=\exp(2a)Y^2-4\exp(a)Y+2$, so $\E_a[R_a^2]=4$. Gaussian linear regression with known error variance also belongs to this class, allowing for a faster detection rate in high dimensions. 

To apply \Cref{prop:general_lower_bound}, we construct mixtures over sparse coefficient vectors with randomly chosen supports and signs, and bound the divergence between the corresponding proxy mixture and a fixed null. The directional contribution to the divergence is governed by the overlap between two independent draws from this prior, capturing uncertainty about both the support and direction of the signal. The residual marginal contribution has leading order $n\E_a[R_a^2]r^4$ and does not depend on sparsity. Combining these bounds with the $1/(n\epsilon)$ lower bound from two-point hypothesis testing gives the lower bounds in \Cref{thm:main-result}. The full details can be found in \Cref{app:lower_bounds}. 

\section{Upper bounds}
\label{sec:upper_bounds}
Guided by the lower bounds, we use two sources of information about the regression signal: dependence between the covariates and response, and changes in the response marginal. 

To detect dependence, a natural starting point is $\E[XY]$ -- under \eqref{eq:linear_regression_model}, we have $\E[XY] = \E[X X^\top] \beta$ provided that $\E[|\xi|] < \infty$. However, our error class does not assume finite moments, and $XY$ may have much heavier tails than the covariates. We therefore replace $Y$ by $\tanh(Y)$. The bound $|\tanh| \leq 1$ ensures that $X\tanh(Y)$ inherits the sub-Gaussian concentration of the covariates, while the strict monotonicity of $\tanh$ ensures that the transformed mean preserves the regression signal, uniformly over the distributions in each model class. 

\begin{lemma}
\label{lemma:transform-choice}
\label{lemma:transformation_properties}
Suppose $(X,Y)\sim P$ for some $P\in\mathcal P_{\mathrm{Lin}}$ or $P\in\mathcal P_{\mathrm{GLM}}(\mathcal F)$, where $\mathcal F$ is a fixed GLM response family satisfying \Cref{assp:glm-regularity}. Set $Z:=X\tanh(Y)$.
Then
\begin{equation*}
 \E[ZZ^\top]\preceq\E[XX^\top], \qquad \|\E[Z]\|_2 \leq \sqrt{\overline{\lambda}},
 \qquad \|Z-\E[Z]\|_{\psi_2}\lesssim C_X.
\end{equation*}
Under the null, $\E[Z]=0$, while there exists a fixed constant $c_\mu>0$ such that under every alternative,
\begin{equation}
 \|\E[Z]\|_2\geq c_\mu\min\{\|\beta\|_2,1\}.
 \label{eq:mean-to-regression}
\end{equation}
\end{lemma}

\Cref{lemma:transform-choice} reduces regression signal detection to private mean testing, giving the first test in \Cref{subsec:generalised_mean_test}. In \Cref{sec:sparse-cov}, we consider sparse alternatives by constructing an adaptive test based on empirical covariances between $\tanh(v^\top X)$ and $\tanh(Y)$ over sparse unit directions $v$. Finally, \Cref{sec:radial_test} constructs response-only tests for GLMs with more than two response support points, which target the second-order marginal change \eqref{eq:projection-resdiual-def} identified in the lower bounds.

\subsection{A mean-reduction test}
\label{subsec:generalised_mean_test}
We start by setting up some notation to define the mean test. For a transformed dataset $z_{1:n}:=(z_1, \dots, z_n) \in (\mathbb{R}^d)^n$, we define
\begin{equation*}
    u_i(z_{1:n}) := \sum_{j \neq i} z_i^\top z_j, \qquad Q_S(z_{1:n}) := \sum_{i, j \in S, \, i \neq j} z_i^\top z_j.
\end{equation*}
For $s, L>0$, let $\mathcal{D}_{n,s,L}$ be the set of datasets
$z_{1:n}\in(\mathbb R^d)^n$ satisfying
\begin{align}
    |u_i(z_{1:n})| &\leq L\{\sqrt{nd}+sn\}, && i\in[n],
    \label{eq:row-bound}\\
    |Q_S(z_{1:n})| &\leq L\{|S|\sqrt{|S|d}+|S|^2\}, && S\subseteq[n].
    \label{eq:subset-bound}
\end{align}
For the mean testing problem, a natural statistic is $Q_{[n]}(Z_{1:n}) := \sum_{i \neq j} Z_i^\top Z_j$, which is unbiased for $n(n-1)\|\mathbb{E}[Z]\|_2^2$. The central idea from \cite{narayanan2022private} is that although $Q_{[n]}$ has high global sensitivity, its variation is substantially smaller on typical datasets. More precisely, datasets in $\mathcal{D}_{n, s, L}$ satisfy the following stability property. 

\begin{lemma}
\label{lemma:restricted-sensitivity}
If $z_{1:n},z'_{1:n}\in \mathcal{D}_{n,s,L}$ differ in $K$ records, then 
\begin{equation*}
    |Q_{[n]}(z_{1:n})-Q_{[n]}(z'_{1:n})| \leq 6KL\{\sqrt{nd}+ns+K\}.
\end{equation*}
\end{lemma}
Comparing this bound with the expectation of $Q_{[n]}(Z_{1:n})$ under the alternative at scale $s$ already suggests the two privacy-dependent terms in \eqref{eq:lower-linear}: the $\sqrt{nd}$ term gives the dimension-dependent contribution, while the $ns$ term gives the dimension-free $1/(n\epsilon)$ contribution.

To apply this stability property to regression testing, we establish that the transformed observations $Z_{1:n}$ belong to $\mathcal D_{n,s,L}$ with high probability and that $Q_{[n]}(Z_{1:n})$ concentrates under the null and alternative; see \Cref{app:typical-set}. Together with our transformation in \Cref{lemma:transform-choice}, these concentration properties allow us to adapt the Lipschitz extension approach of \citet{canonne2020private} and \citet{narayanan2022private} to the regression models considered here.

For $n \geq 2$ and $\epsilon_0>0$, let 
\begin{equation}
 G(z_{1:n}) := \left[Q_{[n]}(z_{1:n}) \right]_0^{B_n}, \quad B_n:=s^2n(n-1), \quad D := \frac{\epsilon_0 B_n}{8 \log(2/\eta)}.
 \label{eq:Lipschitz_test_statistic}
\end{equation}
For the choices of $s$ and $L$ specified below, \Cref{lemma:restricted-sensitivity} implies that $G$ is $D$-Lipschitz with respect to Hamming distance on $\mathcal{D}_{n, s, L}$. Let $\widehat{G}$ be the McShane extension \citep{mcshane1934extension}
\begin{equation}
        \widehat{G}(z_{1:n}) := \inf_{\tilde{z}_{1:n} \in \mathcal{D}_{n, s, L}} \left\{ G(\tilde{z}_{1:n}) + D \Ham(z_{1:n}, \tilde{z}_{1:n})\right\}. 
        \label{eq:mcshane-extension}
\end{equation}
By construction, $\widehat{G}$ is globally $D$-Lipschitz and agrees with $G$ on $\mathcal{D}_{n, s, L}$. Finally, define the test 
\begin{equation}
     \phi_s(Z_{1:n}) := \mathbbm{1} \left\{\widehat{G}(Z_{1:n}) + W > \frac{B_n}{2} \right\}, \quad W \sim \mathrm{Lap}\left(\frac{D}{\epsilon_0} \right).
    \label{eq:lipschitz_test}
\end{equation}

The next lemma gives the privacy and testing guarantees of $\phi_s$ at a fixed scale $s$.

\begin{lemma}[Single-scale mean test]
Let $\{(X_i, Y_i)\}_{i=1}^n$ be i.i.d.~ from $P\in\mathcal P_{\mathrm{Lin}}$ or $P \in \mathcal P_{\mathrm{GLM}}(\mathcal F)$, where $\mathcal{F}$ is a fixed GLM response family satisfying \Cref{assp:glm-regularity}. Set $Z_i := X_i \tanh(Y_i)$ for $i \in [n]$. Let $0<s\leq \sqrt{\overline{\lambda}}$, $0<\eta\leq1/4$, and choose the typical set parameter $L$ according to \eqref{eq:typical-L} with failure probability $\eta/2$ to define the extension $\widehat{G}$. Then the test $\phi_s$ from \eqref{eq:lipschitz_test} is $(\epsilon_0, 0)$-DP. Furthermore, there exists a constant $C$ such that whenever 
\begin{align}
    s \geq \rho(n,d,\epsilon_0, \eta)
    &:= C\left\{\frac{l_{n,\eta}d^{1/4}}{\sqrt n}+\frac{l_{n,\eta}d^{1/4}}
        {n^{3/4}\epsilon_0^{1/2}} +\frac{l_{n,\eta}^2}{n\epsilon_0}\right\}, 
        \quad l_{n,\eta}:= 1+\log \left(\frac{24(n+1)}{\eta}\right),
    \label{eq:IT-test_s-radius}
\end{align} 
the type-I and type-II errors of $\phi_s$ in testing for
\begin{equation*}
    H_0:\mathbb{E}[Z]=0
    \qquad\text{against}\qquad
    H_1(s):s\leq\|\mathbb{E}[Z]\|_2\leq2s
\end{equation*}
are each at most $\eta$. 
\label{lemma:single_scale_test}
\end{lemma}

\Cref{lemma:transformation_properties} transfers the mean testing guarantee to regression signal detection. With $c_\mu>0$ as the constant in \eqref{eq:mean-to-regression}, every alternative with $\|\beta\|_2\geq\tilde r/c_\mu$ satisfies $\|\E[Z]\|_2\geq\tilde r$, provided $\tilde r\leq c_\mu$. The remaining step is to remove the dependence on the unknown scale of $\|\E[Z]\|_2$. 

To cover the full alternative, we apply \Cref{lemma:single_scale_test} over a dyadic grid of signal strengths. Splitting the privacy budget across scales and controlling the errors by a union bound only inflates the detection boundary by logarithmic factors; see \Cref{app:multiscale_test} for details. Combining this multiscale construction with \Cref{lemma:transformation_properties} gives the dense directional component in the minimax rate. 

Computing the extension involves optimisation over typical datasets $\mathcal{D}_{n, s, L}$, which is not directly polynomial-time implementable. In \Cref{sec:polytime}, we construct a polynomial-time test attaining the same directional rate, up to logarithmic factors, when $\epsilon\gtrsim n^{-1/2}$.

\subsection{Sparse upper bound}
\label{sec:sparse-cov}
The test in the previous section uses the transformed mean $\E[X \tanh(Y)]$, which may not be sparse even when $\beta$ is. To use the sparsity of $\beta$ directly, we consider sparse covariate projections and apply $\tanh$ to obtain statistics with uniformly bounded sensitivity. We test whether $\cov(\tanh(v^\top X),\tanh(Y))=0$. Under the alternative, we have that 
\begin{equation*}
    \left|\cov \left(\tanh(v^\top X), \tanh(Y) \right) \right| \gtrsim \min\{\|\beta\|_2, 1\}, \quad v:= \beta / \|\beta\|_2.
\end{equation*}
This signal persists for sufficiently close directions, allowing us to restrict attention to a fixed set of candidate directions. 

To adapt to unknown sparsity, we combine tests over the dyadic grid $\mathcal K:=\{2^j:0\leq j\leq\lfloor\log_2d\rfloor\}\cup\{d\}$. For a covering radius $\eta_{\mathrm{net}}\in(0,1)$, and $S\subseteq[d]$ with $|S| \in \mathcal{K}$, fix a Euclidean $\eta_{\mathrm{net}}$-net $\mathcal V_S$ of $\{v\in\mathbb R^d:\|v\|_2=1,\ \operatorname{supp}(v)\subseteq S\}$ such that $|\mathcal V_S|\leq(1+2/\eta_{\mathrm{net}})^{|S|}$; such nets exist by Corollary 4.2.11 of \cite{vershynin2026}. For each $k \in \mathcal{K}$, collect the candidate directions in $\mathcal V_k:=\bigcup_{S\subseteq[d], |S|=k}\mathcal V_S$, and define, for $n \geq 2$, 
\begin{equation*}
 T_k:=\max_{v\in\mathcal V_k}\frac{1}{2m}\sum_{j=1}^m \{\tanh(v^\top X_{2j-1})-\tanh(v^\top X_{2j})\} \{\tanh(Y_{2j-1})-\tanh(Y_{2j})\}, \quad m := \lfloor n/2 \rfloor.
\end{equation*}
For each fixed $v$, the paired differences give an unbiased covariance estimator with sensitivity at most $4/m$. Maximising over directions and sparsity levels preserves this bound, so a single Laplace perturbation ensures privacy without splitting the budget over sparsity levels:
\begin{equation}
 \phi_{\mathrm{sparse}}:= \ind\left\{\max_{k\in\mathcal K} \left\{T_k-\sqrt{\frac{8\log(2|\mathcal{K}||\mathcal{V}_k|/\eta)}{m}}\right\}
 +\frac{4W}{m\epsilon} >\frac{4\log(1/\eta)}{m\epsilon} \right\},
 \quad W\sim\operatorname{Lap}(1).
 \label{eq:adaptive-covariance-test}
\end{equation}

\begin{proposition}
\label{prop:sparse-upper}
Let $\{(X_i, Y_i)\}_{i=1}^n$ be i.i.d.~from $P\in\mathcal P_{\mathrm{Lin}}$ or $P \in \mathcal P_{\mathrm{GLM}}(\mathcal F)$, where $\mathcal{F}$ is a fixed GLM response family satisfying \Cref{assp:glm-regularity}, and fix $\eta \in (0, 1/4)$. The test in \eqref{eq:adaptive-covariance-test} is $(\epsilon,0)$-DP and has type-I error at most $\eta$. There exist constants $\eta_{\mathrm{net}}\in(0,1)$ and $c>0$ such that the type-II error of $\phi_{\mathrm{sparse}}$ is at most $\eta$ provided that $\|\beta\|_0\leq s^*$ and
\begin{equation}
 \sqrt{\frac{s^*\log(ed/s^*)+\log (|\mathcal{K}|/\eta)}{n}} +\frac{\log(1/\eta)}{n\epsilon} \leq c\min\{\|\beta\|_2,1\}.
 \label{eq:adaptive-scan-power}
\end{equation}
\end{proposition}

The test $\phi_{\mathrm{sparse}}$ is therefore optimal in the sparse regime, up to logarithmic factors. The number of candidate test directions affects the non-private rate but not the scale of the privacy noise. The proof is given in \Cref{app:sparse}. Checking all candidate directions to adapt over sparsity levels up to $s^*=d$ requires $O(nd(2+2/\eta_{\mathrm{net}})^d)$ operations. Restricting attention to levels $k\leq s_{\max}$, where $s_{\max}$ is a fixed constant, reduces the computational cost to $O(nd^{s_{\max}})$, with adaptation over these levels.

\subsection{Response marginal test}
\label{sec:radial_test}
For GLMs with more than two response support points, the lower bounds identify second-order information in the response marginal that cannot be absorbed by adjusting the nuisance intercept. This residual information, quantified by $\E_a[R_a^2]$ in \eqref{eq:projection-resdiual-def}, motivates a test based on the responses alone.

We construct a bounded  statistic that preserves a signal of order $\|\beta\|_2^2$ for sufficiently small $\|\beta\|_2$. Specifically, we truncate $R_a(Y)$ from \eqref{eq:projection-resdiual-def}, centre the resulting function under the null, and then re-orthogonalise it to the nuisance score. The resulting functions have the following properties. 

\begin{lemma}
\label{lemma:bounded-orthog-contrast}
    Let $\mathcal{F}$ be any response family that satisfies \Cref{assp:glm-regularity} and has more than two response support points. Let $\Gamma_1$ be a compact interval satisfying $\Gamma \subset \intr(\Gamma_1) \subset\Gamma_1 \subset U$. There exist measurable functions $\{g_a\}_{a \in \Gamma_1}$ and constants $c_g, B_g >0$ such that, uniformly over $a \in \Gamma_1$, we have 
    \begin{align}
    \E_a[g_a(Y)]&=0,\label{eq:g-centred}\\
    \E_a[g_a(Y)S_a(Y)]&=0,\label{eq:g-orthogonal}\\
    \E_a[g_a(Y)H_a(Y)]&\ge c_g,\label{eq:g-curvature}\\
    \sup_y|g_a(y)|&\le B_g.\label{eq:g-bounded}
\end{align}
\end{lemma}
Property \eqref{eq:g-centred} ensures that the statistic has mean zero under the null hypothesis, when evaluated at the true intercept, while \eqref{eq:g-orthogonal} gives nuisance orthogonality. In particular, since
\begin{equation*}
    \left.\frac{\mathrm d}{\mathrm dt}\E_{a+t}[g_a(Y)]\right|_{t=0} = \E_a[g_a(Y) S_a(Y)] = 0,
\end{equation*}
the statistic is locally insensitive to errors in estimating the nuisance parameter, in the spirit of Neyman's $C(\alpha)$ principle \citep[e.g.][]{gu2016neyman}.  The condition in \eqref{eq:g-curvature} ensures a positive second-order signal, while the uniform bound in \eqref{eq:g-bounded} controls the global sensitivity of the empirical mean of $g_a(Y)$. The construction of $g_a$, the corresponding constants $c_g$ and $B_g$, and the proof of \Cref{lemma:bounded-orthog-contrast} are given in \Cref{app:nuisance-orthog-test}. 

The testing procedure consists of two stages, which are described in \Cref{alg:generic-radial-test}. In the first stage, we use the first half of the samples to obtain a private estimate of the intercept $\widehat{\gamma}$. In the second stage, we use the remaining samples to privately estimate $\E[g_{\widehat{\gamma}}(Y)]$ and reject the null hypothesis if the estimate exceeds a threshold.

\begin{algorithm}[h!]
\caption{Response-only test}
\label{alg:generic-radial-test}
\begin{algorithmic}[1]
\Require Data $\{Y_i\}_{i=1}^n$ with $n \geq 2$, privacy parameter $\epsilon>0$, error probability $\eta\in(0,1/4)$, rejection threshold constant $c_T>0$. 
\State Let $n_1=\lfloor n/2\rfloor$, $n_2=n-n_1$. 

\State  Set
\begin{equation}
    \widetilde{m} := \frac1{n_1}\sum_{i=1}^{n_1} \tanh(Y_i)
    +W_1, \qquad W_1 \sim \mathrm{Lap} \left(\frac{2}{n_1 \epsilon}\right)
    \label{eq:tildem-def}
\end{equation}
\State Obtain the private intercept estimate 
\begin{equation}
    \widehat\gamma := \argmin_{a \in \Gamma_1}|m(a) - \widetilde{m}|, \qquad m(a):=\E_a[\tanh(Y)].
    \label{eq:private-intercept-def}
\end{equation}

\State Set the private test statistic 
\begin{equation}
    \widetilde{T} :=    \frac{1}{n_2}\sum_{i=n_1+1}^{n}
        g_{\widehat\gamma}(Y_i)
    +W_2, \qquad \qquad W_2 \sim \mathrm{Lap} \left(\frac{2B_g}{n_2 \epsilon}\right), \,\ W_2 \perp W_1
\end{equation}

\State Set $a_n(\epsilon, \eta)
:=\sqrt{\log(8/\eta)/n}+\log(8/\eta)/(n\epsilon)$.
\State \Return $\ind\{\widetilde{T}\geq c_T\,a_n(\epsilon, \eta)\}$.
\end{algorithmic}
\end{algorithm}

\Cref{alg:generic-radial-test} satisfies the following privacy and error guarantees.
\begin{proposition}
    \label{prop:generic-radial-test}
    The test in \Cref{alg:generic-radial-test} is $(\epsilon, 0)$-DP. Suppose further that the observations are drawn from some $P \in \mathcal{P}_{\mathrm{GLM}}(\mathcal{F})$ where $\mathcal{F}$ satisfies \Cref{assp:glm-regularity} and has more than two support points. Fix $\eta \in (0, 1/4)$. For any sufficiently large fixed $c_T>0$, there exist constants $C,r_0,a_*>0$ such that when \Cref{alg:generic-radial-test} is run with threshold constant $c_T$, it has type-I error at most $\eta$ whenever $a_n(\epsilon, \eta) \leq a_*$ and type-II error at most $\eta$ for every alternative satisfying
\begin{equation}
    0<r\leq\|\beta\|_2\le r_0,\qquad r^2 \geq
  C\left\{\sqrt{\frac{\log(8/\eta)}{n}}+\frac{\log(8/\eta)}{n\epsilon} \right\}.
  \label{eq:radial-separation}
\end{equation}
\end{proposition}

\Cref{prop:generic-radial-test} shows that, for GLMs with more than two response support points, the regression signal can be detected at rate
\begin{equation*}
r^* \lesssim_{\log} n^{-1/4}+(n\epsilon)^{-1/2}
\end{equation*}
using the responses alone. This gives the marginal contribution to the GLM minimax rate of \Cref{thm:main-result}. Unlike the mean-reduction test, which detects dependence between the covariates and response, the response-only test detects a second-order change in the response marginal without searching for the unknown signal direction. Its rate is therefore independent of $d$ and can improve on the mean-reduction test's rate when $d\gtrsim n$. The two rates are comparable when $d\asymp n$, although the sparse covariance test may give a smaller rate than either.

\begin{remark}
\label{remark:exponential-kernel-test}
The restriction $\|\beta\|_2\leq r_0$ in \Cref{prop:generic-radial-test} comes from the local likelihood expansion used to bound the type-II error and carries over to the upper bound in \Cref{thm:main-result}. For some response families, invariance to the intercept allows us to avoid estimating it and to detect signals beyond this local neighbourhood. For example, in exponential regression with log rate $\gamma+x^\top\beta$, the intercept scales independent responses $Y_1$ and $Y_2$ by the same factor. The kernel
\begin{equation}
h_{\mathrm{Exp}}(y_1,y_2):=\ind\{\max(y_1,y_2)>2\min(y_1,y_2)\}-2/3
\label{eq:exponential-test-kernel}
\end{equation}
is invariant under this common scaling, and is bounded and satisfies $\E_{0,\gamma}[h_{\mathrm{Exp}}(Y_1,Y_2)]=0$ and $\E_{\beta,\gamma}[h_{\mathrm{Exp}}(Y_1,Y_2)]\gtrsim\min\{\|\beta\|_2^2,1\}$ for all $\gamma\in\mathbb R$. Privatising the corresponding U-statistic therefore gives the same separation rate as the generic response-only test, without requiring $\|\beta\|_2\leq r_0$. Further details are given in \Cref{app:exponential_regression}.
\end{remark}

To complete the upper bound in \Cref{thm:main-result}, we combine the multiscale mean-reduction test (\Cref{subsec:generalised_mean_test}) and the sparse covariance test (\Cref{sec:sparse-cov}), adding the response-only test for GLMs with more than two response support points. We omit any component whose separation condition cannot be met within the alternative class, allocate a fixed fraction of the privacy budget to each remaining test, and reject if any component rejects. The combined test is $(\epsilon,0)$-DP by composition (\Cref{lemma:composition}). Its type-I error is controlled by a union bound and its type-II error is controlled whenever at least one component's separation condition holds. The combined test therefore attains the minimum of the component rates, up to polylogarithmic factors, without requiring knowledge of the sparsity level. Full details are given in \Cref{app:upper_assembly}.

\section{Computation and numerical experiments}
\label{sec:polytime}
The mean-reduction test in \Cref{subsec:generalised_mean_test} and the sparse covariance test in \Cref{sec:sparse-cov} attain the dense directional and sparse rates in \Cref{thm:main-result}, respectively, up to logarithmic factors. However, the former relies on a Lipschitz extension with no polynomial-time guarantee, while the latter checks exponentially many directions. We now give computationally efficient alternatives and quantify their additional privacy costs, and then present numerical experiments in \Cref{subsec:numerical}.

\subsection{A polynomial-time mean-reduction test}
We consider the same quadratic statistic as in \Cref{subsec:generalised_mean_test}, but use a clipped version of the inner product and a different privatisation mechanism. For some fixed constant $C_\mathrm{clip}>0$ and a failure level $\alpha \in (0, 1)$, define the kernel for scale $s$ as 
\begin{equation}
    h_{s}(z_1, z_2) := [z_1^\top z_2]_{-R(s)}^{R(s)}, \quad R(s) := C_\mathrm{clip}(\sqrt{d l_s} + l_s), \quad l_s:= 1 +  \log\{(d+1) (n+1)(s^{-1}+1)/\alpha\}.
    \label{eq:clipped_kernel_def}
\end{equation}
We privatise the resulting U-statistic using the local H\'ajek mechanism of \cite{chaudhuri2024DPUstats}, which downweights influential observations and calibrates noise to a smooth sensitivity bound. In \Cref{lemma:clipped_kernel}, we bound the clipping bias, variance and deviations of the kernel averages associated with individual observations to establish the type-I and type-II error of the test. As in \Cref{subsec:generalised_mean_test}, we combine tests over a dyadic grid of signal strengths to cover the full alternative. The full testing procedure is specified in \Cref{alg:polytime_test}. 

\begin{algorithm}[h!]
\caption{Polynomial-time mean-reduction test}
\label{alg:polytime_test}
\begin{algorithmic}[1]
\Require Data $\{Z_i\}_{i=1}^n$ with $n \geq 2$, smallest separation scale $0<\tilde{r}\leq \sqrt{\overline{\lambda}}$, privacy budget $\epsilon>0$, failure probability $\eta\in (0, 1/4)$. 

\State Set $J \gets \max\{1, \lceil \log_2(\sqrt{\overline{\lambda}}/\tilde{r}) \rceil \}$ and $\alpha \gets\eta/ (3J)$.
\For{$j=0,\ldots,J-1$}
\State Set $s_j \gets 2^j \tilde{r}$ and $\xi_{s_j} \gets C_\mathrm{tol}\, l_{s_j}\left(s_j+\sqrt{d / n}\right)$
\State Set $\widetilde U_{j}\gets\textsc{LocalH\'ajekMean}(Z_{1:n},h_{s_j},2R(s_j),\xi_{s_j},\epsilon/J)$.  
\If{ $\widetilde{U}_{j} > s_j^2/2$}
\State \Return $1$. 
\EndIf
\EndFor
\State \Return $0$.
\end{algorithmic}
\end{algorithm}

The \textsc{LocalH\'ajekMean} algorithm in \Cref{alg:polytime_test} can be found as \Cref{alg:local-hajek-primitive} in \Cref{app:polytime_implementation}; it specialises Algorithm 1 of \cite{chaudhuri2024DPUstats} to U-statistics of order two. The test has the following guarantees. 

\begin{proposition} 
    \label{prop:polytime_test}
The test in \Cref{alg:polytime_test} is $(\epsilon, 0)$-DP and requires  $O(Jdn^2)$ arithmetic operations. Let $\{(X_i, Y_i)\}_{i=1}^n$ be i.i.d.~from $P\in\mathcal P_{\mathrm{Lin}}$ or $P \in \mathcal P_{\mathrm{GLM}}(\mathcal F)$, where $\mathcal{F}$ is a fixed GLM response family satisfying \Cref{assp:glm-regularity}, and fix $\eta \in (0, 1/4)$. For sufficiently large fixed choices of the kernel clipping constant $C_{\mathrm{clip}}$ and H\'ajek tolerance constant $C_{\mathrm{tol}}$, there exist constants $c, C, C_1>0$ such that whenever 
\begin{equation}
    c \eta ^{-1/2} \ell_{n, d, \epsilon}^3 \left\{ \frac{d^{1/4}}{n^{1/2}} + \frac{d^{1/4}}{n^{3/4} \epsilon^{1/2}}  + \frac{d^{1/4} }{n\epsilon} + \frac{d^{1/4}}{n^{3/2} \epsilon^{3/2}}\right\} \leq \tilde{r} \leq C, 
    \label{eq:polytime-rate}
\end{equation}
the sum of the type-I and type-II errors of \Cref{alg:polytime_test} for the hypothesis test in \eqref{eq:hypotheses-def} with $r=C_1\tilde{r}$ and $R=\infty$ is bounded by $2\eta$. With this $\tilde{r}$, the running time of the test is $O(dn^2\log(n+1))$.
\end{proposition}
\Cref{prop:polytime_test} shows that \Cref{alg:polytime_test} provides a polynomial-time alternative to the mean-reduction test in \Cref{subsec:generalised_mean_test}, at the cost of additional privacy-dependent terms. The final term in \eqref{eq:polytime-rate} is dominated by the third term whenever $n\epsilon\gtrsim 1$. In the complementary regime, there exists a constant $c_0>0$ such that when $n \epsilon < c_0$, group privacy arguments already rule out nontrivial private signal detection. The privacy-dependent term of $1/(n\epsilon)$ of the mean-reduction test's detection rate becomes $d^{1/4}/(n\epsilon)$ in \eqref{eq:polytime-rate}. The two rates match when we have $\epsilon \gtrsim n^{-1/2}$.

\subsection{A coordinate test for sparse alternatives}
\label{sec:sparse-polytime}
We next consider a computationally efficient alternative to the sparse test in \Cref{sec:sparse-cov}, based on the largest coordinate of the transformed sample mean from \Cref{lemma:transformation_properties}. Although $\E[X \tanh(Y)]$ is not necessarily sparse, we show in \Cref{app:sparse-polytime} that its largest coordinate still retains a signal 
\begin{equation*}
    \|\E[X\tanh(Y)]\|_\infty \geq \frac{c_\mu}{\sqrt{s^*}} \min\{\|\beta\|_2, 1\}. 
\end{equation*}

We therefore target this quantity by clipping each covariate coordinate at $B:=C_X \sqrt{\log(2n)}$ and privatising the maximum absolute transformed coordinate mean.

\begin{proposition} \label{prop:sparse-polytime}
Suppose that $\{(X_i,Y_i)\}_{i=1}^n$ are i.i.d.~from $P\in\mathcal P_{\mathrm{Lin}}$ or $P\in\mathcal P_{\mathrm{GLM}}(\mathcal F)$, where $\mathcal F$ satisfies \Cref{assp:glm-regularity}. Fix $\epsilon>0$ and $\eta\in(0,1/4)$. Then the test
\begin{equation*}
\phi := \ind\left\{T+\frac{2B}{n\epsilon}W>b_B+q_B+\frac{2B\log(1/\eta)}{n\epsilon}\right\}, \qquad T:=\max_{j\in[d]}\left|\frac1n\sum_{i=1}^n[X_{ij}]_{-B}^{B}\tanh(Y_i)\right|,
\end{equation*}
with $W\sim\operatorname{Lap}(1)$ satisfies $(\epsilon,0)$-DP and requires $O(nd)$  operations. If $b_B:=C_X^2\exp(-B^2/C_X^2)/B$ and  $q_B:=\sqrt{2\overline\lambda\log(4d/\eta)/n}+ 4B\log(4d/\eta)/(3n)$, the test has type-I error at most $\eta$. Furthermore, there exists a constant $C>0$ such that, for every $1\leq s^*\leq d$, the test's type-II error is at most $\eta$ provided that $\|\beta\|_0\leq s^*$ and
\begin{equation}
C\sqrt{s^*}\left\{\sqrt{\frac{\log(4d/\eta)}n}+\frac{\sqrt{\log(2n)}\log(4d/\eta)}n+\frac{\sqrt{\log(2n)}\log(1/\eta)}{n\epsilon}\right\}\leq\min\{\|\beta\|_2,1\}.
\label{eq:coordinate-power}
\end{equation}
\end{proposition}

The test in \Cref{prop:sparse-polytime} does not require knowledge of $s^*$ and matches the sparse test's separation rate, up to polylogarithmic factors, when $\epsilon\gtrsim n^{-1/2}$. Under stronger privacy constraints, its upper bound incurs an additional factor of $\sqrt{s^*}$ in the privacy term. 

The response-only test in \Cref{alg:generic-radial-test} also has an efficient Monte Carlo implementation for many commonly used response families. As long as we can efficiently sample from $P_a$ and evaluate $S_a(y)$ and $R_a(y)$ (from \Cref{lemma:projection}) for every $a\in\Gamma_1 $, the required one-dimensional expectations can be approximated in $O(n^2/\eta)$ time while preserving the separation rate; see \Cref{app:mc-approx} for details. Combining Propositions \ref{prop:generic-radial-test}, \ref{prop:polytime_test}, and \ref{prop:sparse-polytime} then gives a polynomial-time procedure across all regimes. When $\epsilon\gtrsim n^{-1/2}$, this procedure attains the minimax rates in \Cref{thm:main-result}, up to polylogarithmic factors. This is the same privacy regime in which the efficient private test of a dense Gaussian mean in \cite{narayanan2022private} attains the optimal rate up to polylogarithmic factors. We leave the question of whether the additional privacy costs in \eqref{eq:polytime-rate} and \eqref{eq:coordinate-power} reflect fundamental computational-statistical gaps or are an artefact of our test constructions to future work.

\subsection{Numerical experiments}
\label{subsec:numerical}
In this subsection, we evaluate the numerical performance of our computationally efficient tests.  The code for these experiments can be found at \url{https://github.com/EKHung/DP_regression_testing}. For each experimental configuration, we use 250 repetitions under each of the null and alternative hypotheses and report the estimated sum of the type-I and type-II errors, with the error bars showing Monte Carlo standard errors.

We set the tuning parameters of \Cref{alg:polytime_test} to be $C_{\mathrm{clip}}=0.5$, $C_{\mathrm{tol}}=0.05$ and $\eta=0.05$. For the smallest scale $\tilde r$, we use the sum of the four polynomial terms in \eqref{eq:polytime-rate}, multiplied by $C_{\mathrm{sep}}=5$. For the test in \Cref{prop:sparse-polytime}, the rejection threshold depends on the covariate population quantities $C_X$ and $\overline{\lambda}$;  the sensitivity analyses in \Cref{app:numerical} suggest reasonable robustness to their misspecification. We also examine how the performance of \Cref{alg:polytime_test} depends on its tuning choices.

We first compare our tests with the following private procedures:
\begin{itemize}
    \item \textbf{dpHSIC}: a DP analogue of the permutation HSIC independence test, from \cite{kim2026permutation}.
    \item \textbf{ToT-F}: the test-of-tests framework of \cite{pmlr-v202-kazan23a}, based on the subsample-and-aggregate approach of \cite{nissim2007smooth}. We use the regression F-test as the non-private test within each block.
    \item \textbf{PB-LRT}: the private likelihood ratio test of \citet{pena2024uncertainty} for Gaussian linear regression. This test also splits the sample into blocks, but then computes and clips the blockwise log-likelihood ratios, averages them, and adds Laplace noise.
\end{itemize}
\Cref{alg:polytime_test} and all the tests listed above satisfy pure differential privacy, i.e.~$\delta=0$ in \Cref{def-cdp}. We provide further implementation details, tuning choices and their sensitivity analyses for the competing tests in \Cref{app:numerical}.

For the comparison, we vary $n \in \{2000, 2300, \dots, 5000\}$, fix $\epsilon=1$ and consider the linear regression model \eqref{eq:linear_regression_model}, with independent covariates $X_i\sim\mathcal{N}(0,I_{500})$ and errors $\xi_i\sim\mathcal{N}(0,1)$. This ensures that the Gaussian regression models underlying the F-test and likelihood ratio competitors are correctly specified. Our experiments for the alternative use $\beta = (1, 0, 0, \dots, 0)$. 
\begin{figure}[h!]
    \centering
    \includegraphics[width=0.7\linewidth]{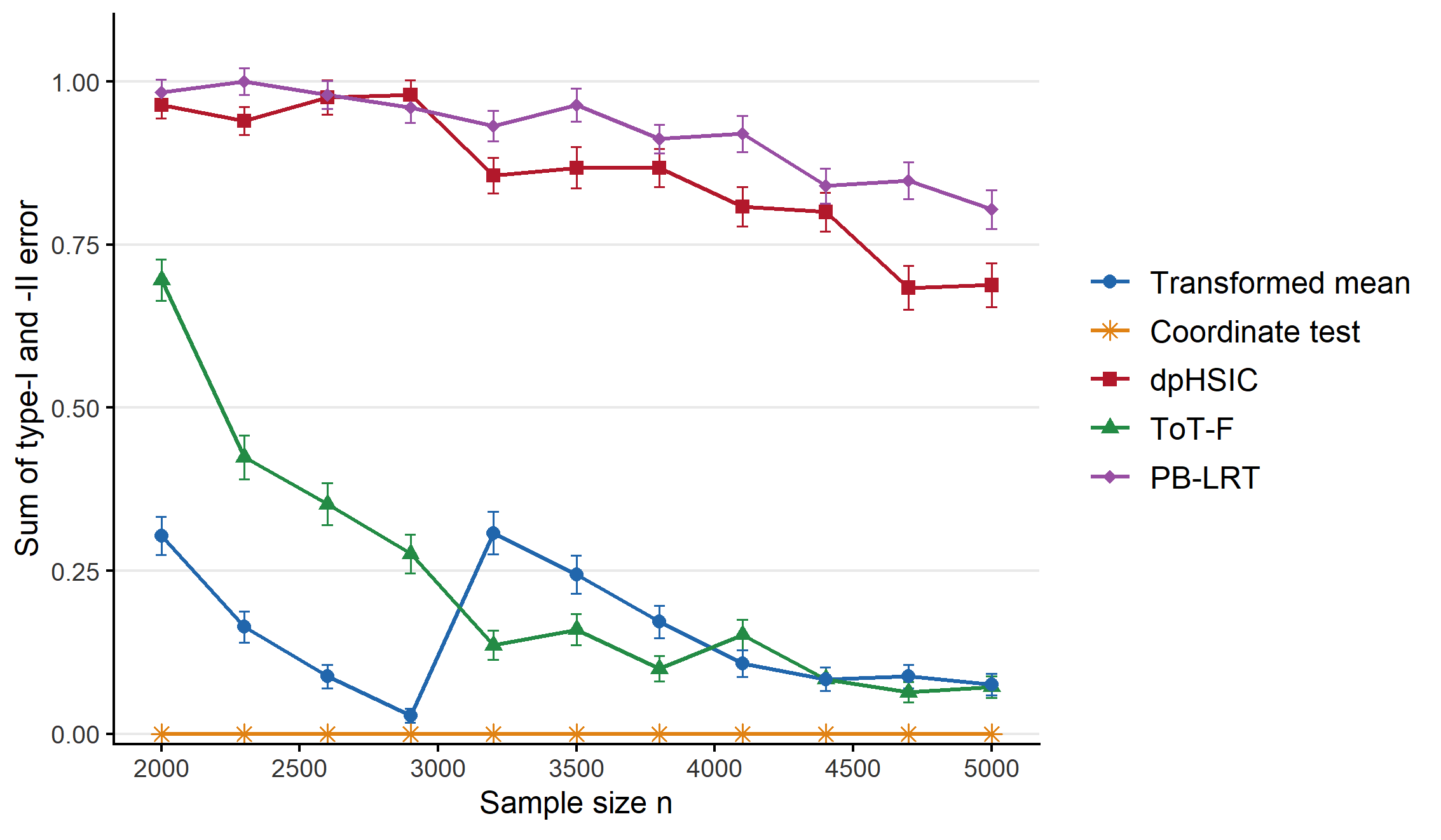}
    \caption{Comparison of private tests under Gaussian linear regression with $d=500$ and $\epsilon=1$.}
    \label{fig:Gaussian_comparison}
\end{figure}

\Cref{fig:Gaussian_comparison} shows that the coordinate test from \Cref{prop:sparse-polytime} makes no error across the simulated repetitions. Since the alternative satisfies $\|\beta\|_0=1$, it is favourable for the sparse test, but we also see that \Cref{alg:polytime_test} performs competitively, especially for smaller samples. The increase in total error of \Cref{alg:polytime_test} between $n=2900$ and $n=3200$ coincides with the addition of a scale, which reduces the privacy budget allocated to each scale. The error subsequently decreases with additional samples. The smaller jumps in the ToT-F curve coincide with changes in the number of aggregation blocks.

\begin{figure}[h!]
    \centering
    \includegraphics[width=0.95\linewidth]{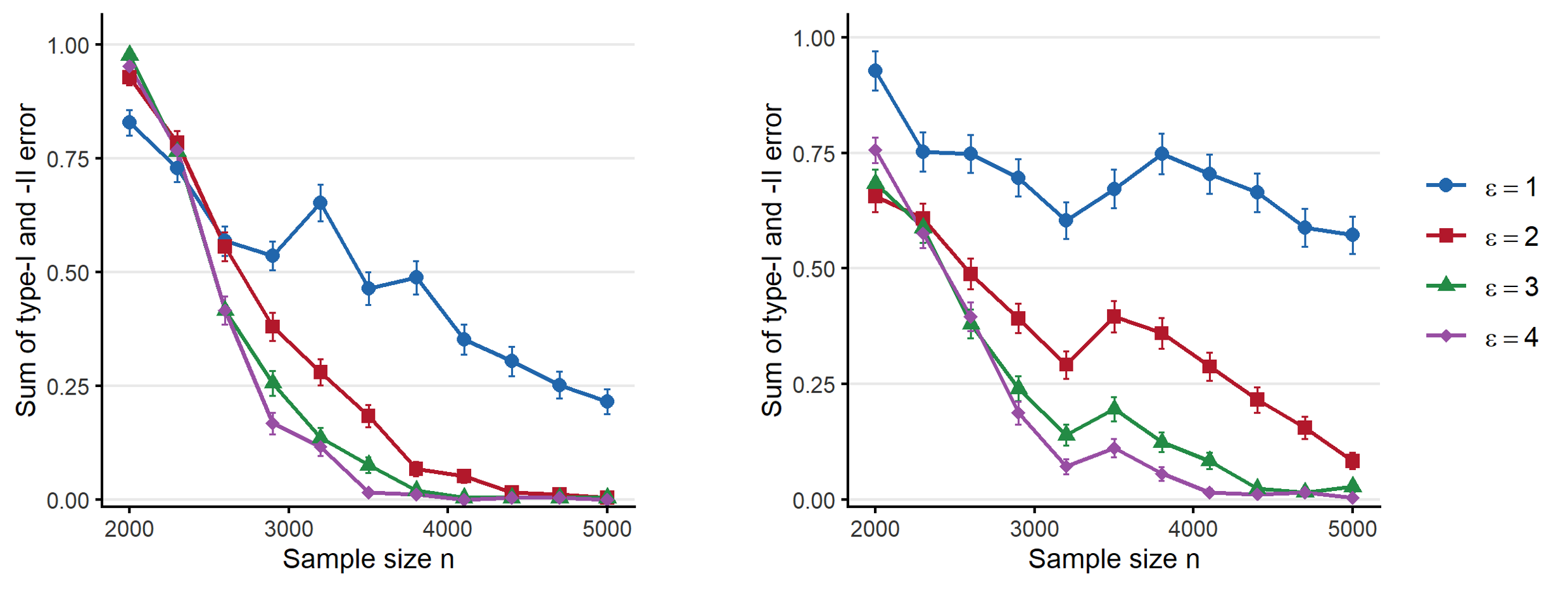}
    \caption{Performance of the polynomial-time mean-reduction test under Cauchy linear regression with $d=500$ (left) and exponential regression with $d=50$ (right), for privacy parameters $\epsilon \in \{1, 2, 3, 4\}$.}
        \label{fig:Cauchy_exponential}
\end{figure}

We also consider two models with independent Rademacher covariate coordinates, and set $\beta = (1/\sqrt{d}, 1/\sqrt{d}, \cdots, 1/\sqrt{d})$, using the same grid of sample sizes as in the Gaussian comparison. The left panel of \Cref{fig:Cauchy_exponential} shows results for a linear regression model \eqref{eq:linear_regression_model} with $d=500$, and errors distributed as $\mathrm{Cauchy}(0, 1)$. This setting illustrates the performance of our test when the errors have neither a finite mean nor a finite variance. The right panel of \Cref{fig:Cauchy_exponential} shows results for an exponential regression model, with $d=50$ and a conditional rate $\exp(-(0.2 + X_i^\top \beta))$. In both settings, we vary $\epsilon \in \{1, 2, 3, 4\}$. We see that the estimated testing error generally decreases as the sample size and privacy budget increase. 

\begin{figure}[h!]
    \centering
    \includegraphics[width=0.7\linewidth]{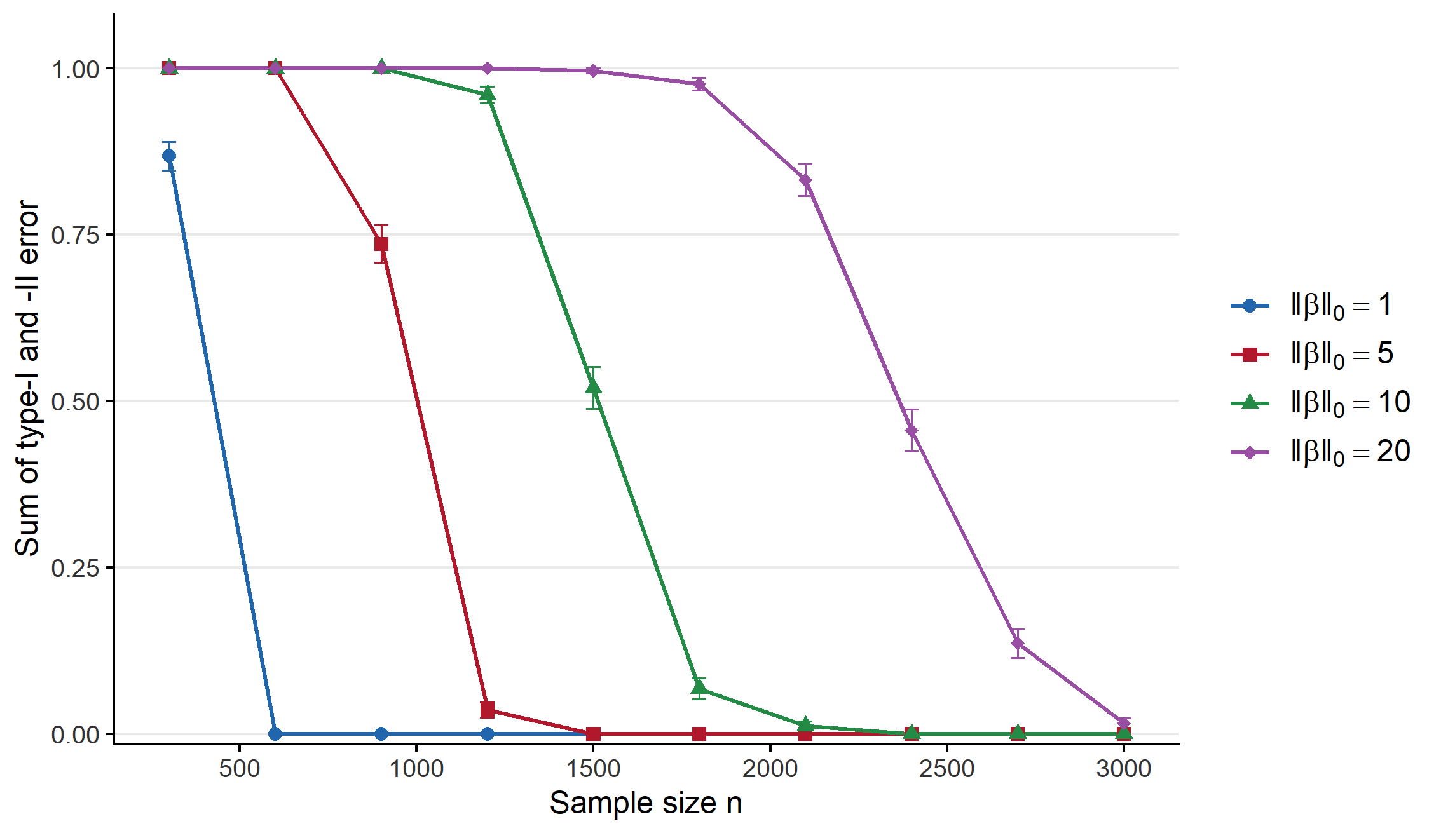}
    \caption{Estimated sum of type-I and type-II errors for the coordinate test under Gaussian linear regression with $d=500$, $\epsilon=1$ and $\|\beta\|_2=1$, for sparsity levels $s^*\in\{1,5,10,20\}$.}
    \label{fig:sparsity}
\end{figure}

Next, we evaluate how our coordinate test from \Cref{sec:sparse-polytime} performs at four sparsity levels $\|\beta\|_0 \in \{1, 5, 10, 20\}$. We fix $\epsilon=1$ and $\|\beta\|_2=1$, and set each non-zero coefficient to $1/\sqrt{\|\beta\|_0}$. To capture the transition from high to low testing error, we consider sample sizes $n \in \{300, 600, \dots, 3000\}$. The results are reported in \Cref{fig:sparsity}. Consistent with \Cref{prop:sparse-polytime}, sparser alternatives require fewer observations to achieve low testing error, illustrating how the test benefits from sparsity without requiring knowledge of the sparsity level.

\begin{figure}[h!]
    \centering    \includegraphics[width=0.7\linewidth]{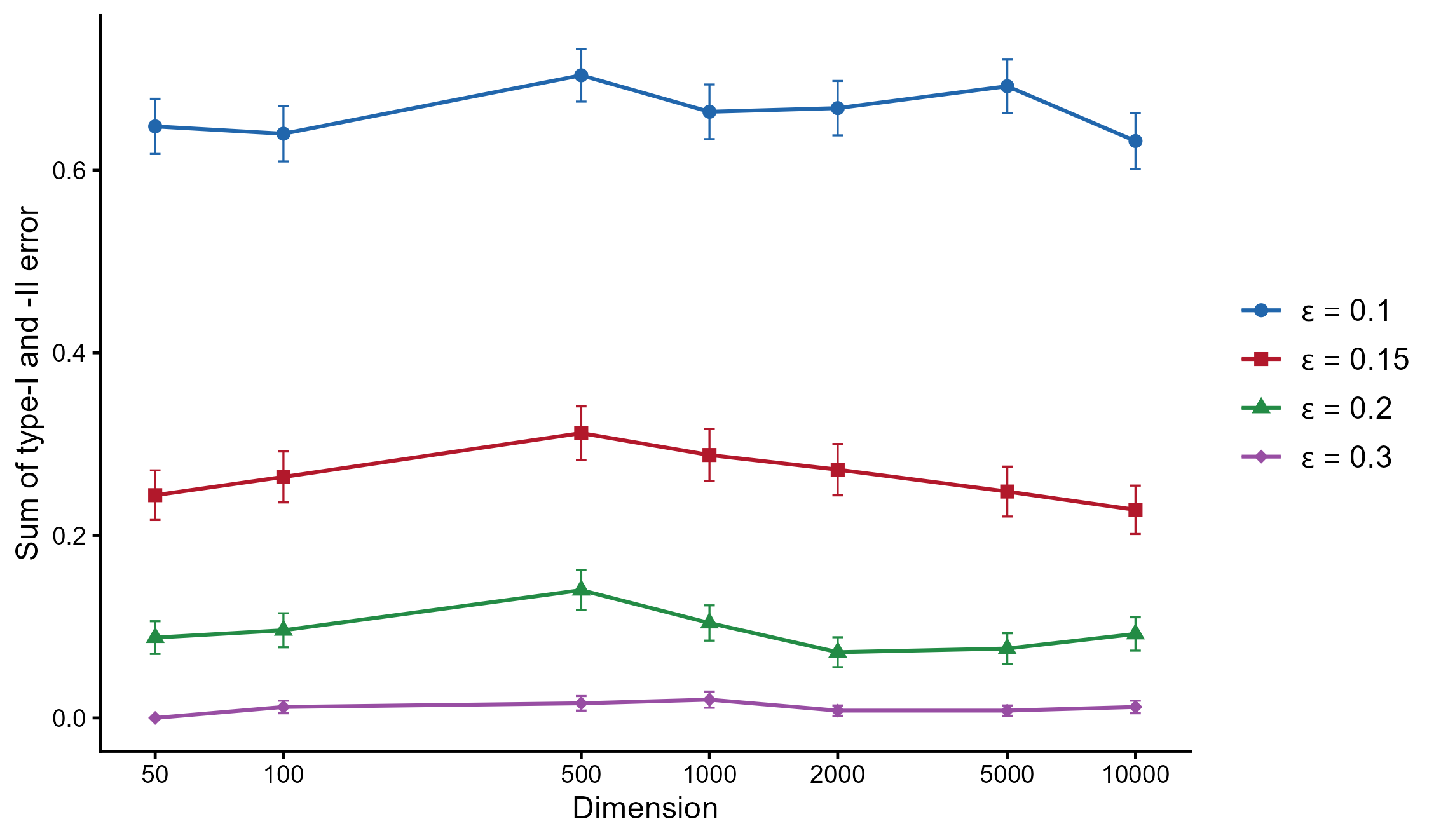}
    \caption{Estimated sum of type-I and type-II errors from numerical experiments using the test in \Cref{remark:exponential-kernel-test}, with $n=2000$ and $\|\beta\|_2=1$, for $\epsilon \in \{0.1, 0.15, 0.2, 0.3\}$. We vary the dimension from 50 to 10,000 along the x-axis on a logarithmic scale.}
\label{fig:exponential_dimension}
\end{figure}

To illustrate a response-only mechanism, we also consider the exponential response-only test described in \Cref{remark:exponential-kernel-test} and \Cref{app:exponential_regression}, setting $\eta=0.05$ in \eqref{eq:exponential-test}. We use the same exponential regression model as above, fix $n=2000$, and vary $d$ from 50 to 10,000 and $\epsilon\in\{0.1,0.15,0.2,0.3\}$. \Cref{fig:exponential_dimension} shows that testing error remains broadly stable as the dimension increases, including beyond the sample size. 

\section{Discussion}
We established minimax separation rates for private signal detection in linear regression and GLMs with unknown nuisance parameters. The rates depend on whether nuisance adjustment can absorb the change in the response marginal: when it cannot, response-only testing can improve detection in high dimensions. We also characterised how sparsity improves the separation rates and constructed tests that attain these rates up to polylogarithmic factors without knowing the sparsity level.

There are several directions for future work. First, it remains open whether polynomial-time tests can attain the minimax separation rates in the very high privacy regime $\epsilon\lesssim n^{-1/2}$. One possible approach for the dense regime, following \cite{hopkins2023robustness}, is to investigate whether sum-of-squares methods can approximate the test in \eqref{eq:lipschitz_test} while preserving its privacy and power guarantees. Second, we could consider testing a subset of regression coefficients, treating the remaining coefficients as potentially high-dimensional nuisance parameters. The nuisance orthogonality used in our response-only test motivates exploring orthogonal test statistics for this setting, drawing on related ideas from double machine learning \citep{chernozhukov2018double}. A central question is how nuisance complexity and privacy constraints jointly affect the minimax separation rates.

\bigskip \noindent \textbf{Acknowledgements.} The author is supported by the Chancellors' Scholarship scheme and the Statistics Centre for Doctoral Training at the University of Warwick. The author thanks Benedict Risebrow for helpful comments. 

\bigskip \noindent \textbf{Use of AI tools.} OpenAI’s ChatGPT was used to assist with mathematical derivations, refining proof arguments, simulation code development, and language polishing. The author reviewed and revised the resulting material and takes full responsibility for the work.

\bibliographystyle{apalike}
\bibliography{bibliography}

\newpage
\appendix
\renewcommand{\contentsname}{Contents of the appendices}
\addtocontents{toc}{\protect\setcounter{tocdepth}{5}}
\tableofcontents
\crefalias{section}{appendix}
\crefalias{subsection}{appendix}

\section{Technical details for \Cref{sec:lower-bounds}}
\label{app:lower_bounds}
\Cref{app:general_lower_bound} provides the proof of our general lower bound. Appendices \ref{app:local-likelhood} and \ref{app:mixture_calculations} develop the constructions for obtaining the signal detection lower bounds. The proof of the lower bounds in \Cref{thm:main-result} is in \Cref{app:lower-assembly}. 

\subsection{Proof of \Cref{prop:general_lower_bound}}
\label{app:general_lower_bound}
\begin{proof}
By the triangle inequality, we have 
\begin{align*}
    \TV(P_0^{\otimes n}, \bar{P}^{(n)}) 
    &\leq \TV(P_0^{\otimes n}, \bar{Q}^{(n)}) + \TV(\bar{Q}^{(n)}, \bar{P}^{(n)})\\
    &\leq \sqrt{\frac{1}{2} n \omega^2} + n\tau . 
\end{align*}
The second inequality is due to conditions (b) and (c) of \Cref{prop:general_lower_bound}, Pinsker's inequality, and convexity and tensorisation of the total variation distance. The non-private rate then follows from Le Cam's two point method. 

To obtain the private rate , we consider two couplings. First, since $\TV(P_0,P_1)\le\alpha$, we couple each pair of observations maximally and independently across indices. The resulting datasets have marginals $P_0^{\otimes n}$ and $P_1^{\otimes n}$, with expected Hamming distance $nD_{\mathrm{TV}}(P_0,P_1)\le n\alpha$.

We next consider a coupling for the mixture distribution to the fixed null. First sample $\beta\sim\pi$. Conditional on $\beta$, for each $i\in[n]$, draw $(D_i, \tilde{D}_i)$ from the maximal coupling of $P_\beta, Q_\beta$. Taking expectation over $\beta$, we have 
    \begin{equation*}
        \mathbb{E}[\Ham (D, \tilde{D})] = \sum_{i=1}^n \Prob (D_i \neq \tilde{D}_i) \leq n \tau, 
    \end{equation*}
    where $D:=(D_1, \dots, D_n)$ and $\tilde{D} := (\tilde{D}_1, \dots, \tilde{D}_n)$. 
    
    Next, by \Cref{lemma:marton_coupling}, there exists a coupling between $\bar{Q}^{(n)}$ and $P_0^{\otimes n}$ such that 
    \begin{equation*}
        \E [\Ham(\tilde{D}, D_0)] \leq \sqrt{\frac{n}{2} \KL \left( \bar{Q}^{(n)} \parallel P_0^{\otimes n} \right)} \leq n \omega / \sqrt{2},
    \end{equation*}
where the final inequality is due to condition (c) of \Cref{prop:general_lower_bound}. 

 We can realise both couplings on a common probability space by first sampling $\tilde{D}$ from its marginal distribution, and then sampling $D$, $D_0$ according to the corresponding conditional coupling distributions. By the triangle inequality, we therefore have a coupling with marginals $P_0^{\otimes n}$ and $\bar{P}^{(n)}$ that satisfies 
 \begin{equation*}
     \mathbb{E}[\Ham (D_0, D)] \leq \mathbb{E}[\Ham (D_0, \tilde{D})] + \mathbb{E}[\Ham (\tilde{D}, D)] \leq n (\tau + \omega). 
 \end{equation*}
Applying \Cref{lemma:DP_le_cam} with the smaller of the Hamming-cost bounds concludes the proof.  
\end{proof}

\begin{lemma}[Hamming cost; \cite{marton1986simple}] Let $\mathcal{X}$ be a Polish space equipped with its Borel $\sigma$-field, and let $\mathcal{P}(\mathcal{X})$ denote the set of all Borel probability measures on $\mathcal{X}$. Let $Q \in \mathcal{P}(\mathcal{X})$ and $P \in \mathcal{P}(\mathcal{X}^n)$. There exists a coupling $(D, D')$ with $D \sim Q^{\otimes n}$ and $D' \sim P$ such that 
\begin{equation*}
    \mathbb{E}[\Ham(D, 
    D')] \leq \sqrt{\frac{n}{2} \KL(P \parallel Q^{\otimes n})}.
\end{equation*}
\label{lemma:marton_coupling}
\end{lemma}
\begin{proof} For completeness, we recall the sequential coupling argument. Construct $(D,D')$ sequentially by maximally coupling, conditional on the previously generated $(D_{<i}, D'_{<i})$, the distributions $Q$ and $P_i(\cdot\mid D'_{<i})$. Then $D\sim Q^{\otimes n}$ and $D'\sim P$, and
\begin{equation*}
\Prob(D_i\neq D_i'\mid D'_{<i})=\TV \left(P_i(\,\cdot\mid D'_{<i}),Q\right).
\end{equation*}
Hence, by Pinsker's inequality and Cauchy--Schwarz,
\begin{align*}
\E[\Ham(D,D')]
&=\sum_{i=1}^n\E\TV \left(P_i(\,\cdot\mid D'_{<i}),Q\right)\\
&\leq\sum_{i=1}^n\sqrt{\frac12\E\KL \left(P_i(\,\cdot\mid D'_{<i}) \parallel Q\right)}\\
&\leq\sqrt{\frac n2\sum_{i=1}^n
\E\KL \left(P_i(\,\cdot\mid D'_{<i}) \parallel Q\right)} = \sqrt{\frac{n}{2} \KL(P \parallel Q^{\otimes n})}. \qedhere
\end{align*}
\end{proof}

\begin{lemma}[Theorem 1 of \cite{acharya2021differentially}] Let $\mathcal{P}_0$ and $\mathcal{P}_1$ denote two classes of distributions on $\mathcal{X}^n$, and let $P_0 \in \mathrm{conv}(\mathcal{P}_0)$ and $P_1 \in \mathrm{conv}(\mathcal{P}_1)$. Let $(D, D')$ be a coupling between $P_0$ and $P_1$ with $H := \mathbb{E}[\Ham(D, D')]$. Then any $(\epsilon, \delta)$-DP testing algorithm $\phi$ must satisfy 
\begin{equation}
    \max_{i \in \{0, 1\}} \sup_{P \in \mathcal{P}_i} \mathbb{P}(\phi(X) \neq i \mid X \sim P)  \geq \frac{1}{2} \max \left\{1 - \TV(P_0, P_1), 0.9 \exp(-10 \epsilon H) - 10H\delta  \right\}. 
\end{equation}
\label{lemma:DP_le_cam}
\end{lemma}

\subsection{Local likelihood expansions and nuisance adjustment}
\label{app:local-likelhood}
\textbf{Notation.} For a linear predictor $a$ in a neighbourhood of $\Gamma$, let $P_a$ denote the response law with natural parameter $\eta_a:=\vartheta(a)$, and write
\begin{equation*}
Z_a:=Y-b'(\eta_a), \qquad\kappa_j(a):=b^{(j)}(\eta_a), \quad j=2,3,4.
\end{equation*}
 We denote the first two likelihood derivatives by
\begin{equation*}
S_a:=\vartheta'(a)Z_a,
\qquad H_a:=\vartheta''(a)Z_a +\vartheta'(a)^2\{Z_a^2-\kappa_2(a)\}.
\end{equation*}
Let $c^*(a)$ denote the least favourable intercept direction, $R_a:=H_a+2c^*(a)S_a$ the corresponding residual after nuisance projection, and $\E_a[R_a^2]$ its squared $L_2(P_a)$ norm. We will also write $\vartheta^{-1}$ for the inverse of $\vartheta$, restricted to $U$. 

We now collect some useful properties that will be used in the lower and upper bound proofs. 

\begin{lemma}[Local likelihood expansion]
\label{lemma:local-likelihood-expansion}
Under \Cref{assp:glm-regularity}, there exist compact intervals
\begin{equation}
\Gamma\subset\intr(\Gamma_1), \qquad \Gamma_1\subset\intr(\Gamma_2), \qquad \Gamma_2\subset U,
    \label{eq:gamma_inclusions}
\end{equation}
a compact interval $K^+\subset\vartheta(U)$, and $\rho_L>0$ such that, for $a\in\Gamma_1$ and $|t_1|,|t_2|\le\rho_L$,
\begin{equation}
    a+t_1,a+t_2\in\Gamma_2,
    \qquad   \vartheta(a+t_1)+\vartheta(a+t_2)-\vartheta(a)\in K^+.
    \label{eq:local-pairwise-parameter}
\end{equation}
Uniformly over $a\in\Gamma_1$, we have the expansion
\begin{equation}
    L_{a,t}(Y) :=\frac{\dint P_{a+t}}{\dint P_a}(Y)
    =1+tS_a(Y)+\frac{t^2}{2}H_a(Y)+ O_{L_2(P_a)}(|t|^3),
    \qquad |t|\le\rho_L.
    \label{eq:likelihood-expansion}
\end{equation}
Moreover the first three derivatives of $t\mapsto L_{a,t}$ exist in $L_2(P_a)$ on this neighbourhood and have uniformly bounded $L_2(P_a)$ norms. In particular, for any fixed $\gamma_0 \in \intr(\Gamma)$ there is a $\kappa>0$ such that $\gamma_0+[-3\kappa,3\kappa]\subset\Gamma_2$ and
\begin{equation}   \vartheta(\gamma_0+u_1)+\vartheta(\gamma_0+u_2)-\vartheta(\gamma_0)\in K^+,
    \qquad |u_1|,|u_2|\le3\kappa.
    \label{eq:pairwise-natural-parameter}
\end{equation}
\end{lemma}

\begin{proof}
Since $\Gamma$ is compact and $\Gamma\subset U$, we can choose compact intervals $\Gamma_1,\Gamma_2$ such that \eqref{eq:gamma_inclusions} holds. The map
\begin{equation*}
    (a,t_1,t_2) \mapsto \vartheta(a+t_1)+\vartheta(a+t_2)-\vartheta(a)
\end{equation*}
is continuous and equals $\vartheta(a)$ at $(t_1,t_2)=(0,0)$. Since $\vartheta(\Gamma_1)$ is a compact subset of the open interval $\Theta$, compactness and uniform continuity give a common $\rho_L>0$ and compact interval $K^+\subset \vartheta(U)$ for which \eqref{eq:local-pairwise-parameter} holds. Decreasing $\rho_L$ if necessary also ensures $a+t_j\in\Gamma_2$. Since $\gamma_0 \in \intr(\Gamma)$, we can pick $\kappa$ sufficiently small such that $3\kappa\le\rho_L$ and $\gamma_0 + [-3\kappa, 3\kappa] \subset \Gamma_2$, which gives \eqref{eq:pairwise-natural-parameter}. All derivatives of $\vartheta$ through order three are uniformly bounded on $\Gamma_2$.

Let $\eta_{a,t}:=\vartheta(a+t)$ and
\begin{equation*}
    \ell_{a,t}(Y):=\log L_{a,t}(Y) =Y(\eta_{a,t}-\eta_a)-\{b(\eta_{a,t})-b(\eta_a)\}.
\end{equation*}
The derivatives of $\ell_{a,t}$ through order three are polynomials of degree at most one in $Y$, with coefficients bounded uniformly over $a\in\Gamma_1$ and $|t|\le\rho_L$. The corresponding derivatives of $L_{a,t}=e^{\ell_{a,t}}$ are $L_{a,t}$ times polynomials of degree at most three in $Y$.

Fix any polynomial $P$, then we have that 
\begin{align*}
\E_a[L_{a,t}^2(Y)P(Y)^2] &= \int P(y)^2
\exp\left( 2\{(\eta_{a,t}-\eta_a)y-b(\eta_{a,t})+b(\eta_a)\}\right) \dint P_a(y)\\
&= \exp\{b(2\eta_{a,t}-\eta_a)-2b(\eta_{a,t})+b(\eta_a)\}
\E_{\vartheta^{-1}(2\eta_{a,t}-\eta_a)}[P(Y)^2].
\end{align*}
We see that squaring the likelihood ratio tilts the natural parameter from $\eta_a$ to $2\eta_{a,t}-\eta_a$, which lies in $K^+$. Since $K^+$ is a compact interval contained in $\Theta$, the exponential-family moment generating function is finite on a uniform neighbourhood of $K^+$, so we have $\sup_{\eta\in K^+}\E_{\vartheta^{-1}(\eta)}[P(Y)^2]<\infty$. 

Since the first three derivatives of $L_{a, t}$ are $L_{a, t}$ multiplied by polynomials in $Y$ of degree at most three, \eqref{eq:local-pairwise-parameter} and the preceding display imply that their $L_2(P_a)$ norms are uniformly bounded over $a \in \Gamma_1, \, |t| \leq \rho_L$. We can therefore apply Taylor's theorem in $L_2(P_a)$ to obtain 
\begin{equation*}
    L_{a, t} = L_{a, 0} + t L'_{a, 0} + \frac{t^2}{2} L''_{a, 0} + O_{L_2(P_a)}(|t|^3).
 \end{equation*}
Noting that $L_{a, 0}'=\ell'_{a,0}=\vartheta'(a)Z_a$ and 
\begin{equation*}
L_{a, 0}''=(\ell'_{a, 0})^2 +\ell''_{a,0}=\vartheta''(a)Z_a +\vartheta'(a)^2\{Z_a^2 - \kappa_2(a)\}
\end{equation*}
proves \eqref{eq:likelihood-expansion}. 
\end{proof}

Our next lemma provides an explicit form for $\E_a[R_a^2]$ from \eqref{eq:projection-resdiual-def}. 

\begin{lemma}[Least favourable nuisance direction]
\label{lemma:projection}
For $a\in\Gamma_1$, define
\begin{equation*}
    c^*(a) :=-\frac{\vartheta''(a)}{2\vartheta'(a)} -\frac{\vartheta'(a)\kappa_3(a)}{2\kappa_2(a)}
\end{equation*}
and 
\begin{equation}
    R_a := H_a + 2 c^*(a)S_a = \vartheta'(a)^2 \left[Z_a^2-\kappa_2(a) -\frac{\kappa_3(a)}{\kappa_2(a)}Z_a \right].
        \label{eq:residual-form}
\end{equation}
Then $\E_a[R_a]=0$, $\E_a[R_a S_a]=0$, and 
\begin{equation}
    \E_a[R_a^2] =\vartheta'(a)^4\left\{\kappa_4(a)+2\kappa_2(a)^2 -\frac{\kappa_3(a)^2}{\kappa_2(a)}\right\}.
    \label{eq:Delta}
\end{equation}
\end{lemma}

\begin{proof}
The exponential family cumulant identities give
\begin{align*}
    \E_aZ_a=0, 
    \qquad\E_aZ_a^2=\kappa_2(a),
    \qquad \E_aZ_a^3=\kappa_3(a),
    \qquad \E_aZ_a^4=\kappa_4(a)+3\kappa_2(a)^2.
\end{align*}
Expanding $H_aS_a=\vartheta'(a)\vartheta''(a)Z_a^2 + \vartheta'(a)^3 Z_a^3 - \vartheta'(a)^3 \kappa_2(a) Z_a$, we have 
\begin{equation*}
\E_a[H_aS_a] = \vartheta'(a) \vartheta''(a) \kappa_2(a) + \vartheta'(a)^3\kappa_3(a),
\end{equation*}
and we also have that $\E_a[S_a^2]=\vartheta'(a)^2\kappa_2(a)$. 
Solving for $\E_a[R_aS_a] = \E_a[H_aS_a] + 2 c \E_a[S_a^2]=0$ therefore gives 
\begin{equation*}
    2c^*(a)=-\frac{\E_a[H_aS_a]}{\E_a[S_a^2]}
    =-\frac{\vartheta''(a)}{\vartheta'(a)}-\vartheta'(a)\frac{\kappa_3(a)}{\kappa_2(a)},
\end{equation*}
which proves the orthogonality and \eqref{eq:residual-form}. Squaring \eqref{eq:residual-form} and applying the cumulant identities shows \eqref{eq:Delta}. 
\end{proof}

Finally, we provide the proof that $\E_a[R_a^2]$ is zero if and only if the GLM response distribution contains exactly two support points. 

\begin{proof}[Proof of \Cref{lemma:two-point-iff}]
If $\E_a[R_a^2]=0$, then
\begin{equation}
    Z_a^2-\kappa_2(a) =\frac{\kappa_3(a)}{\kappa_2(a)}Z_a, \qquad P_a\text{-almost surely},
    \label{eq:Delta-eq-0-condition}
\end{equation}
so $Z_a$ is supported on the roots of a quadratic polynomial. Since $b''(\eta_a)=\mathrm{Var}_a(Y)>0$, there are exactly two support points. 

For the converse statement, suppose that $Y$ has two support points, with $Z_a \in \{z_1, z_2\}$, so that $Z_a^2 - (z_1 + z_2)Z_a + z_1z_2=0$, $P_a$-a.s.. Taking expectations and using $\E_a[Z_a]=0$ gives $z_1z_2 = -\kappa_2(a)$. Multiplying the quadratic identity by $Z_a$ and taking expectations leads to $\eqref{eq:Delta-eq-0-condition}$, which implies $\E_a[R_a^2]=0$. 

Since each of the functions that make up \eqref{eq:Delta} is continuous and $\kappa_2(a)>0$, the map $a\mapsto \E_a[R_a^2]$ is continuous on $U$, and positive whenever $|\operatorname{supp}(Y)| > 2$. Since $\Gamma$ is a compact interval, the extreme value theorem gives a point $\gamma^* \in \Gamma$ such that 
\begin{equation*}
    \inf_{\gamma \in \Gamma} \E_\gamma[R_\gamma^2] = \E_{\gamma^*}[R_{\gamma^*}^2] > 0. \qedhere
\end{equation*}

\end{proof}

\subsection{Distribution constructions and divergence bounds}
\label{app:mixture_calculations}
Throughout this subsection, we impose \Cref{assp:glm-regularity} for the GLM constructions. We can take $X \sim \mathcal{N}(0, \underline{\lambda}I_d)$ for both hypotheses to satisfy \Cref{assp:covariate-dist}. For simplicity, by replacing $X$ and $\beta$ with $X/\sqrt{\underline{\lambda}}$ and $\sqrt{\underline{\lambda}}\beta$, respectively, we may carry out the calculations with a standard Gaussian design. This rescales the separation radius by a fixed model-dependent factor. 

The notation described next will be used throughout this subsection, unless stated otherwise. Fix some $\gamma_0 \in \intr(\Gamma)$, and write $P_{\gamma_0}$ for the null response distribution, and $P_0$ for the null joint law. We denote $c^*:=c^*(\gamma_0)$ and residual $R:=R_{\gamma_0}=H_{\gamma_0} + 2 c^*S_{\gamma_0}$, as defined in \Cref{lemma:projection}. 

We choose a cutoff function by fixing an even smooth function $w:\mathbb R\to[0,1]$ with
\begin{equation*}
    w(u)=1\quad\text{for }|u|\leq\kappa,
    \qquad w(u)=0\quad\text{for }|u|\geq2\kappa,
\end{equation*}
where $\kappa$ is chosen as in \Cref{lemma:local-likelihood-expansion}. For $\beta \in \mathbb{R}^d$ and admissible intercept $\gamma_r$, define the alternative proxy distribution as
\begin{equation}
    Q_{\beta,\gamma_r}(\dint y\mid X=x) := w(x^\top\beta)P_{\beta,\gamma_r}(\dint y\mid X=x) +\{1-w(x^\top\beta)\}P_0(\dint y\mid X=x),
    \label{eq:smooth-proxy-definition}
\end{equation}
where we write $P_{\beta, \gamma_r}(\dint y \mid X=x) = P_{\gamma_r + x^\top \beta} (\dint y)$. 

Let $\pi_{r,s^*}$ be the uniform distribution on vectors with exactly $s^*$ non-zero entries, each equal to $\pm r/\sqrt{s^*}$. For the alternative, we consider the mixture distribution 
\begin{equation}
    \bar{Q}_{r,s^*}^{(n)}:=\int Q_{\beta, \gamma_r}^{\otimes n}\,\pi_{r,s^*}(\dint \beta).
    \label{eq:sparse-alt-mixture}
\end{equation}

Let 
\begin{equation*}
\Lambda_u := \frac{\dint P_{\gamma_0 + u}}{\dint P_{\gamma_0}}, \quad  L_r(u,y):= 1+w(u)\{\Lambda_{u+\gamma_r - \gamma_0}(y)-1\},
\end{equation*}
so that for $X \sim \mathcal{N}(0, I_d)$ and $\|\beta\|_2=r$, we have 
\begin{equation*}
    \frac{\dint Q_{\beta,\gamma_r}}{\dint P_0}(X, Y) = L_r(X^\top \beta, Y). 
\end{equation*}

For $\beta,\beta'\in\mathbb R^d$ with $\|\beta\|_2= \|\beta'\|_2=r$, let $\rho := \beta^\top \beta'$, $U:=X^\top\beta$, $V:=X^\top\beta'$. Then 
\begin{equation*}
    (U,V)\sim\mathcal N\left( 0,
    \begin{pmatrix}
        r^2 & \rho\\
        \rho & r^2
    \end{pmatrix}
    \right).
\end{equation*}
We define the pairwise likelihood ratio overlap by 
\begin{equation}
     M_r(\rho) := \E_\rho\E_{\gamma_0}[L_r(U,Y)L_r(V,Y)],
     \label{eq:pairwise-overlap-rho}
\end{equation}
where $\E_\rho$ denotes expectation under this bivariate Gaussian law and $\E_{\gamma_0}$ is with respect to $Y\sim P_{\gamma_0}$, independently of $(U,V)$. Where it is obvious, we will suppress dependence on $P_{\gamma_0}$.

\begin{lemma}
\label{lemma:gaussian-mixture}
Let $M_r(\rho)$ be as defined in \eqref{eq:pairwise-overlap-rho}. Suppose that, for some fixed $C>0$ and every $0<r\leq r_0$, there exist $a_r\in\mathbb R$ and $B_r\geq0$, with $|a_r|\leq C$, such that
\begin{equation}
    \log M_r(\rho)\leq a_r\rho+C\rho^2+B_r, \qquad |\rho|\leq r^2.
    \label{eq:second-moment}
\end{equation}
There exists a constant $L$, depending only on $C,r_0$, such that, for $\bar{Q}_{r, {s^*}}^{(n)}$, or any other similar mixture alternative construction satisfying \eqref{eq:second-moment}, we have 
\begin{equation*}
 \KL(\bar Q_{r,s^*}^{(n)}\|P_0^{\otimes n})
 \leq nB_r+\frac{(s^*)^2}{d}\left\{\cosh\left(\frac{Lnr^2}{s^*}\right)-1\right\},
 \qquad 0<r\leq r_0.
\end{equation*}
\end{lemma}

\begin{proof}
For independent $\beta,\beta'\sim\pi_{r,s^*}$, write $K:=|S\cap S'|$, where $S, S'$ are the respective supports and $\rho:=\beta^\top\beta'$. Conditional on $S$ and $S'$, $\rho$ is distributed as $(r^2/s^*)\sum_{j=1}^K\zeta_j$, where the $\zeta_j$ are independent Rademacher variables. We can bound 
\begin{align}
    \KL(\bar Q_{r, s^*}^{(n)}\|P_0^{\otimes n})
    &\leq\log\{1+\chi^2(\bar Q_{r, s^*}^{(n)}\|P_0^{\otimes n})\}\notag\\
    &=\log\E_{\beta,\beta'} \left[M_r(\beta^\top\beta')^n \right] \notag \\
    &\leq \log \E \left[\exp(na_r\rho + Cn\rho^2) \right] + nB_r 
\end{align}
where the first inequality is due to Theorem 2.7 of \cite{tsybakov2009introduction} and the second is due to \eqref{eq:second-moment}.

Since $\rho$ is symmetrically distributed, we have 
\begin{align*}
\E[\exp(na_r\rho+Cn\rho^2)]&=\E \left[\exp(Cn\rho^2) \left\{\frac{\exp(na_r\rho)+\exp(-na_r\rho)}{2} \right\}\right]\\
&=\E[\exp(Cn\rho^2) \cosh(na_r\rho)].
\end{align*}

We claim that, for $L\geq\max\{1, \,2C , \,2\sqrt{C/\tanh(1)}, \,2Cr_0^2/\tanh(1)\}$, we can bound 
\begin{equation}
\exp(Cnx^2)\cosh(n|ax|)\leq\cosh(Lnx),
\qquad |a|\leq C,\quad |x|\leq r_0^2.
\label{eq:cosh_bound}
\end{equation}

We can see this from 
\begin{align*}
    \log \cosh (Ln|x|) - \log \cosh (n|ax|) &\geq \int_{Ln|x|/2}^{Ln|x|} \tanh(t) \dint t  \geq \frac{Ln|x|}{2} \tanh(L n|x|/2) \\
    & \geq \tanh(1) \min \left\{\frac{L^2n^2x^2}{4},
                         \frac{Ln|x|}{2} \right\} \geq Cnx^2                        
\end{align*}
where the first inequality is due to $n|ax| \leq Ln|x|/2$ and the second line is due to $\tanh(t)\geq \tanh(1)\min(t,1)$ for $t\geq0$. Rearranging and exponentiating proves \eqref{eq:cosh_bound}.  Since $|\rho|\leq r^2\leq r_0^2$, it follows that
\begin{equation}
\E[\exp(na_r\rho+Cn\rho^2)] \leq\E[\cosh(Ln\rho)]
=\E\left[\left\{\cosh(Lnr^2/s^*)\right\}^K\right].
\label{lower:exp-to-cosh}
\end{equation}

For fixed $k \in \{0, 1, \dots, s^*\}$, we have 
\begin{align*}
    \E \left[ \binom{K}{k} \right] &= \E_S\left\{\sum_{A\subset S,\,|A|=k}\Prob \left(A\subset S' \mid S \right) \right\} = \E_S\left\{\E\left[\sum_{A\subset S,\,|A|=k}\ind\{A\subset S'\}\,\middle| \,S\right]\right\}\\
    &= \binom{s^*}{k} \frac{\binom{d-k}{s^*-k}}{\binom{d}{s^*}} =\binom{s^*}{k} \prod_{i=0}^{k-1}\frac{s^*-i}{d-i} \leq \binom{s^*}{k} \left(\frac{s^*}{d} \right)^k,
\end{align*}
where for $k=0$ we consider the product to evaluate to one. 

Expanding $z^K=\sum_{k=0}^{s^*}\binom Kk(z-1)^k$ with $z=\cosh(Lnr^2/s^*)$, we therefore have that 
\begin{align*}
    \E\left[\left\{\cosh(Lnr^2/s^*)\right\}^K\right] &\leq \sum_{k=0}^{s^*} \left(\cosh(Lnr^2/s^*) - 1  \right)^k  \binom{s^*}{k} \left(\frac{s^*}{d} \right)^k \\
    &= \left\{ 1 + \frac{s^*}{d} \left(\cosh(Lnr^2/s^*) - 1  \right)\right\}^{s^*} \\
    &\leq\exp\left\{\frac{(s^*)^2}{d}\left(\cosh(Lnr^2/s^*) - 1  \right)\right\}.
\end{align*}
Substituting the final expression into \eqref{lower:exp-to-cosh} concludes the proof. 
\end{proof}

\begin{lemma}
     Under \Cref{assp:glm-regularity}, there exist constants $r_0,C>0$ such that, for $0<r\leq r_0$ and $\gamma_r = \gamma_0 + c^*r^2$, we have 
\begin{equation}
\KL(\bar Q_{r,s^*}^{(n)}\|P_0^{\otimes n})
\le Cnr^4+\frac{(s^*)^2}{d}\left\{\cosh\left(\frac{Lnr^2}{s^*}\right)-1\right\}.
\label{eq:glm-mixture-KL}
\end{equation}
and, for each fixed alternative component,
\begin{equation}    \TV(P_{\beta,\gamma_r},Q_{\beta,\gamma_r})\lesssim 2\exp\{-\kappa^2/(2r^2)\}.
    \label{eq:smooth-proxy-TV}
\end{equation}
Moreover, we have that 
\begin{equation}
\TV(P_{\beta, \gamma_r}, P_0) \lesssim r + \exp\{-\kappa^2/(2r^2)\}.
    \label{eq:alternative-fixed-direction-TV}
\end{equation}
\label{lemma:glm-lower}
\end{lemma}

\begin{proof}
By \Cref{lemma:local-likelihood-expansion}, uniformly where $|u| \leq \kappa$, we have 
\begin{equation*}
\Lambda_u = 1+uS_{\gamma_0}+\frac{u^2}{2}H_{\gamma_0}+O_{L_2(P_0)}(|u|^3).
\end{equation*}
Let $A:=U+c^* r^2$, where $U\sim\mathcal N(0,r^2)$. 
Decreasing $\kappa$ and $r_0$ if necessary so that $3\kappa\le\rho_L$ and $|c^*|r_0^2\le\kappa$, by \Cref{lemma:local-likelihood-expansion} we have uniform $L_2(P_0)$ bounds on the first two derivatives of $L_r(\cdot,Y)$. Since $w(U)\neq1$ only when $|U|>\kappa$, which occurs with probability $O(e^{-c/r^2})$, we obtain
\begin{align*}
\E_U[L_r(U,Y)] &= 1+\E[A]S_{\gamma_0}+\frac{\E[A^2]}{2}H_{\gamma_0} +o_{L_2}(r^2)\\
&= 1+\frac{r^2}{2}R +o_{L_2}(r^2).
\end{align*}
Recalling \eqref{eq:pairwise-overlap-rho}, independence of
$U$ and $V$ at $\rho=0$ gives
\begin{align*}
    M_r(0)-1
    &=\E_{\gamma_0}\left[\{\E_U L_r(U,Y)\}^2\right]-1\\
    &=\E_{\gamma_0}\left[\{\E_U L_r(U,Y)-1\}^2\right]\\
    &=\frac{r^4}{4}\E_{\gamma_0}[R^2]+o(r^4)
    =\frac{\E_{\gamma_0}[R^2]}{4}r^4+o(r^4),
\end{align*}
where the second equality uses $\E_{\gamma_0}\E_U[L_r(U,Y)]=1$, and the third follows from the preceding expansion.

By the Cauchy--Schwarz inequality and \Cref{lemma:local-likelihood-expansion}, we have $|\E_0[\partial_u L_r(u, Y)\, \partial_v L_r(v, Y)]| \leq C$ for some constant $C$. By Price's identity \citep[e.g.][Remark 2]{voigtlaender2021price}, we have for $|\rho|<r^2$ that 
\begin{equation*}
M_r'(\rho) = \E_\rho\E_{\gamma_0}\left[\partial_uL_r(U,Y)\partial_vL_r(V,Y) \right],
\end{equation*}
so averaging over $(U, V)$ gives $|M_r'(\rho)| \leq C$. Similar arguments give $|M_r''(\rho)| \leq C$. Since $M_r(0)=1+O(r^4)$, it follows that $M_r(\rho)$ is uniformly bounded away from zero for $|\rho|< r^2$ and sufficiently small $r$. Hence
\begin{equation*}
\sup_{|\rho|< r^2}
|(\log M_r)''(\rho)| = \sup_{|\rho| < r^2} \left| \frac{M_r''(\rho)}{M_r(\rho)}  - \left\{\frac{M_r'(\rho)}{M_r(\rho)}\right\}^2 \right| \leq C, 
\end{equation*}
so taking a Taylor expansion at $\rho=0$ leads to 
\begin{equation}
\log M_r(\rho) \leq a_r\rho+C\rho^2+Cr^4,
\qquad a_r:=(\log M_r)'(0), \qquad |a_r|\leq C.
\label{eq:M-taylor-expansion}
\end{equation}
The choice of smooth cutoff function $w$ and \Cref{lemma:local-likelihood-expansion} imply that $(u,v)\mapsto\E_{\gamma_0}[L_r(u,Y)L_r(v,Y)]$ is bounded and continuous. Since the Gaussian pair $(U,V)$ converges in distribution as $\rho\to\pm r^2$, $M_r(\rho)$ extends continuously to these endpoints. Its lower bound therefore extends as well, so \eqref{eq:M-taylor-expansion} holds for $|\rho|\le r^2$. Applying \Cref{lemma:gaussian-mixture} with $B_r=Cr^4$ proves \eqref{eq:glm-mixture-KL}.

Since $P_{\beta,\gamma_r}$ and $Q_{\beta,\gamma_r}$ have the same marginal distribution for $X$, we have that 
\begin{align*}
\TV(P_{\beta,\gamma_r},Q_{\beta,\gamma_r})
&\leq \E_X\left[ \TV\left( P_{\beta,\gamma_r}(\cdot\mid X), Q_{\beta,\gamma_r}(\cdot\mid X) \right) \right]\\
&\leq \E_X[1-w(X^\top\beta)]\\
&\leq \Prob_X(|X^\top \beta| > \kappa)\\
&\leq 2 \exp(-\kappa^2 / 2r^2).
\end{align*}
where in the third inequality we used that $w(u) = 1$ for $|u| \leq \kappa$. 

Finally, to show the total variation distance bound in \eqref{eq:alternative-fixed-direction-TV}, we start with the triangle inequality 
\begin{equation}
\TV (P_{\beta, \gamma_r} , P_0) \leq \TV(P_{\beta, \gamma_r}, Q_{\beta, \gamma_r}) + \TV(Q_{\beta, \gamma_r}, P_0).
\label{eq:TV-alt-triangle}
\end{equation}
The first term can be bounded by \eqref{eq:smooth-proxy-TV}. To bound the second term, for $t:=X^\top\beta+ c^*r^2$, write $\eta_t:=\vartheta(\gamma_0+t)$ and $\eta_0:=\vartheta(\gamma_0)$.
For a natural exponential family, we have that for all $|t|< 3\kappa $
\begin{align*}
\KL(P_{\gamma_0+t} \parallel P_{\gamma_0})
&= \E_{\gamma_0+t}\left[(\eta_t-\eta_0)Y -{b(\eta_t)-b(\eta_0)}\right]\\
&=(\eta_t-\eta_0)b'(\eta_t) - {b(\eta_t)-b(\eta_0)} \\
&=\int_{\eta_0}^{\eta_t} \{b'(\eta_t )-b'(s)\} \dint s.
\end{align*}
By the choice of $\kappa$, we have $\gamma_0+t\in\Gamma_2$ for $|t|\leq3\kappa$. Since $\vartheta'$ and $b''$ are bounded on the compact intervals $\Gamma_2$ and $\vartheta(\Gamma_2)\subset\Theta$, respectively, where $\Gamma_2$ is as defined in \Cref{lemma:local-likelihood-expansion}, the preceding identity and the mean value theorem give
\begin{equation*}
    \KL(P_{\gamma_0+t}\parallel P_{\gamma_0})
    \leq C|\eta_t-\eta_0|^2
    \leq Ct^2.
\end{equation*}
We can now bound 
\begin{align*}
\KL(Q_{\beta,\gamma_r} \parallel P_0)
&\leq \E_X\left[ w(X^\top\beta) \KL(P_{\gamma_0+t} \parallel P_{\gamma_0})\right]\\
&\leq C\E_X\left[ w(X^\top\beta)(X^\top\beta+ c^* r^2)^2 \right]\\
& \lesssim r^2 +(c^*)^2 r^4 \lesssim r^2.
\end{align*}
Therefore, using Pinsker's inequality and substituting into \eqref{eq:TV-alt-triangle} finishes the proof for \eqref{eq:alternative-fixed-direction-TV}. 
\end{proof}

\begin{lemma}[Construction for GLMs with two support points]
\label{lemma:bernoulli-lower-bounds}
Consider the Gaussian design and prior in \Cref{app:mixture_calculations}. Suppose that the response family has two-point support and satisfies \Cref{assp:glm-regularity}, and fix $\gamma_0\in\intr(\Gamma)$. For all sufficiently small $r>0$, there exists an admissible intercept $\gamma_r$ with $|\gamma_r-\gamma_0|\lesssim r^2$ such that each proxy $Q_{\beta, \gamma_r}$ in \eqref{eq:smooth-proxy-definition}, formed using this intercept, has response marginal $P_{\gamma_0}$. With this choice, there exists $L>0$ such that, for all $n\geq1$ and $1\leq s^*\leq d$,
\begin{equation*}
\KL(\bar Q_{r,s^*}^{(n)}\|P_0^{\otimes n})
\le\frac{(s^*)^2}{d}\left\{\cosh\left(\frac{Lnr^2}{s^*}\right)-1\right\}.
\end{equation*}
Moreover, uniformly over $\|\beta\|_2=r$,
\begin{equation}
    \TV(P_{\beta,\gamma_r},Q_{\beta, \gamma_r}) \leq 2\exp\left(-\frac{\kappa^2}{2r^2}\right),
    \qquad \TV(P_{\beta,\gamma_r},P_0)\lesssim r.
    \label{eq:two-point_proxy}
\end{equation}
\end{lemma}
\begin{proof}
Since the support has two points, after a linear mapping, we may assume that $Y\in\{0,1\}$. Write $p(a):=P_a(Y=1)$ and $p_0:=p(\gamma_0)$, and define
\begin{equation*}
    F(r,\gamma) :=\E\left[w(rG)\{p(\gamma+rG)-p_0\}\right],
    \qquad G\sim\mathcal N(0,1).
\end{equation*}
The proxy has the same response marginal as $P_{\gamma_0}$  when $F(r,\gamma)=0$. By \Cref{assp:glm-regularity}, $p$ and its first two derivatives are bounded on the compact neighbourhood selected by the cutoff function $w$. Differentiating the integrand with respect to $r$ up to twice introduces at most two powers of $G$. The derivatives through order two are therefore bounded by $C(1+|G|+G^2)$, which is integrable. Differentiation under the expectation thus shows that $F$ is twice continuously differentiable near $(0,\gamma_0)$. At $r=0$, we have $F(0,\gamma_0)=0$ and $\partial_\gamma F(0,\gamma_0)=p'(\gamma_0)\neq0$. Hence, the implicit function theorem gives a locally unique twice continuously differentiable function $\gamma(r)$ satisfying
\begin{equation*}
    \gamma(0)=\gamma_0, \qquad F(r,\gamma(r))=0.
\end{equation*}
By the symmetry of $G$, we have that 
\begin{equation*}
    F(-r,\gamma)=\E\left[w(-rG)\{p(\gamma-rG)-p_0\}\right] =\E\left[w(rG)\{p(\gamma+rG)-p_0\}\right]=F(r,\gamma),
\end{equation*}
and by local uniqueness, $\gamma(-r)=\gamma(r)$, so $\gamma'(0)=0$. Setting $\gamma_r:=\gamma(r)$, Taylor's theorem shows that $|\gamma_r-\gamma_0|\lesssim r^2$. Since $\gamma_0\in\intr(\Gamma)$, $\gamma_r$ remains admissible for all sufficiently small $r$.

Let $L_r(u,y):= 1+w(u)\{\Lambda_{u+\gamma_r-\gamma_0}(y)-1\}$. We have that $\E_U[L_r(U, y)]=1$, and thus, when $\rho=0$, recalling the notation in \eqref{eq:pairwise-overlap-rho}, independence of $U$ and $V$ gives
\begin{equation*}
    M_r(0) =\E_{\gamma_0} \left[\{\E_U[L_r(U,Y)]\}^2\right] =1.
\end{equation*}
The derivative bounds used to establish
\eqref{eq:M-taylor-expansion} hold uniformly for
$|\gamma_r-\gamma_0|\lesssim r^2$. Since the exact marginal
matching gives $M_r(0)=1$, the same Taylor expansion argument leads to 
\begin{equation*}
    \log M_r(\rho)\leq a_r\rho+C\rho^2,
    \qquad |a_r|\leq C,\qquad |\rho|\leq r^2.
\end{equation*}
Applying \Cref{lemma:gaussian-mixture} with $B_r=0$ gives the claimed KL bound. The claims in \eqref{eq:two-point_proxy} follow from similar arguments to the proof of \Cref{lemma:glm-lower} with the specified alternative intercept.
\end{proof}

\begin{lemma}[Construction for linear regression]
    Write $P_{\beta, \sigma}$ for the joint law of $X\sim \mathcal{N}(0, I_d)$ and $Y \mid X=x \sim \mathcal{N}(\beta^\top x, \sigma^2)$. Choose a fixed $\sigma_0>0$, depending only on $M_\xi$ and $p_\xi$, such that $\Prob(|\sigma_0G|\leq M_\xi)\geq p_\xi$ for $G\sim\mathcal N(0,1)$, and set $\sigma_r^2:=\sigma_0^2-r^2$. For $r^2<\sigma_0^2/2$, both error distributions $\mathcal N(0,\sigma_0^2)$ and $\mathcal N(0,\sigma_r^2)$ satisfy \Cref{assp:linear-error} and the response marginal agrees with that under $P_{0, \sigma_0}$.
For $n\ge1$ and $1\le s^*\le d$, let
\begin{equation*}
\bar P_{r,s^*}^{(n)}:=\int P_{\beta,\sigma_r}^{\otimes n}\pi_{r,s^*}(\mathrm d\beta).
\end{equation*}
There exists a constant $L>0$, depending only on $\sigma_0$, such that
\begin{equation}
\KL(\bar P_{r,s^*}^{(n)}\|P_{0,\sigma_0}^{\otimes n})
\le\frac{(s^*)^2}{d}\left\{\cosh\left(\frac{Lnr^2}{s^*}\right)-1\right\},
\qquad 0<r^2<\frac{\sigma_0^2}{2}.
\label{eq:linear-mixture-KL}
\end{equation}
Moreover for every fixed $\beta$ with $\|\beta\|_2=r$, we have that $\TV(P_{\beta, \sigma_r}, P_{0, \sigma_0}) \lesssim r$. 
\label{lemma:gaussian-lower}
\end{lemma}

\begin{proof}
Since $\sigma_r\leq\sigma_0$, for $\sigma\in\{\sigma_0,\sigma_r\}$ we have
\begin{equation*}
    \Prob(|\sigma G|\leq M_\xi) \geq \Prob(|\sigma_0G|\leq M_\xi) \geq p_\xi,
\end{equation*}
which verifies \Cref{assp:linear-error}.

For every $\|\beta\|_2=r$, we have $X^\top\beta\sim\mathcal N(0,r^2)$, and hence under $P_{\beta,\sigma_r}$,
\begin{equation*}
    Y\sim\mathcal N(0,r^2+\sigma_r^2)=\mathcal N(0,\sigma_0^2), 
\end{equation*}
i.e.~the response marginal agrees exactly with that under $P_{0, \sigma_0}$.

Writing $U=X^\top\beta$, $V=X^\top\beta'$ and $\rho=\beta^\top\beta'$, we have
\begin{equation*}
\E_0\left[ \frac{\mathrm dP_{\beta, \sigma_r}}{\mathrm dP_{0, \sigma_0}} \frac{\mathrm dP_{\beta', \sigma_r}}{\mathrm dP_{0, \sigma_0}}
    \,\middle|\,X\right]=
\frac{\sigma_0^2}{\sqrt{(\sigma_0^2-r^2)(\sigma_0^2+r^2)}}
\exp\left\{\frac{(U+V)^2}{4(\sigma_0^2+r^2)}
    -\frac{(U-V)^2}{4(\sigma_0^2-r^2)}\right\}.
\end{equation*}
Taking expectations with respect to $X$ gives
\begin{align}
\E_0\left[\frac{\mathrm dP_{\beta, \sigma_r}}{\mathrm dP_{0, \sigma_0}}
    \frac{\mathrm dP_{\beta', \sigma_r}}{\mathrm dP_{0, \sigma_0}} \right] 
&= \frac{\sigma_0^2}{\sqrt{(\sigma_0^2-r^2)(\sigma_0^2+r^2)}}
\left(\frac{\sigma_0^2+r^2}{\sigma_0^2-\rho}\right)^{1/2}
\left(\frac{\sigma_0^2-r^2}{\sigma_0^2-\rho}\right)^{1/2} \notag \\
&=\frac{\sigma_0^2}{\sigma_0^2-\rho}
=\frac{1}{1-\rho/\sigma_0^2},
\label{eq:gauss-pairwise-overlap}
\end{align}
where we used that $U+V$ and $U-V$ are independent centred Gaussians with variances
$2(r^2+\rho)$ and $2(r^2-\rho)$, respectively. Since $|\rho|\leq r^2<\sigma_0^2/2$, we can bound 
\begin{equation*}
    \log \frac{1}{1-\rho/\sigma_0^2}  \leq \frac{\rho}{\sigma_0^2}+C\rho^2.
\end{equation*}
The expectation in \eqref{eq:gauss-pairwise-overlap}  defines the analogous version of $M_r(\rho)$ for this construction. Thus, applying \Cref{lemma:gaussian-mixture} with $Q_\beta=P_{\beta,\sigma_r}$, $P_0=P_{0,\sigma_0}$ and $B_r=0$ proves \eqref{eq:linear-mixture-KL}.

Finally, for any fixed $\beta$ with $\|\beta\|_2=r$, we have that 
\begin{equation*}
    \KL \left(P_{\beta, \sigma_r} \parallel P_{0, \sigma_0} \right) = \frac{1}{2}\log\left( \frac{\sigma_0^2}{\sigma_0^2-r^2} \right) \lesssim r^2,
\end{equation*}
so by Pinsker's inequality, we have $\TV(P_{\beta, \sigma_r}, P_{0, \sigma_0}) \lesssim r$. 
\end{proof}

\subsection{Proof of lower bounds}
\label{app:lower-assembly}
\begin{proof}[Proof of the lower bounds in \Cref{thm:main-result}] Our constructions in \Cref{app:mixture_calculations} lie in the model class after a rescaling of the Gaussian covariates and a corresponding inverse scaling of the regression coefficients, which only changes the separation radius by a fixed constant. 

First, suppose $0 < \epsilon \leq 1$, and define 
\begin{equation}
r_{\mathrm{np}}:=\left\{\frac{s^*}{n} \operatorname{arcosh}\left(1+\frac{d}{(s^*)^2}\right)\right\}^{1/2}  \qquad 
r_{\mathrm{priv}}:=\left\{\frac{s^*}{n} \operatorname{arcosh}\left(1+\frac{d}{n\epsilon^2(s^*)^2}\right)\right\}^{1/2}.
\label{eq:r_np-def}
\end{equation}
For $0\leq\lambda\leq1$, the power series of $\cosh$ gives $\cosh(\lambda x)-1\leq\lambda^2\{\cosh(x)-1\}$. Thus, for a sufficiently small constant $c>0$, \Cref{lemma:gaussian-lower,lemma:bernoulli-lower-bounds}
imply for linear regression and two-point GLMs that 
\begin{equation*}
\KL(\bar Q_{r,s^*}^{(n)}\|P_0^{\otimes n}) \lesssim
\begin{cases}
c^4,&r=c r_{\mathrm{np}},\\
c^4/(n\epsilon^2),&r=c r_{\mathrm{priv}},
\end{cases}
\end{equation*}
where for linear regression we take
$Q_\beta=P_{\beta,\sigma_r}$ and $P_0=P_{0,\sigma_0}$ and use \eqref{eq:sparse-alt-mixture} for the two-point GLM. 

Recalling the notation of \eqref{eq:separation-notation}, note that we can bound 
\begin{align*}
\min\left\{\frac{d^{1/4}}{\sqrt n},
\sqrt{\frac{s^*}{n}}\right\} \lesssim r_{\mathrm{np}} \lesssim \min\left\{\frac{d^{1/4}}{\sqrt n},
\sqrt{\frac{s^*\ell_{n, d, \epsilon}}{n}}\right\}.
\end{align*}
The lower bound follows from $\operatorname{arcosh}(1+x) \gtrsim \min\{\sqrt{x}, 1\}$. As for the upper bound, we first have 
\begin{equation*}
    r_{\mathrm{np}}^2 \leq \frac{s^*}{n} \sqrt{ \frac{2d}{(s^*)^2}} = \frac{\sqrt{2d}}{n},
\end{equation*}
where we used $\cosh(t) \geq 1 + t^2/2$ to bound $\mathrm{arcosh}(1+x) \leq \sqrt{2x}$. 

Secondly, we have that 
\begin{align*}
\frac{s^*}{n}\operatorname{arcosh}\left(1+\frac{d}{(s^*)^2}\right)
&= \frac{s^*}{n}\log\left\{1+\frac{d}{(s^*)^2}
+\sqrt{\frac{d}{(s^*)^2} \left(2+\frac{d}{(s^*)^2}\right)}\right\}\\
&\leq \frac{s^*}{n}\log\left\{2\left(1+\frac{d}{(s^*)^2}\right)\right\}\\
&\leq \frac{s^*}{n}\log\{2(1+d)\}\lesssim \frac{s^* \ell_{n, d, \epsilon}}{n}.
\end{align*}

Similar arguments show that 
\begin{equation*}
   \min\left\{\frac{d^{1/4}}{n^{3/4} \epsilon^{1/2}}, \sqrt{\frac{s^*}{n}}\right\} \lesssim  r_{\mathrm{priv}} \lesssim \min\left\{\frac{d^{1/4}}{n^{3/4}\epsilon^{1/2}}, \sqrt{\frac{s^*\ell_{n, d, \epsilon}}{n}}\right\}.
\end{equation*}

Therefore, for two-point GLMs and linear regression, by \Cref{lemma:bernoulli-lower-bounds} and \Cref{lemma:gaussian-lower} respectively, there exist distributions such that we can bound the quantities in \Cref{prop:general_lower_bound} by 
\begin{equation*}
    \alpha\lesssim r,\qquad
\tau\lesssim\exp\left(-\frac{\kappa^2}{2r^2}\right), \qquad \omega \lesssim
\begin{cases}
c^2/\sqrt n,&r=c r_{\mathrm{np}},\\
c^2/(n\epsilon),&r=c r_{\mathrm{priv}}.
\end{cases}
\end{equation*}

The bounds on $r_{\mathrm{np}}$ and $r_{\mathrm{priv}}$, together with the sample size condition in \Cref{thm:main-result}, ensure that the chosen radii give valid constructions. Taking $c_P$ sufficiently small also makes $r^2\log(n+1)$ sufficiently small, so that $n\tau$ can be made negligibly small.

Taking $r=c_1 r_{\mathrm{np}}$ with $c_1$ sufficiently small allows us to bound 
\begin{equation*}
\frac{1}{2}\left(1-\sqrt n\,\omega-n\tau\right)>\frac{1}{3},
\end{equation*}
proving the lower bound term $r_{\mathrm{np}}$.

For the private contributions, \Cref{prop:general_lower_bound} depends on $\alpha\wedge(\omega+\tau)$. Taking $r=c_2r_{\mathrm{priv}}$ makes $n\epsilon\omega$ sufficiently small. Since $\epsilon\leq1$ and $n\tau$ is sufficiently small, this also controls $n\epsilon(\omega+\tau)$. Alternatively, taking $r=c_3/(n\epsilon)$ and using $\alpha\lesssim r$ makes $n\epsilon\alpha$ sufficiently small. In either case, $\delta\leq c_\delta\epsilon$ gives the corresponding control of $n\delta\{\alpha\wedge(\omega+\tau)\}$. For sufficiently small constants, the private part of \Cref{prop:general_lower_bound} therefore gives the lower bound $\max\{r_{\mathrm{priv}},(n\epsilon)^{-1}\}$.

Combining these three lower bounds, we have that 
\begin{equation*}
r^*(n,d, s^*,\epsilon,\delta,\mathcal P,\infty)
\gtrsim \max\left\{\frac{d^{1/4}}{\sqrt n} \wedge 
\sqrt{\frac{s^*}{n}} , \, \frac{d^{1/4}}{n^{3/4} \epsilon^{1/2}} \wedge \sqrt{\frac{s^*}{n}}, \,\frac{1}{n\epsilon}\right\}.
\end{equation*}

For GLMs with more than two response support points, by \Cref{lemma:glm-lower}, there exist distributions such that we can bound the quantities in \Cref{prop:general_lower_bound} by 
\begin{equation*}
    \alpha\lesssim r,\qquad
\tau\lesssim\exp\left(-\frac{\kappa^2}{2r^2}\right), \qquad \omega \lesssim
\begin{cases}
c^2/\sqrt n,&r=c \min\{r_{\mathrm{np}},n^{-1/4}\},\\
c^2/(n\epsilon),&r=c \min\{r_{\mathrm{priv}},(n \epsilon)^{-1/2}\}.
\end{cases}
\end{equation*}
The same testing
arguments therefore yield
\begin{equation*}
r^*(n,d, s^*, \epsilon,\delta, \mathcal P_{\mathrm{GLM}}(\mathcal F),r_0) \gtrsim \max\left\{r_{\mathrm{np}} \wedge n^{-1/4},\, r_{\mathrm{priv}} \wedge (n\epsilon)^{-1/2},\, \frac{1}{n\epsilon} \right\}.
\end{equation*}
We scale $r_0$ to enforce the construction to be admissible, and the sample size condition ensures that all chosen coefficient norms lie below $r_0$ and that $n\tau$ is sufficiently small.

Finally, if $\epsilon>1$, then, for $D\geq1$,
\begin{equation*}
\frac{D^{1/4}}{n^{3/4}\epsilon^{1/2}} \leq\frac{D^{1/4}}{\sqrt n},
\qquad \frac1{n\epsilon} \leq\min\left\{\frac{D^{1/4}}{\sqrt n}, \sqrt{\frac{s^*}{n}}\right\}.
\end{equation*}
Consequently, the non-private lower bound alone gives the
claimed rate, with $D=d$ or $D=d\wedge n$ as appropriate.
\end{proof}

\section{Technical details of the mean-reduction test}
This appendix contains the proof of \Cref{lemma:transform-choice} and further details of \Cref{subsec:generalised_mean_test}.  

\subsection{Proof of \Cref{lemma:transform-choice}}
We first establish properties of the conditional response mean. Let $s_{\vartheta}=1$ for linear regression and GLMs of form \eqref{eq:GLM-model-def} with non-decreasing $\vartheta$  and set $s_{\vartheta}=-1$ otherwise. Write 
\begin{equation}
    m_\nu(t) := s_{\vartheta}\, \E_{\beta, \nu}[\tanh(Y) \mid X^\top\beta =t],
    \label{eq:conditional_mean}
\end{equation}
where $\nu$ denotes the error distribution in \eqref{eq:linear_regression_model} or the nuisance intercept in \eqref{eq:GLM-model-def}. The conditional distribution is defined through the model response law at predictor $t$. 

\begin{lemma}
    Let $(X, Y) \sim P$, where $P \in \mathcal{P}_{\mathrm{Lin}}$ or $P \in \mathcal{P}_{\mathrm{GLM}}(\mathcal{F})$ for some $\mathcal{F}$ satisfying \Cref{assp:glm-regularity}. We have, uniformly over all nuisance parameters $\nu$, that $|m_\nu|\leq1$ and $m_\nu$ is non-decreasing on $\mathbb R$.
Moreover, there exist fixed constants $t_0,c_m,C_m>0$ such that
\begin{equation*}
    c_m(t-s)\leq m_\nu(t)-m_\nu(s)\leq C_m(t-s), \qquad -t_0\leq s\leq t\leq t_0.
\end{equation*}
    \label{lemma:conditinal-mean}
\end{lemma}

\begin{proof}
\underline{Linear regression transformation.} Here, $m_\nu(t) := \E[\tanh(t + \xi)]$. Since $\tanh$ has bounded derivative, differentiating under the expectation gives 
\begin{equation*}
    m_\nu'(t) = \E[\operatorname{sech}^2(t+\xi)] \leq 1.
\end{equation*}
Therefore, for $|t| \leq t_0$, 
\begin{align*}
    m_\nu'(t) &\geq \operatorname{sech}^2(M_\xi + t_0) \Prob(|\xi| \leq M_\xi) \\
    &\geq p_\xi \operatorname{sech}^2(M_\xi + t_0) =: C_\xi >0. 
\end{align*}
It follows that $m_\nu$ is non-decreasing on $\mathbb{R}$.  

\medskip \noindent \underline{GLM transformation.}
Fix $\eta \in \Theta$ and let $Y$ follow the response distribution with natural parameter $\eta$. Since $\Theta$ is open, exponential moments exist in a neighbourhood of $\eta$, allowing us to interchange differentiation and integration to obtain 
\begin{align*}
    \frac{\dint}{\dint\eta}\E[\tanh(Y)]
    &= \int \tanh(y)\{y-b'(\eta)\} \exp\{\eta y-b(\eta)\}h(y) \dint y\\
    &= \E[Y\tanh(Y)]
    - b'(\eta)\E[\tanh(Y)]\\
    &= \E[Y\tanh(Y)]
    - \E[Y]\E[\tanh(Y)]\\
    &= \mathrm{Cov}(\tanh(Y),Y)\\
    &= \frac{1}{2} \E[(\tanh(Y) - \tanh(Y'))(Y-Y')]>0, 
\end{align*}
where in the final line $Y'$ is an independent copy of $Y$. Since $\vartheta$ is monotone, the choice of $s_\vartheta$ therefore makes $m_\gamma$ non-decreasing on $\mathbb{R}$. Letting $\Gamma_2$ be the compact interval from \Cref{lemma:local-likelihood-expansion}, we have that 
\begin{align*}
    m_\gamma'(t) &= s_\vartheta\vartheta'(\gamma+t) \operatorname{Cov}_{\gamma +t}(\tanh(Y), Y) \\
    &=|\vartheta'(\gamma+t)| \operatorname{Cov}_{\gamma +t}(\tanh(Y), Y).
\end{align*}
Since both factors are continuous and strictly positive, compactness of $\Gamma_2$ implies there exist some $c_m, C_m>0$ such that 
\begin{equation}
    c_m \leq m_\gamma'(t) \leq C_m \qquad \forall (\gamma, t) \text{ with } \gamma+t \in \Gamma_2.
    \label{eq:m_gamma_deriv_lower}
\end{equation}
Choose $t_0 >0$ such that $\Gamma + [-t_0, t_0] \subset \Gamma_2$. Combining both model cases, we obtain fixed constants $t_0,c_m,C_m>0$ such that $|m_\nu|\leq1$, $m_\nu$ is non-decreasing on $\mathbb{R}$, and we have, by the mean value theorem, that 
\begin{equation*}
c_m(t-s)\leq m_\nu(t)-m_\nu(s)\leq C_m(t-s),
 \qquad -t_0\leq s\leq t\leq t_0,
\end{equation*}
uniformly over the allowed nuisance parameters.
\end{proof}

We can now proceed to prove \Cref{lemma:transformation_properties}. 
\begin{proof}[Proof of \Cref{lemma:transformation_properties}]
Since $|\tanh(Y)|\leq1$ almost surely, we have $ZZ^\top\preceq XX^\top$ pointwise, and hence $\E[ZZ^\top]\preceq\Sigma := \E[XX^\top]$. Furthermore, for any unit vector $u$,
\begin{equation*}
    |u^\top\E[Z]| \leq \{\E[(u^\top Z)^2]\}^{1/2}
    \leq \sqrt{u^\top \Sigma u} \leq \sqrt{\overline\lambda}.
\end{equation*}
Taking the supremum over $u$ proves $\|\E[Z]\|_2 \leq \sqrt{\overline{\lambda}}$. The sub-Gaussian bound follows from $|\tanh(Y)|\leq 1$, the sub-Gaussian assumption on the covariates, and the centring inequality in Lemma 2.7.8 of \cite{vershynin2026}.  

Under the null, independence of $X$ and $Y$ gives 
\begin{equation*}
    \E[Z]=\E[X\E[\tanh(Y)\mid X]]= s_{\vartheta} m_\nu(0)\E[X]=0.
\end{equation*}

To show \eqref{eq:mean-to-regression}, let $u:=\beta/\|\beta\|_2$. Since $\E[u^\top X]=0$, we obtain
\begin{align}
    s_\vartheta u^\top\E[Z]
    &=\E\left[u^\top X\{m_\nu(\|\beta\|_2u^\top X)-m_\nu(0)\}\right] \notag \\
    &\geq c_m\|\beta\|_2\,\E\left[(u^\top X)^2\ind\{\|\beta\|_2|u^\top X|\leq t_0\}\right] \notag \\
    &\geq c_m\|\beta\|_2\left\{u^\top\Sigma u-\frac{\|\beta\|_2^2}{t_0^2}\E[(u^\top X)^4]\right\} \notag \\
    &\geq c_m\|\beta\|_2\left\{\underline\lambda-\frac{C\|\beta\|_2^2C_X^4}{t_0^2}\right\}.
    \label{eq:transformed-mean-directional}
\end{align}
The first inequality is due to \Cref{lemma:conditinal-mean} and the final inequality follows from the sub-Gaussian assumption on the covariates. Thus, for a sufficiently small absolute constant $c_0>0$, we have that $s_\vartheta u^\top\E[Z]\geq c_m\underline\lambda/2\|\beta\|_2$ provided that $\|\beta\|_2\leq c_0t_0\sqrt{\underline\lambda}/C_X^2$. 

Since $m_\nu$ is non-decreasing, the function $t\mapsto u^\top x\{m_\nu(tu^\top x)-m_\nu(0)\}$ is non-decreasing for every fixed $x$. Therefore, when $\|\beta\|_2>c_0t_0\sqrt{\underline\lambda}/C_X^2$, evaluating the expectation at this threshold gives the required constant lower bound. Combining the two cases gives
\begin{equation*}
    \|\E[Z]\|_2\geq s_\vartheta u^\top\E[Z]
    \geq c\min\left\{\|\beta\|_2,\frac{t_0\sqrt{\underline\lambda}}{C_X^2}\right\} \geq c_\mu \min\{\|\beta\|_2, 1\}
\end{equation*}
where $c_\mu>0$ depends only on the model constants.
\end{proof}

\subsection{Proof of \Cref{lemma:restricted-sensitivity}}
\begin{proof}[Proof (adapted from Lemma 44 of \citealp{narayanan2022private})]
Let $S\subseteq[n]$ be the set of changed records. We can decompose 
\begin{align*}
    Q_{[n]}(z_{1:n}) &= Q_{S^c}(z_{1:n}) + 2\sum_{i \in S, j \in S^c} z_i ^\top z_j + Q_S(z_{1:n})\\
    &=Q_{S^c}(z_{1:n})+2\sum_{i\in S}u_i(z_{1:n})-Q_S(z_{1:n}).
\end{align*}
This can be used to bound 
\begin{align*}
    |Q_{[n]}(z_{1:n})-Q_{[n]}(z_{1:n}')| &\leq
    2\sum_{i\in S}\{|u_i(z_{1:n})|+|u_i(z_{1:n}')|\} +|Q_S(z_{1:n})|+|Q_S(z_{1:n}')|\\
    &\leq 4KL\{\sqrt{nd}+ns\} +2L\{K\sqrt{Kd}+K^2\}\\
    &\leq 6KL\{\sqrt{nd}+ns+K\},
\end{align*}
where $K\leq n$ was used in the last step.
\end{proof}

\subsection{Proof of typical set probability and concentration}
\label{app:typical-set}
\begin{lemma}
\label{lemma:typical_set_prob}
Let $(X_1,Y_1),\ldots,(X_n,Y_n)$ be independent observations with common distribution $P$ belonging either to $\mathcal P_{\mathrm{Lin}}$ or to $\mathcal P_{\mathrm{GLM}}(\mathcal F)$ for a fixed response family $\mathcal F$ satisfying
\Cref{assp:glm-regularity}. Fix $0<\eta\leq1/4$. Suppose either that $\mathbb{E}[Z]=0$ or that
$s\leq\|\mathbb{E}[Z]\|_2\leq2s$, for some $0<s \leq \sqrt{\overline{\lambda}}$. There exists a constant $C$, depending only on the fixed model constants, such that, for
\begin{equation}
    L = C\left\{1+\log\left(\frac{12(n+1)}{\eta}\right)\right\},
    \label{eq:typical-L}
\end{equation}
we have, with probability at least $1-\eta$, that 
$(Z_1,\ldots,Z_n) \in\mathcal{D}_{n,s,L}$ and 
\begin{align*}
    \left|Q_{[n]}(Z)-n(n-1)\|\mathbb{E}[Z]\|_2^2\right| \leq L\{n\sqrt d+n^{3/2}s\}.
\end{align*}
\end{lemma}

\begin{proof}
    Denoting $R_i:=Z_i - \E[Z]$, we can write 
    \begin{equation*}
    u_i=(n-1)\|\mathbb{E}[Z]\|_2^2+(n-1)\langle R_i,\mathbb{E}[Z]\rangle +\mathbb{E}[Z]^\top \sum_{j \neq i} R_j+ R_i^\top \sum_{j \neq i} R_j.
\end{equation*}
Since $R_i$ and $\sum_{j\neq i} R_j$ are independent with $\|R_i\|_{\psi_2} \leq c$ and $\|\sum_{j\neq i} R_j\|_{\psi_2} \leq c\sqrt{n}$, we have by \Cref{lemma:sub-Gaussian-inner-product} that 
\begin{equation*}
\Prob \left(|u_i-(n-1)\|\E Z\|_2^2| >C\left\{\|\E Z\|_2 n\sqrt t +\|\E Z\|_2\sqrt{nt}+\sqrt n(\sqrt{dt}+t)\right\} \right)\leq3e^{-t}.
\end{equation*}
Taking $t=\log(12n/\eta)$ and applying a union bound over $i\in[n]$ shows that \eqref{eq:row-bound} holds for every $i\in[n]$ with probability at least $1-\eta/4$.
    
We next consider the probability of the event in \eqref{eq:subset-bound}. For any fixed set $S$ such that $|S|=k$, we have  
    \begin{equation}
        Q_S = k(k-1)\|\mathbb{E}[Z]\|_2^2 + 2(k-1) \mathbb{E}[Z] ^\top \sum_{i \in S} R_i + \sum_{i, j \in S, \, i \neq j} R_i^\top R_j
        \label{eq:typical-terms}
    \end{equation}
    We introduce a sequence of Bernoulli random variables $\zeta := \{\zeta_i : i \in S\}$, $\zeta_i \overset{\text{i.i.d.}}{\sim}$ Bernoulli($1/2$), in order to use a decoupling argument to bound the third term in \eqref{eq:typical-terms}. Let 
    \begin{align*}
        C_\zeta := A_\zeta ^\top B_\zeta := \left \langle \sum_{i\in S} \zeta_i R_i, \sum_{j\in S} (1-\zeta_j) R_j \right \rangle;
    \end{align*}
    note that we have 
    \begin{equation*}
        D:= \sum_{i, j \in S, i \neq j} R_i^\top R_j = 4\mathbb{E}_\zeta [C_\zeta \mid \{R_i : i \in S\}]. 
    \end{equation*} 
    
    Fix a realisation of $\zeta$.  The vectors $A_\zeta$ and $B_\zeta$ are centred, and are independent conditional on $\zeta$. Proposition 2.7.1 of \cite{vershynin2026} gives
 \begin{equation*}
    \|A_\zeta\|_{\psi_2} \leq C\left(\sum_{i\in S:\zeta_i=1}\|R_i\|_{\psi_2}^2\right)^{1/2}
    \leq CC_X\sqrt k,
\end{equation*}
and similarly $\|B_\zeta\|_{\psi_2}\leq CC_X\sqrt k$.
Lemma~\ref{lemma:sub-Gaussian-inner-product} therefore implies, for every $p\geq1$, that
\begin{equation*}
 \|C_\zeta\|_{L^p(R\mid\zeta)} \leq C_1k(\sqrt{dp}+p),
\end{equation*}
uniformly over $\zeta$. Therefore, by Jensen's inequality, we have 
\begin{align*}
 \E_R[|D|^p]
 &= 4^p\E_R \left| \E_\zeta[C_\zeta\mid \{R_i : i \in S\}]\right|^p \\
 &\leq 4^p\E_R\E_\zeta  \left[ |C_\zeta|^p\mid \{R_i : i \in S\}\right]\\
 &= 4^p\E_\zeta\E_R|C_\zeta|^p\\
 &\leq\left\{C_1k(\sqrt{dp}+p)\right\}^p.
\end{align*}
Applying Markov's inequality for $t>1$ with the bound $\|D\|_{L^t(R)}
 \leq C_1k(\sqrt{dt}+t)$ leads to 
\begin{equation}
 \Prob\left(
 \left|\sum_{i, j\in S: i\ne j} R_i^\top R_j \right|
 > C_2 k\{\sqrt{dt}+t\}\right)\leq e^{-t}.
 \label{eq:typical-decoupled_tail}
\end{equation}
The second term in \eqref{eq:typical-terms} has sub-Gaussian norm at most $C\|\E[Z]\|_2k^{3/2}$. Since $\|\E[Z]\|_2\leq2s$, combining its tail bound with \eqref{eq:typical-decoupled_tail} gives
\begin{equation}
    \Prob\left(
    |Q_S-k(k-1)\|\E[Z]\|_2^2|
    >C\{k\sqrt{dt}+kt+sk^{3/2}\sqrt t\}
    \right)\leq3e^{-t},\qquad t>1.
    \label{eq:typical-centred-subset}
\end{equation}
Adding the deterministic term $k(k-1)\|\E[Z]\|_2^2\leq4k^2s^2$ yields
\begin{equation}
 \Prob \left(|Q_S|\leq c_3
 \{k\sqrt{dt}+kt+ s k^{3/2}\sqrt t+k^2s^2\} \right) \geq 1-3e^{-t}
 \label{eq:typical-fixed_subset}
\end{equation}
for all $t \geq 1$. 

Now for each $k$, there are $\binom{n}{k} \leq (en/k)^k$ subsets of $[n]$ that are of size $k$. Taking $t_k:=k\log(en/k)+\log(12n/\eta)$, we have $t_k\lesssim k\{1+\log(12(n+1)/\eta)\}$. Since $s\leq\sqrt{\overline\lambda}$, the bound in \eqref{eq:typical-fixed_subset} is therefore at most $L\{k\sqrt{kd}+k^2\}$ for sufficiently large $C$ in \eqref{eq:typical-L}. A union bound over all subsets gives
\begin{equation*}
    \Prob(\mathcal A^c) \leq 3\sum_{k=1}^n\left(\frac{en}{k}\right)^k e^{-t_k} \leq \frac{\eta}{4},
\end{equation*}
where $\mathcal A$ denotes the event that \eqref{eq:subset-bound} holds for every $S\subseteq[n]$.

Finally, applying \eqref{eq:typical-centred-subset} with $S=[n]$ and $t=\log(12/\eta)$ gives
\begin{equation*}
    |Q_{[n]}-n(n-1)\|\E[Z]\|_2^2| \leq C\{n\sqrt{dt}+nt+sn^{3/2}\sqrt t\}
    \leq L\{n\sqrt d+n^{3/2}s\}
\end{equation*}
except on an event of probability at most $\eta/4$. Combining the three failure probabilities proves the claim.
\end{proof}

\subsection{Proof of \Cref{lemma:single_scale_test}}
\begin{proof}
\noindent \underline{Privacy.} For any $z, z'$ that are different in one observation, we have, by the triangle inequality, that 
\begin{align*}
\widehat{G}(z) &= \inf_{\tilde{z} \in \mathcal{D}_{n, s, L}} \left\{ G(\tilde{z}) + D \Ham(z, \tilde{z}) \right\} \\
&\leq \inf_{\tilde{z} \in \mathcal{D}_{n, s, L}} \left\{G(\tilde{z}) + D\Ham(z, z') + D \Ham(z', \tilde{z}) \right\}\\
= \widehat{G}(z') + D \Ham(z, z').
\end{align*}
Interchanging the roles of $z$ and $z'$ gives $|\widehat{G}(z) - \widehat{G}(z')| \leq D \Ham(z, z')$, so the privacy guarantee follows from \Cref{lemma:laplace-test,lemma:post-processing}. 

\medskip \noindent \underline{Error guarantee.} We will first show that $G$ is $D$-Lipschitz on $\mathcal{D}_{n, s, L}$ with respect to Hamming distance. 
    
Note that the definition of $D$, with the lower bound on $s$ in \eqref{eq:IT-test_s-radius}, implies we have 
    \begin{equation}
        6L(\sqrt{nd} + sn + 8 \log(2/\eta)/\epsilon_0) \leq D. 
        \label{eq:lip-const-bound}
    \end{equation}
    Suppose two datasets $z_{1:n}, z_{1:n}' \in \mathcal{D}_{n, s, L}$ differ in $K$ records. By \Cref{lemma:restricted-sensitivity}, we have that 
    \begin{equation*}
        |Q(z_{1:n})-Q(z_{1:n}')| \leq 6KL(\sqrt{nd} + sn + K).
    \end{equation*}
    First suppose $K \leq (8 \log(2/\eta))/\epsilon_0$. Since clipping is 1-Lipschitz, we have by \eqref{eq:lip-const-bound} that $$|G(z_{1:n}) - G(z_{1:n}')| \leq |Q(z_{1:n})-Q(z_{1:n}')| \leq KD.$$ Otherwise, if $K > (8 \log(2/\eta))/\epsilon_0$, the clipping allows us to bound $$|G(z_{1:n}) - G(z_{1:n}')| \leq B_n \leq KD.$$ Combining both cases, we have that $|G(z_{1:n}) - G(z_{1:n}')| \leq D \Ham(z_{1:n}, z_{1:n}')$ for all $z_{1:n}, z_{1:n}' \in \mathcal{D}_{n, s, L}$. 

    Since $\mathcal{D}_{n, s, L}$ is non-empty, by the McShane extension theorem (\Cref{thm:mcshane-extension}), $\widehat{G}$ as defined in \eqref{eq:mcshane-extension} satisfies $\widehat{G}(z_{1:n}) = G(z_{1:n})$ for all $z_{1:n} \in \mathcal{D}_{n, s, L}$ and has global Lipschitz constant $D$.

    By \Cref{lemma:typical_set_prob}, under either hypothesis, we have that $Z \in \mathcal{D}_{n, s, L}$ with probability at least $1-\eta/2$, and that 
    \begin{equation*}
        \left|Q(Z_{1:n}) - n(n-1)\|\mathbb{E}[Z]\|_2^2 \right| \leq L \left(n\sqrt{d} + n^{3/2} s \right).
    \end{equation*}
    Therefore, on this typical event under $H_0$, we have that $G(Z_{1:n}) \leq B_n/4$, while under $H_1(s)$, we have $G(Z_{1:n}) \geq 3B_n/4$; this is due to the separation condition on $s$ in \eqref{eq:IT-test_s-radius}. In addition, we have that 
    \begin{equation*}
        \Prob(|W|\geq B_n/4)=\exp\{-\epsilon_0B_n/(4D)\}=(\eta/2)^2\leq\eta/2,
    \end{equation*}
    so a union bound proves the error guarantee. 
\end{proof}

\subsection{Multiscale test}
\label{app:multiscale_test}
\begin{lemma}[Multiscale test]
\label{lemma:private-multiscale-reduction}
Let $\mathcal P$ be a class of distributions equipped with a mean parameter $\mathbb{E}[Z](P)\in\mathbb R^d$, and suppose that $\|\mathbb{E}[Z](P)\|_2\leq B$ for all $P\in\mathcal P$
for some known $B>0$.  Suppose that, for every $s\in(0,B]$, $0<\epsilon_0$ and $0<\eta_0\leq1/4$, there exists an $(\epsilon_0,0)$-differentially private test $\phi_{s,\epsilon_0,\eta_0}$ of
\begin{equation*}
    H_0:\mathbb{E}[Z](P)=0
    \qquad\text{against}\qquad
    H_1(s):s\leq\|\mathbb{E}[Z](P)\|_2\leq2s
\end{equation*}
such that, whenever $s\geq\rho(n, d, \epsilon_0, \eta_0)$, its errors satisfy
\begin{equation*}
    \sup_{\substack{P\in\mathcal P\\\mathbb{E}[Z](P)=0}}
    \Prob_P(\phi_{s,\epsilon_0,\eta_0}=1) \leq\eta_0,
    \qquad \sup_{\substack{P\in\mathcal P\\ s\leq\|\mathbb{E}[Z](P)\|_2\leq2s}}
    \Prob_P(\phi_{s,\epsilon_0,\eta_0}=0) \leq\eta_0.
\end{equation*}

Fix $0<r\leq B$, $\epsilon > 0$ and $0<\eta\leq1/4$. Provided that 
\begin{equation}
    r \geq \rho \left(n, d, \frac{\epsilon}{J}, \frac{\eta}{J} \right), \quad J:=\max\left\{1,\left\lceil\log_2\left(\frac {B}{r}\right)\right\rceil\right\},
    \label{eq:multiscale-condition-generic}
\end{equation}
there exists an $(\epsilon,0)$-differentially private test $\Phi_r$ of
\begin{equation*}
    H_0:\mathbb{E}[Z]=0 \qquad\text{against}\qquad H_1:\|\mathbb{E}[Z]\|_2\geq r
\end{equation*}
whose type-I and type-II errors are each at most $\eta$.
\end{lemma}
\begin{proof}
Set $s_j:=2^jr$, $\epsilon_j:=\epsilon/J$ and $\eta_j:=\eta/J$ for $j=0,\ldots,J-1$, and reject whenever at least one $\phi_{s_j,\epsilon_j,\eta_j}$ rejects. By the composition property of DP (\Cref{lemma:composition}), the test is $(\epsilon,0)$-DP. Under the null, a union bound gives type-I error at most $\sum_j\eta_j=\eta$. If $\|\E[Z]\|_2\geq r$, there is a $j\in\{0,\ldots,J-1\}$ such that
\begin{equation*}
    s_j\leq\|\E[Z]\|_2\leq2s_j,
\end{equation*}
because $\|\E[Z]\|_2\leq B$. Since $s_j\geq r$ and \eqref{eq:multiscale-condition-generic} holds, the corresponding single-scale test has type-II error at most $\eta_j\leq\eta$.
\end{proof}

\section{Technical details of the sparse upper bound}
\label{app:sparse}
\begin{proof}[Proof of \Cref{prop:sparse-upper}] We take the choice of $\eta_{\mathrm{net}}:= \min\{1/2, c_c/(2C_c)\}$, where $c_c, C_c$ are the constants from \Cref{lemma:sparse-cov}. For $v\in\mathcal V_k$, write
    \begin{equation*}
 K_{j,v}:=\frac{1}{2}\left\{\tanh(v^\top X_{2j-1})-\tanh(v^\top X_{2j})\right\}\left\{\tanh(Y_{2j-1})-\tanh(Y_{2j})\right\}, \qquad \bar{K}_v:=\frac1m\sum_{j=1}^m K_{j,v},
\end{equation*}
and 
\begin{equation*}
    A(D) := \max_{k \in \mathcal K} \left\{\max_{v \in \mathcal{V}_k} \bar{K}_v - q_k \right\}, \qquad q_k := \sqrt{\frac{8\log(2|\mathcal{K}| |\mathcal{V}_k|/\eta)}{m}}.
\end{equation*}

\noindent \underline{Privacy.} Since $K_{j, v} \in [-2, 2]$, the sensitivity of $\bar{K}_v$ is $4/m$ for any $v$. Let $(k_*,v_*)$ be the values that attain the maximum in $A(D)$. For neighbouring datasets $D,D'$,
\begin{align*}
 A(D)=\overline K_{v_*}(D)-q_{k_*}
 \leq\overline K_{v_*}(D')-q_{k_*}+\frac{4}{m}
 \leq A(D')+\frac{4}{m}.
\end{align*}
Interchanging $D$ and $D'$ gives $|A(D)-A(D')|\leq 4/m$. \Cref{lemma:laplace-test} and \Cref{lemma:post-processing}
 therefore imply that $\phi_{\mathrm{sparse}}$ is $(\epsilon, 0)$-DP. 

\medskip \noindent \underline{Type-I error.} By independence of the observations, we have that 
\begin{equation*}
    \E[K_{j,v}]=\cov\{\tanh(v^\top X),\tanh(Y)\}.
\end{equation*}
Under the null, this expectation is zero, so Hoeffding's inequality and a union bound give 
\begin{equation*}
    \Prob(A(D)>0) \leq\sum_{k\in\mathcal K} |\mathcal{V}_k|\exp\left(-\frac{mq_k^2}{8}\right) =\frac{\eta}{2}.
\end{equation*}
Since $\Prob(W > \log(1/\eta)) = \eta/2$, the type-I error is at most $\eta$.

\medskip \noindent \underline{Type-II error.}  Suppose that $\|\beta\|_2>0$, and $\|\beta\|_0 \leq s^*$. Let $k^* \in \mathcal{K}$ be such that $s^*\leq k^* \leq 2s^*$. By \Cref{lemma:sparse-cov}, one of $u^* = \beta / \|\beta\|_2$ or $u^*=-\beta / \|\beta\|_2$ satisfies  
\begin{equation*}
    \cov \left(\tanh((u^*)^\top X),\tanh(Y) \right) \geq c_c\min\{\|\beta\|_2,1\}.
\end{equation*}

By construction, there exists $v^* \in \mathcal{V}_{k^*}$ such that $\|v^* - u^*\|_2 \leq \eta_{\mathrm{net}}$. By \Cref{lemma:sparse-cov}, we have that 
\begin{align}
    \E [K_{1,v^*}] &=\cov \left(\tanh((v^*)^\top X),\tanh(Y) \right) \notag \\
 &\geq\cov\left(\tanh((u^*)^\top X),\tanh(Y) \right) -C_c\|v^*-u^*\|_2\min\{\|\beta\|_2,1\} \notag \\
 &\geq(c_c-C_c\,\eta_{\mathrm{net}})\min\{\|\beta\|_2,1\} \notag \\
 &\geq\frac{c_c}{2}\min\{\|\beta\|_2,1\},
 \label{eq:sparse-sep-lower}
\end{align}
by noting that since $\tanh$ is odd, \eqref{eq:sparse-cov-continuity} applies around either sign of $\beta/\|\beta\|_2$. 

Hoeffding's inequality implies that  
\begin{equation}
    \Prob_\beta  \left(\bar{K}_{v^*} < \E[K_{1, v^*}] - q_{k^*}  \right) \leq \exp\left(-\frac{m q_{k^*}^2}{8}\right) =\frac{\eta}{2|\mathcal{K}||\mathcal{V}_{k^*}|}.
    \label{eq:scan-type-ii}
\end{equation}

By construction of the nets $\mathcal{V}_S$, we have 
\begin{align}
 \log|\mathcal V_{k^*}|
 &\leq k^*\log(ed/k^*)+k^*\log(1+2/\eta_{\mathrm{net}}) \notag\\
 &\leq2s^*\log(ed/s^*)+2s^*\log(1+2/\eta_{\mathrm{net}}) \lesssim s^*\log(ed/s^*).
 \label{eq:net-card-bounds}
\end{align}

Applying the separation condition in \eqref{eq:adaptive-scan-power} with $c$ sufficiently small to \eqref{eq:sparse-sep-lower} therefore gives 
\begin{align*}
    \E[K_{1, v^*}] > 2 q_{k^*} + \frac{8 \log(1/\eta)}{m\epsilon}.
\end{align*}

Thus, combining \eqref{eq:scan-type-ii} and that $\mathbb{P}(W < - \log(1/\eta)) = \eta/2$, the type-II error is upper bounded by $\eta$. 
\end{proof}

\begin{lemma} Suppose $(X,Y)$ follows a model in $\mathcal P_{\mathrm{Lin}}$ or $\mathcal P_{\mathrm{GLM}}(\mathcal F)$, where $\mathcal{F}$ satisfies \Cref{assp:glm-regularity}. 
There exist model-dependent constants $c_c,C_c>0$ such that, for every $\beta= \|\beta\|_2u$ with $\|\beta\|_2>0$, $\|u\|_2=1$, and every unit vector $v$, we have that 
    \begin{equation}
        \left|\cov \left(\tanh(u^\top X), \tanh(Y) \right) \right| \geq c_c \min\{\|\beta\|_2,  1\}
        \label{eq:sparse-cov-signal}
    \end{equation}
and 
\begin{equation}
    \left|\cov \left(\tanh(u^\top X)-\tanh(v^\top X),\tanh(Y) \right) \right| \leq C_c \|v-u\|_2 \min\{\|\beta\|_2,1\}.
 \label{eq:sparse-cov-continuity}
\end{equation}
\label{lemma:sparse-cov}
\end{lemma}

\begin{proof}
    Let $V = u^\top X$. Throughout, $c_m, C_m, t_0$ are the constants from \Cref{lemma:conditinal-mean}. 
    
    By \Cref{assp:covariate-dist}, we have that $\E[V]=0$, $\underline{\lambda} \leq\E [V^2] \leq\overline\lambda$, and $\E[V^4] \leq M_4$ for a fixed constant $M_4$. 
    Since $\E[V]=0$, by H\"older's inequality we have that 
    \begin{equation*}
        \E \left[\{V\}_+\right] = \E[\{V\}_-] = \frac{1}{2} \E[|V|] \geq \frac{1}{2} \sqrt{\frac{\E[V^2]^3}{\E[V^4]}} \geq c_V :=\frac{\underline{\lambda}^{3/2}}{2\sqrt{M_4}} > 0.
    \end{equation*}
    We can bound 
    \begin{align*}
        \E \left[ \{V\}_+ \ind \left\{V \geq c_V/2 \right \}\right] = \E[\{V\}_+] - \E \left[ \{V\}_+ \ind \left\{V< c_V/2 \right \}\right]  \geq c_V/2,
    \end{align*}
    and by Cauchy--Schwarz, we have that 
    \begin{align*}
         \E \left[ \{V\}_+ \ind \left\{V \geq c_V/2 \right \}\right]^2 \leq \E[\{V\}_+^2] \Prob \left(V \geq  c_V/2 \right),
    \end{align*}
    which implies that 
    \begin{equation*}
    \Prob \left( V \geq c_V/2 \right) \geq p_V, \qquad p_V := \frac{\underline{\lambda}^3}{16 M_4 \overline{\lambda}}. 
    \end{equation*}
    Likewise, we have $\Prob \left( V \leq -c_V/2 \right) \geq p_V$. 

    With $V'$ as an independent copy of $V$, and recalling the definition of $m_\nu$ from \eqref{eq:conditional_mean}, we have 
    \begin{align*}
     \cov(\tanh(V),m_\nu(\|\beta\|_2V))
     &=\frac{1}{2}\E[(\tanh(V)-\tanh(V'))(m_\nu(\|\beta\|_2 V)-m_\nu( \|\beta\|_2 V'))]\\
     &\geq 2p_V^2 \tanh(c_V/2) \left\{m_\nu(\|\beta\|_2 c_V/2)-m_\nu(-\|\beta\|_2 c_V/2) \right\}\\
     &\geq 4p_V^2 c_m \tanh(c_V/2)\min\{\|\beta\|_2 c_V/2,t_0\} \gtrsim\min(\|\beta\|_2,1),
     \end{align*}
     where the final inequality is due to \Cref{lemma:conditinal-mean}. Conditioning on $X$, we have     
    \begin{equation*}
        \cov \left(\tanh(V), m_\nu(\|\beta\|_2 V) \right) = s_{\vartheta} \cov \left(\tanh(u^\top X), \tanh(Y) \right)
    \end{equation*}
    which proves \eqref{eq:sparse-cov-signal}. 

    Conditioning on $X$, and applying the Cauchy--Schwarz inequality, we have 
    \begin{align}
       &\left|\cov \left(\tanh(u^\top X)-\tanh(v^\top X),\tanh(Y) \right) \right| = \left|\cov \left(\tanh(u^\top X)-\tanh(v^\top X), m_\nu(\|\beta\|_2 V )\right) \right| \notag\\
       & \hspace{4cm} \leq  \var \left(\tanh(u^\top X)-\tanh(v^\top X)  \right) ^{1/2} \var\left(m_\nu(\|\beta\|_2 V) \right)^{1/2}.
       \label{eqref:covariance-cauchy}
    \end{align}

    Since $\tanh$ is $1$-Lipschitz, we have that 
    \begin{align*}
        \var \left(\tanh(u^\top X)-\tanh(v^\top X) \right) \leq \E \left[ \left((u-v)^\top X \right)^2 \right] \leq \overline{\lambda}\|u-v\|_2^2. 
    \end{align*} 
    
   It remains to show that $\var(m_\nu(\|\beta\|_2 V)) \lesssim \min\{1, \|\beta\|_2\}^2$.  If $\|\beta\|_2 \leq 1$, then \Cref{lemma:conditinal-mean} implies that
    \begin{align*}
         \E|m_\nu(\|\beta\|_2V)-m_\nu(0)|^2
         &\leq C_m^2 \|\beta\|_2^2\E[V^2]+4\Pr(\|\beta\|_2|V|>t_0)\\
         &\leq(C_m^2+4/t_0^2)\overline\lambda \|\beta\|_2^2.
    \end{align*}
    Otherwise, if $\|\beta\|_2 \ge 1$, we use  $|m_\nu| \leq 1$ to obtain that the left-hand side of the preceding display is at most four. Since 
    \begin{equation*}
         \var\left(m_\nu(\|\beta\|_2V)\right)\leq\E|m_\nu(\|\beta\|_2V)-m_\nu(0)|^2,
    \end{equation*}
    combining the two cases yields the second variance bound in \eqref{eqref:covariance-cauchy}, which proves \eqref{eq:sparse-cov-continuity}. 
\end{proof}

\section{Technical details on response-only tests}
\Cref{app:nuisance-orthog-test} states the precise form of the $g_a$ function used within \Cref{alg:generic-radial-test} and proves its properties. \Cref{app:radial-test} proves the privacy and error guarantees of \Cref{alg:generic-radial-test}. \Cref{app:mc-approx} shows how a Monte Carlo approximation can be implemented to have negligible effect on the rates. Finally, \Cref{app:exponential_regression} provides a formal justification for the claims in \Cref{remark:exponential-kernel-test}. 

\subsection{Proof of \Cref{lemma:bounded-orthog-contrast}}
\label{app:nuisance-orthog-test}
For $a \in \Gamma_1$, with $\Gamma_1$ as defined in \Cref{lemma:local-likelihood-expansion}, we define our test function as follows. 
Let 
\begin{equation}
    A_{a,M}(y):=[R_a(y)]_{-M}^{M},
    \qquad  M:=\max\left\{1,\,2\sqrt{\frac{K_4}{\Delta_{\min}^2}},\,\frac{8K_SK_R\sqrt{K_4}}{c_m\Delta_{\min}^2}\right\},
    \label{eq:truncated-residual-def}
\end{equation}
where 
\begin{equation*}
    K_4 := \sup_{a\in\Gamma_1}\E_a[R_a^4], 
    \qquad K_R := \sup_{a\in\Gamma_1}\|R_a\|_{L_2(P_a)}, 
    \qquad K_S  := \sup_{a\in\Gamma_1}\|S_a\|_{L_2(P_a)},
\end{equation*}
and, defining $m(a) := \E_a[\tanh(Y)]$ and $\bar{\psi}_a(y) := \tanh(y)-m(a)$, the constants
\begin{equation*}
    c_m := \inf_{a\in\Gamma_1} |m'(a)|>0,
    \qquad \Delta_{\min}^2 :=\inf_{a\in\Gamma_1}\E_a[R_a^2]>0,
\end{equation*}
where $\Gamma_1$ is the compact interval defined in \Cref{lemma:local-likelihood-expansion}. The first three constants are well-defined due to \Cref{lemma:exp-moments}. The positivity of $c_m$ follows from \eqref{eq:m_gamma_deriv_lower}. The positivity of $\Delta_{\min}^2$ follows from \Cref{lemma:two-point-iff} applied on the compact interval $\Gamma_1 \subset U$. 

Finally, define 
\begin{align*}
    \lambda_{a, M} &:=\frac{\E_a[A_{a,M}(Y)S_a(Y)]}{m'(a)}, \qquad 
    g_{a}(y):=A_{a,M}(y)-\E_a[A_{a,M}(Y)]-\lambda_{a,M}\bar\psi_a(y).
\end{align*}

\begin{proof}
We have that $\E_a[g_a(Y)]=0$ by construction, since  $\E_a[\bar\psi_a(Y)]=0$. Further, from 
\begin{align*}
    m'(a)&=\int \tanh(y) p_a(y) \partial_a \log(p_a(y)) \dint y \\
    &= \E_a[\{\tanh(Y) - m(a)\} S_a(Y)] + m(a) \E_a[S_a(Y)]\\
    &=\E_a[\bar{\psi}_a(Y)S_a(Y)],
\end{align*}
we see that the orthogonality claim in $\eqref{eq:g-orthogonal}$ follows from 
\begin{align*}
    \E_a[g_a(Y)S_a(Y)] = \E_a[A_{a,M}(Y)S_a(Y)] -\lambda_{a,M}\E_a[\bar\psi_a(Y)S_a(Y)]=0.
\end{align*}

We now show that this choice of $g_a$ preserves a positive second-order component. Recalling from \Cref{lemma:projection} that $R_a = H_a + 2c^*(a)S_a$ and since $\E_a[g_a(Y) S_a(Y)]=0$, we have 
\begin{align*}
    \E_a[g_a(Y) H_a(Y)] = \E_a[g_a(Y)R_a(Y)]&= \E_a[R_a^2] + \E_a[(g_a - R_a)R_a]\\
    & =\E_a[R_a^2] +\E_a[(A_{a,M}-R_a)R_a]-\lambda_{a,M}\E_a[\bar{\psi}_a R_a]\\
    &\geq \Delta_{\min}^2 + \E_a[(A_{a, M} - R_a)R_a] - \lambda_{a, M} \E_a[\bar{\psi}_a R_a],
\end{align*}
where the second equality is due to $\E_a[R_a]=0$. 

The truncation error satisfies 
\begin{align*}
    \E_a[(A_{a, M}-R_a)R_a] &\geq - \E_a[R_a^2 \ind\{|R_a|>M\}] \\
    &\geq -\frac{\E_a[R_a^4]}{M^2} \geq \frac{-K_4}{M^2}.
\end{align*}

For the correction term, we can bound  
\begin{align*}
    |\lambda_{a, M}||\E_a[\bar{\psi}_a R_a]| & \leq \frac{|\E_a[A_{a, M} S_a]|}{|m'(a)|} \times  2 \E_a[|R_a|] \\
    & \leq 2 K_R \frac{|\E_a[(A_{a, M} - R_a)S_a]|}{c_m}\\
    & \leq \frac{2 K_R K_S}{c_m} \E_a[(A_{a,M}-R_a)^2]^{1/2} \leq \frac{2K_RK_S \sqrt{K_4}}{c_m M}
\end{align*}
where the first inequality is from $|\tanh| \leq 1$ implying that $|\bar{\psi}_a| \leq 2$, the second from recalling that $\E_a[R_aS_a]=0$ by \Cref{lemma:projection}, and the final inequality is from bounding 
\begin{align*}
    \E_a[(A_{a,M}-R_a)^2]^{1/2} &\leq \E_a[R_a^2 \ind\{|R_a|>M\}]^{1/2}\\
    &\leq \E_a[R_a^4/M^2]^{1/2} \leq \sqrt{K_4}/M.
\end{align*}

We therefore have 
\begin{align*}
    \E_a[g_aH_a] \geq \Delta_{\min}^2 - \frac{K_4}{M^2} - \frac{2K_RK_S \sqrt{K_4}}{c_m M} \geq \frac{\Delta_{\min}^2}{2}
\end{align*}
uniformly over $a \in \Gamma_1$, by our choice of $M$ in \eqref{eq:truncated-residual-def}, which proves \eqref{eq:g-curvature} with $c_g := \Delta_{\min}^2/2$. 

Finally, we have 
\begin{align*}
    |g_a(y)| &\leq |A_{a, M}(y)| +  |\E_a[A_{a, M}(Y)]| + 2 |\lambda_{a, M}| \\
    &\leq 2M + \frac{2K_S \sqrt{K_4}}{c_mM}
\end{align*}
uniformly over $y$ and $a \in \Gamma_1$, so we can take the final expression to be $B_g$ for the bound in \eqref{eq:g-bounded}. 
\end{proof}

\begin{lemma}
\label{lemma:exp-moments}
Under Assumption~\ref{assp:glm-regularity}, and with $\Gamma_2$ as defined in \Cref{lemma:local-likelihood-expansion}, there exists $s_0>0$ such that
\begin{equation}
    \sup_{a\in\Gamma_2}\E_a\exp\{s_0|Y|\}<\infty.
    \label{eq:uniform-exp-moment}
\end{equation}
Consequently, for every fixed integer $k\ge1$, $\sup_{a\in\Gamma_2}\E_a|Y|^k<\infty$. 
\end{lemma}

\begin{proof}
The compact set $K_0:=\vartheta(\Gamma_2)$ lies in the interior of the open interval $\Theta$. Hence there exists $s_0>0$ such that $K_0+[-s_0,s_0]\subset\Theta$. 
For $Y\sim P_{\vartheta^{-1}(\eta)}$ and $\eta\in K_0$,
\begin{equation*}
    \E_{\vartheta^{-1}(\eta)} e^{s_0Y} =\exp\{b(\eta+s_0)-b(\eta)\},
    \qquad \E_{\vartheta^{-1}(\eta)} e^{-s_0Y} =\exp\{b(\eta-s_0)-b(\eta)\}.
\end{equation*}
The two right-hand sides are uniformly bounded because $b$ is finite and continuous on the compact set $K_0+[-s_0,s_0]$. Since $\exp(s_0 |Y|) \leq \exp(s_0 Y)+ \exp(-s_0 Y)$, we obtain \eqref{eq:uniform-exp-moment}. The polynomial moment bound follows from $|y|^k\leq C_{k,s_0} \exp(s_0|y|)$.
\end{proof}

\subsection{Proof of \Cref{prop:generic-radial-test}}
\label{app:radial-test}

\begin{lemma}
    Let $\widehat{\gamma}$ be the private intercept estimator in \eqref{eq:private-intercept-def}. There exist constants $C_{\gamma}, r_0>0$ such that for all $\|\beta\|_2 < r_0$ and $\eta \in (0, 1/4)$, it holds that 
    \begin{equation*}
        \Prob \left(|\widehat{\gamma} - \gamma | > C_{\gamma}  \left\{ \|\beta\|_2^2 + a_n(\epsilon, \eta) \right\} \right) \leq \frac{\eta}{2}, \qquad a_n(\epsilon, \eta) := \sqrt{\frac{\log(8/\eta)}n} +\frac{\log(8/\eta)}{n\epsilon}.
    \end{equation*}
    \label{lemma:gamma-error}
\end{lemma}

\begin{proof}
By \Cref{lemma:conditinal-mean}, $m(a)$ is strictly monotone with $|m'| \geq c_m>0$ on $\Gamma_1$. Since $m(\widehat\gamma)$ is the projection of $\widetilde m$ onto $m(\Gamma_1)$, we have
\begin{equation*}
    |\widehat\gamma-\gamma| \leq c_m^{-1}|m(\widehat\gamma)-m(\gamma)| \leq c_m^{-1}|\widetilde m-m(\gamma)|.
\end{equation*}
It therefore suffices to obtain a high probability bound on 
    \begin{equation}
        |\widetilde{m} - m(\gamma)| \leq |\widetilde{m}  - \E_{\beta, \gamma}[\tanh(Y)]| + |\E_{\beta, \gamma}[\tanh(Y)] - m(\gamma)|,
        \label{eq:gamma-error-decomp}
    \end{equation}
where we recall the definition of $\widetilde{m}$ from \eqref{eq:tildem-def}. 
We first bound the bias term. By \Cref{lemma:local-likelihood-expansion}, we have, over $a\in\Gamma$ and $|u|\leq t_0$ for some fixed $t_0>0$, that  
\begin{align*}
    m(a+u) &= \E_a[\tanh(Y) L_{a, u}(Y)]\\
    &= \E_a[\tanh(Y)] + u\E_a[\tanh(Y) S_a(Y)] + \frac{u^2}{2} \E_a[\tanh(Y) H_a(Y)] + o(u^2)\\
    &= m(a) + m'(a) u + O(u^2).
\end{align*}
This remainder bound extends to all $u\in\mathbb R$: for $|u|>t_0$, boundedness of $m$ and uniform boundedness of $|m'|$ on $\Gamma$ give
\begin{equation*}
    |m(\gamma+u)-m(\gamma)-m'(\gamma)u|
    \leq 2+C|u|\leq C'u^2.
\end{equation*}
Thus, uniformly over $\gamma\in\Gamma$, the remainder is bounded by $C'u^2$ for every $u\in\mathbb R$.

Since $\E[X^\top\beta]=0$, it follows that
\begin{align*}
    |\E_{\beta,\gamma}[\tanh(Y)]-m(\gamma)|
    &=\left|\E_X\left[m(\gamma+X^\top\beta)-m(\gamma)-m'(\gamma)X^\top\beta\right]\right|\\
    &\leq C\E_X[(X^\top\beta)^2]
    \lesssim\|\beta\|_2^2.
\end{align*}
By Hoeffding's inequality, we have that 
\begin{align*}
    \Prob \left( \left|\frac{1}{n_1} \sum_{i=1}^{n_1} \tanh(Y_i) - \E_{\beta, \gamma}[\tanh(Y)] \right| > \sqrt{\frac{2 \log(8/\eta)}{n_1}} \right) \leq \frac{\eta}{4}.
\end{align*}
We also have that  $\Prob\left(|W_1|>2\log(4/\eta)/(n_1\epsilon)\right)=\eta/4$. A union bound therefore bounds the first term of \eqref{eq:gamma-error-decomp} by $C a_n(\epsilon,\eta)$ with probability at least $1-\eta/2$.
\end{proof}

\begin{proof}[Proof of \Cref{prop:generic-radial-test}] 

\medskip 
\noindent \underline{Privacy guarantee.}
The first stage is $(\epsilon, 0)$-DP since the empirical mean of $\tanh(Y)$ has global sensitivity $2/n_1$. If the replaced observation is in the first block, the second stage can be viewed as a post-processing of the first stage output with the observations for the second stage held fixed. Alternatively, if it is in the second block, then the output from the first stage has the same distribution under both datasets. For every fixed $a\in\Gamma_1$, the second stage empirical mean has sensitivity at most $2B_g/n_2$, and is therefore also $(\epsilon, 0)$-DP on the second half of the data. Hence, the test is $(\epsilon, 0)$-DP with respect to the full dataset. 

\medskip
\noindent\underline{Concentration result.}
Let $\mathcal{F}_1$ contain the first-half observations and the first-stage Laplace noise. Conditional on $\mathcal F_1$, the variables $g_{\widehat{\gamma}}(Y_i)$, $i\in \{n_1+1, \dots, n\}$, are independent and satisfy
$|g_{\widehat{\gamma}}(Y_i)|\le B_g$. Hence Hoeffding's inequality and the independent
Laplace noise imply that, for some fixed constant $C_{\mu}>0$,
\begin{equation}
    \Prob\left(|\widetilde{T}-\mu(\widehat{\gamma})|>C_{\mu} a_n(\epsilon, \eta)  \,\middle|\,\mathcal F_1\right)\le \eta/2,
    \qquad \mu(a):=\E_{\beta,\gamma}[g_{a}(Y)].
    \label{eq:radial-conditional-concentration}
\end{equation}

\medskip
\noindent \underline{Type-I error.} Under the null, we have by \Cref{lemma:gamma-error} that $|\gamma - \widehat{\gamma}| \leq C_\gamma a_n(\epsilon, \eta)$ with probability at least $1-\eta/2$. Under this event, \Cref{lemma:plugin-expansion} gives $|\mu(\widehat{\gamma})| \lesssim a_n(\epsilon, \eta)^2$, so combining with \eqref{eq:radial-conditional-concentration}, we have, for $c_T$ chosen sufficiently large, that $\Prob_{0, \gamma}\{\tilde{T} > c_T a_n( \epsilon, \eta)\} \leq \eta$. 

\medskip
\noindent \underline{Type-II error.} Suppose that $C_1 a_n(\epsilon,\eta)\leq \|\beta\|_2^2 \leq r_0^2$. By \Cref{lemma:gamma-error}, with probability at least $1-\eta/2$, 
\begin{equation}
    |\gamma - \widehat{\gamma}| \lesssim \|\beta\|_2^2 + a_n(\epsilon, \eta) \lesssim \|\beta\|_2^2
    \label{eq:gamma_event}
\end{equation}
where the last inequality holds provided that $C_1$ is sufficiently large. Choose $v_0$ sufficiently small that the remainder bound in \Cref{lemma:plugin-expansion} satisfies $Cv_0\leq c_g/4$. We condition on \eqref{eq:gamma_event} and take $r_0$ sufficiently small that $|\gamma -\widehat\gamma|\leq v_0/2$. Then \Cref{lemma:plugin-expansion} gives
\begin{equation*}
    \E_{\gamma+X^\top\beta} \bigl[g_{\widehat\gamma}(Y)\bigr] \geq \frac{c_g}{4}
    \bigl(\gamma-\widehat\gamma+X^\top\beta\bigr)^2, \qquad  |\gamma-\widehat\gamma+X^\top\beta|
    \le v_0.
\end{equation*}
We can therefore bound 
\begin{align}
    \mu(\widehat\gamma) &:= \E_X \E_{\gamma + X^\top \beta}[g_{\widehat{\gamma}}(Y)] \notag\\
    &\geq \frac{c_g}{4} \E_X\left[\bigl(\gamma-\widehat\gamma+X^\top\beta\bigr)^2
    \ind\left\{|\gamma-\widehat\gamma+X^\top\beta| \le v_0\right\}\right] -B_g \Prob\left( |\gamma-\widehat\gamma+X^\top\beta|>v_0\right)
    \label{eq:mu_gamma_bound}
\end{align}
where the inequality is due to the fact that $|g_{\widehat{\gamma}}(y)|\leq B_g$ by construction. 

We have that 
\begin{equation*}
    \E_X\left[\bigl(\gamma-\widehat\gamma+X^\top\beta\bigr)^2\right] = (\gamma - \widehat{\gamma})^2 + \beta^\top \E[X X^\top] \beta \geq \underline{\lambda}\|\beta\|_2^2,
\end{equation*}
and 
\begin{equation*}
    \Prob_X(|\gamma-\widehat{\gamma} + X^\top \beta| > v_0) \leq \Prob(|X^\top \beta| > v_0/2) \lesssim \exp(-c/\|\beta\|_2^2).
\end{equation*}
Therefore, applying Cauchy--Schwarz to bound 
\begin{equation*}
    \E_X\left[\bigl(\gamma-\widehat\gamma+X^\top\beta\bigr)^2 \ind\left\{|\gamma-\widehat\gamma+X^\top\beta| > v_0\right\}\right] \lesssim \|\beta\|_2^2 \exp(-c/\|\beta\|_2^2)
\end{equation*}
 and substituting into \eqref{eq:mu_gamma_bound}, we obtain 
\begin{equation*}
    \mu(\widehat{\gamma}) \geq \frac{c_g}{4} \left\{ \underline{\lambda}\|\beta\|_2^2  -C \|\beta\|^2_2 \exp(-c/\|\beta\|^2_2)\right\} - C B_g \exp(-c/\|\beta\|_2^2) \gtrsim \|\beta\|_2^2.
\end{equation*}
Since $\|\beta\|_2^2 \geq C_1 a_n(\epsilon, \eta)$, taking $C_1$ sufficiently large ensures that $\mu(\widehat{\gamma}) - C_\mu a_n(\epsilon, \eta) \geq c_T a_n(\epsilon, \eta)$. 
Hence, \eqref{eq:radial-conditional-concentration} together with the probability bound $1-\eta/2$ for \eqref{eq:gamma_event} gives $\Prob_{\beta, \gamma}(\tilde{T} \leq c_T a_n) \leq \eta$. 
\end{proof}

\begin{lemma}
\label{lemma:plugin-expansion}
Fix $a \in \Gamma_1$, where $\Gamma_1$ is as defined in \Cref{lemma:local-likelihood-expansion}, and let $g_a$ be as defined in \Cref{lemma:bounded-orthog-contrast}. There exist fixed constants $v_0,C>0$, such that, uniformly over $a\in\Gamma_1$ and $|v|\le v_0$,
\begin{equation*}
    \left|\E_{a+v}[g_a(Y)] - \frac{1}{2} \E_a[g_a(Y) H_a(Y)]  v^2 \right| \leq  C|v|^3, \quad \E_a[g_a(Y) H_a(Y)] \geq c_g>0.
\end{equation*}
\end{lemma}
\begin{proof}
The uniformly bounded third $L_2(P_a)$ derivative in \Cref{lemma:local-likelihood-expansion} gives
\begin{equation*}
    L_{a, v}(Y) := \frac{\dint P_{a +v}}{\dint P_a}=1+vS_a+\frac{v^2}{2}H_a+ R_{a,v},
    \qquad \|R_{a,v}\|_{L_2(P_a)}\le C|v|^3
\end{equation*}
on a common neighbourhood of zero. Since $\E_{a+v}[g_a(Y)]=\E_a[g_a L_{a,v}]$, \eqref{eq:g-centred} and \eqref{eq:g-orthogonal} give
\begin{equation*}
    \E_{a+v}[g_a] =\frac{v^2}{2}\E_a[g_aH_a]+\E_a[g_a R_{a,v}],
\end{equation*}
where by \eqref{eq:g-bounded}, we have $|\E_a(g_a R_{a,v})| \le B_g\|R_{a,v}\|_{L_2(P_a)}\le CB_g|v|^3$, which completes the proof.
\end{proof}

\subsection{Monte Carlo approximations}
\label{app:mc-approx}
We will show that the inversion step in \Cref{alg:generic-radial-test} and the relevant $g_a$ function can be approximated sufficiently well by simulations. 

Fix a grid $\mathcal{A}_n \subset \Gamma_1$ that is a $n^{-1/2}$ net of $\Gamma_1$, with cardinality $O(\sqrt{n})$. For each $a \in \mathcal{A}_n$,  generate $N=\lceil Cn|\mathcal A_n|/\eta\rceil$ independent draws $V_{a,1},\ldots,V_{a,N}\sim P_a$ and compute
\begin{align*}
\widehat m(a)&=\frac1N\sum_{\ell=1}^N\tanh(V_{a,\ell}),&
\widehat d(a)&=\frac1N\sum_{\ell=1}^N\tanh(V_{a,\ell})S_a(V_{a,\ell}),\\
\widehat c(a)&=\frac1N\sum_{\ell=1}^NA_{a,M}(V_{a,\ell}),&
\widehat q(a)&=\frac1N\sum_{\ell=1}^NA_{a,M}(V_{a,\ell})S_a(V_{a,\ell}),
\end{align*}
where the draws are independent between the different $a$'s and of the observed data and privacy-preserving noise. Denote their respective population versions by $m, d, c$, and $q$. Set 
\begin{equation*}
\widehat{\gamma}_{\mathrm{MC}}\in \argmin_{a \in \mathcal{A}_n} | \widehat{m}(a) - \tilde{m}|,
\end{equation*}
recalling the definition of $\tilde{m}$ in \eqref{eq:tildem-def}. Further define 
\begin{equation*}    \widehat{g}_{\widehat{\gamma}_{\mathrm{MC}}}(y)=\left[A_{\widehat{\gamma}_{\mathrm{MC}},M}(y)-\widehat{c}(\widehat{\gamma}_{\mathrm{MC}})-\frac{\widehat q(\widehat{\gamma}_{\mathrm{MC}})}{\widehat{d}(\widehat{\gamma}_{\mathrm{MC}})}
\{\tanh(y)-\widehat{m}(\widehat{\gamma}_{\mathrm{MC}})\}\right]_{-B_g}^{B_g}
\end{equation*}
if $|\widehat{d}(\widehat{\gamma}_{\mathrm{MC}})| \geq c_m/2$, where $c_m$ is the constant from \Cref{lemma:conditinal-mean}, and set $\widehat{g}_{\widehat{\gamma}_{\mathrm{MC}}} \equiv 0$ otherwise.

\begin{lemma}
    \label{lemma:mc-approx}
    Suppose that independent draws from $P_a$ and evaluations of $S_a(y)$ and $R_a(y)$ are available at unit cost, uniformly over $a\in\Gamma_1$ and $y\in\mathcal{Y}$. Implementing the test in \Cref{alg:generic-radial-test} using $\widehat{\gamma}_{\mathrm{MC}}$ and $\widehat{g}_{\widehat{\gamma}_{\mathrm{MC}}}$ in place of $\widehat{\gamma}$ and $g_{\widehat{\gamma}}$ maintains the same privacy and testing guarantees at the separation rate in \Cref{prop:generic-radial-test} and has a runtime of $O(n^2/\eta)$. 
\end{lemma}

\begin{proof}
For the privacy guarantee, we may condition on the simulation draws, which are independent of the data, so that the release of $\widehat{\gamma}_{\mathrm{MC}}$ can be viewed as a post-processing of $\tilde{m}$, and the $B_g$ bound holds by the construction of $\widehat{g}$, which ensures privacy for the second stage. The runtime for the function evaluations on $\mathcal{A}_n$ is $O(|\mathcal{A}_n|N)=O(n^2/\eta)$, while the grid search and second-stage evaluation cost is $O(|\mathcal{A}_n| + n) = O(n)$. 

    Let $\mathcal{G}$ denote the event 
    \begin{equation*}
        \mathcal{G} := \left\{\max_{a\in\mathcal A_n}\max_{h\in\{m,d,c,q\}} |\widehat h(a)-h(a)| \leq c_0n^{-1/2}\right\}, \quad c_0 := \min\{1, c_m/4\}.
    \end{equation*}
    Since $\tanh$ and $A_{a,M}$ are uniformly bounded and $\sup_{a\in\Gamma_1}\E_a[S_a(Y)^2]<\infty$, the four functions being averaged, evaluated at $Y\sim P_a$, have variances bounded uniformly over $a\in\Gamma_1$. With $C$ in the definition of $N$ chosen to be sufficiently large, we have by Chebyshev's inequality and a union bound that $\Prob(\mathcal{G}) \geq 1 - \eta/4$. 

    We work on $\mathcal{G}$ for the remainder of the proof. Let $\widehat{\gamma}$ be the minimiser in \eqref{eq:private-intercept-def}, and choose $a_{\widehat{\gamma}} \in \mathcal{A}_n$ such that $|a_{\widehat{\gamma}} - \widehat{\gamma}| \leq n^{-1/2}$. Since $m(\widehat\gamma)$ is the projection of $\widetilde m$ onto the interval $m(\Gamma_1)$, and $0<c_m \leq |m'|\leq C_m<\infty $ on $\Gamma_1$, by definition of  $\widehat{\gamma}_\mathrm{MC}$, we have that 
    \begin{align*}
        |\widehat\gamma_{\mathrm{MC}}-\widehat\gamma| 
        &\leq \frac{1}{c_m} \left\{ |m(\widehat\gamma_{\mathrm{MC}})-\widetilde m|-|\widetilde m - m(\widehat\gamma)| \right\} \\
        &\leq \frac{1}{c_m} \left\{ |m(\widehat\gamma_{\mathrm{MC}})-\widehat m(\widehat\gamma_{\mathrm{MC}})|
     +|\widehat m(\widehat\gamma_{\mathrm{MC}})-\widetilde m| -|\widetilde m - m(\widehat\gamma)| \right\}\\
    &\leq \frac{1}{c_m}\left\{\sup_{u\in\Gamma_1}|m'(u)|+2c_0\right\}n^{-1/2} \lesssim n^{-1/2}.
    \end{align*}
    Let $a \in \mathcal{A}_n$ and $y \in \mathcal{Y}$ be arbitrary. Without considering the clipping, we have 
    \begin{align*}
        &\left| A_{a,M}(y)- \widehat c(a)-\frac{\widehat q(a)}{\widehat d(a)}
         \{\tanh(y)-\widehat m(a)\}-g_a(y) \right|\\
        & \qquad \leq \left|c(a)-\widehat c(a) \right| 
         + \left|\left\{\frac{\widehat q(a)}{\widehat d(a)}-\frac{q(a)}{d(a)}\right\}
         \{\tanh(y)-\widehat m(a)\} \right|
        + \left|\frac{q(a)}{d(a)}\{\widehat m(a)-m(a)\}  \right|.
    \end{align*}
    On $\mathcal{G}$, the first and last terms are of order $n^{-1/2}$. We can also bound the middle term of the preceding display by
    \begin{align*}
        \left|\frac{\widehat q(a)}{\widehat d(a)}-\frac{q(a)}{d(a)}\right| &\leq \left|\frac{\widehat q(a)-q(a)}{\widehat d(a)} \right|
        + \left|\frac{q(a)\{d(a)-\widehat d(a)\}}{\widehat d(a)d(a)} \right|\\
        &\leq \frac{4}{3c_m}|\widehat q(a)-q(a)| +\frac{4|q(a)|}{3c_m^2}|\widehat d(a)-d(a)| \lesssim n^{-1/2}
    \end{align*}
    where the second inequality is due to 
    \begin{equation*}
        |\widehat d(a)|\geq |d(a)|-|\widehat d(a)-d(a)| \geq c_m-c_0n^{-1/2}\geq 3c_m/4.
    \end{equation*}

    In particular, the bounds hold for $a =\widehat{\gamma}_{\mathrm{MC}}$, and since $g_a(y) \in [-B_g, B_g]$, we have that 
    \begin{equation}
    \label{eq:estimated-constrast}
        \sup_{y \in \mathcal{Y}} \left| \widehat{g} _{\widehat{\gamma}_{\mathrm{MC}}}(y) - g_{\widehat{\gamma}_{\mathrm{MC}}} (y)\right| \lesssim n^{-1/2}.
    \end{equation}

    By a union bound considering the event $\mathcal{G}$ and \Cref{lemma:gamma-error}, we have that 
    \begin{align*}
        |\widehat\gamma_{\mathrm{MC}}-\gamma|
        &\leq|\widehat\gamma_{\mathrm{MC}}-\widehat\gamma|+|\widehat\gamma-\gamma|\\
        &\lesssim n^{-1/2}+\|\beta\|_2^2+a_n(\epsilon,\eta/2)
        \lesssim \|\beta\|_2^2+a_n(\epsilon,\eta)
    \end{align*}
    with probability at least $1-\eta/2$. Using \eqref{eq:estimated-constrast}, we can apply the same concentration argument in \eqref{eq:radial-conditional-concentration}. Following the remaining arguments for the type-I and type-II error bounds in the proof of \Cref{prop:generic-radial-test} gives the separation rate, with constants adjusted.     
\end{proof}

\subsection{Technical details of the exponential regression test}
\label{app:exponential_regression}
In this subsection, we give a generic $U$-statistic based test for detecting heterogeneous parameters in the response distributions. Then, in \Cref{lemma:exponential-kernel}, we show that the kernel construction in \Cref{remark:exponential-kernel-test} satisfies the requirements needed for this test. 

\begin{proposition}
Let $h$ be a symmetric kernel with range contained within an interval of length $H$. Suppose there exists some constant $c>0$ such that we have
\begin{equation}
    \E_0[h(Y_1, Y_2) ]=0 \quad \mathrm{and} \quad \E_\beta[h(Y_1, Y_2)] \geq c \min\{\|\beta\|_2^2, 1\}
    \label{eq:radial_kernel_condition}
\end{equation}
uniformly over the nuisance class. Define $U := \binom{n}{2}^{-1} \sum_{i<j} h(Y_i, Y_j)$. Then the test 
\begin{equation}
    \phi := \ind\left\{ U + \frac{2H}{n\epsilon}W \geq H \sqrt{\frac{2 \log(2/\eta)}{n}} + \frac{2H \log(1/\eta)}{n\epsilon}\right\}, \quad W \sim \mathrm{Lap}(1)
    \label{eq:exponential-test}
\end{equation}
is $(\epsilon, 0)$-DP, for any $\eta \in (0, 1/4)$. Furthermore, if 
\begin{equation*}
    c \min\{\|\beta\|_2^2, 1\} \geq 4H \left\{ \sqrt{\frac{\log(1/\eta)}{n}} + \frac{\log(1/\eta)}{n\epsilon} \right\},
\end{equation*}
then the sum of the type-I and type-II errors is upper-bounded by $2\eta$. 
\label{proposition:radial_test}
\end{proposition}

\begin{proof}
\underline{Privacy.} If two datasets differ in one observation, then at most $n-1$ of the terms in the sum defining $U$ can change. Since the range of $h$ is contained in an interval of length $H$, each such term changes by at most $H$, and thus the sensitivity of $U$ is $2H/n$. It follows from the Laplace mechanism and the post-processing property of DP that $\phi$ is $(\epsilon, 0)$-DP. 

\medskip \noindent \underline{Type-I error.}
We next control the type-I error. Changing one observation changes $U$
by at most $2H/n$, so McDiarmid's inequality gives
\begin{equation*}
    \Prob\left(U-\E[U]\geq t\right)
    \leq\exp\left(-\frac{nt^2}{2H^2}\right).
\end{equation*}
Under the null, $\E[U]=0$, and hence
\begin{equation*}
    \Prob\left(U \geq H\sqrt{\frac{2\log(2/\eta)}{n}}\right)
    \leq \frac{\eta}{2}.
\end{equation*}
Moreover, since $W\sim\operatorname{Lap}(1)$,
\begin{equation*}
    \Prob\left(\frac{2H}{n\epsilon}W \geq \frac{2H\log(1/\eta)}{n\epsilon} \right) =  \frac{\eta}{2}.
\end{equation*}
A union bound therefore gives $\Prob(\phi=1) \leq \eta$. 

\medskip \noindent \underline{Type-II error.}
 Applying McDiarmid's inequality to the lower tail gives
\begin{equation*}
    \Prob\left(U- \mathbb{E}[U] \leq -H\sqrt{\frac{2\log(2/\eta)}{n}}\right)\leq\frac{\eta}{2}.
\end{equation*}
Combining the above with the Laplace tail bound, we have that the type-II error is at most $\eta$ whenever
\begin{equation*}
    \E[U] \geq 2H\sqrt{\frac{2\log(2/\eta)}{n}}+
    \frac{4H\log(1/\eta)}{n\epsilon},
\end{equation*}
and recalling the assumption in \eqref{eq:radial_kernel_condition} concludes the proof. 
\end{proof}

Given \Cref{proposition:radial_test}, it suffices to show that the kernel in \eqref{eq:exponential-test-kernel} satisfies the properties in \eqref{eq:radial_kernel_condition} for the exponential regression model. 

\begin{lemma}
\label{lemma:exponential-kernel}
Suppose that $(X_1, Y_1)$ and $(X_2, Y_2)$ are independent copies of $(X, Y)$ where $X \sim P_X$ for some $P_X$ satisfying \Cref{assp:covariate-dist} and
\begin{equation*}
Y\mid X=x\sim\operatorname{Exp} \left(\exp\{\gamma+x^\top\beta\} \right),
\end{equation*}
where $\operatorname{Exp}(\lambda)$ denotes the exponential distribution with rate $\lambda$. Recall that 
\begin{equation*}
    h_{\mathrm{Exp}}(y_1,y_2) :=
    \ind\{\max(y_1,y_2)>2\min(y_1,y_2)\}-\frac{2}{3}.
\end{equation*}
Then, for all $\gamma \in \mathbb{R}$, we have that 
\begin{equation}
    \E_{0,\gamma} \left[h_{\mathrm{Exp}}(Y_1,Y_2)\right]=0.
    \label{eq:exponential-kernel-null}
\end{equation}
and 
\begin{equation}
    \E_{\beta,\gamma} \left[h_{\mathrm{Exp}}(Y_1,Y_2)\right] \gtrsim  \min\{\|\beta\|_2^2,1\}.
    \label{eq:exponential-kernel-alternative}
\end{equation}
\end{lemma}

\begin{proof}
Write $\lambda_i:=\exp\{\gamma+X_i^\top\beta\}$ for $i=1,2$.
Conditional on $(X_1,X_2)$, the responses are independent, and
\begin{equation*}
    \Pr(Y_1>2Y_2\mid X_1,X_2) =\int_0^\infty \exp(-2\lambda_1 y)\lambda_2 \exp(-\lambda_2 y) \dint y
    =\frac{\lambda_2}{2\lambda_1+\lambda_2}.
\end{equation*}
Since the events $\{Y_1>2Y_2\}$ and $\{Y_2>2Y_1\}$ are disjoint, setting $u:=\beta^\top(X_1-X_2)$, we obtain
\begin{align*}
    \E[h_{\mathrm{Exp}}(Y_1,Y_2)\mid X_1,X_2]
    &=\frac{1}{1+2e^u}+\frac{e^u}{e^u+2}-\frac23\\
    &=\frac{2(e^u-1)^2}{3(1+2e^u)(e^u+2)} = \frac{2 e^u (e^u + e^{-u} -2)}{3 e^u(2e^u + 5 + 2 e^{-u})}\\
    &=\frac{4(\cosh(u)-1)}{3(4\cosh(u)+5)} =:f(u).
\end{align*}
Under the null, $u=0$ almost surely and $f(0)=0$, proving \eqref{eq:exponential-kernel-null}. 

Since the map $x\mapsto4x/(27+12x)$ is increasing on $[0,\infty)$ and $\cosh(u)-1\geq u^2/2$, we have
\begin{equation*}
    f(u) =\frac{4(\cosh(u)-1)}{27+12(\cosh(u)-1)}
    \geq\frac{2u^2}{27+6u^2}.
\end{equation*}
By the Cauchy--Schwarz inequality, we have that 
\begin{equation*}
    \left(\E[u^2]\right)^2 \leq
    \E\left[\frac{u^2}{27+6u^2}\right]
    \E\left[u^2(27+6u^2)\right].
\end{equation*}
For $\beta$ such that $\|\beta\|_2>0$, \Cref{assp:covariate-dist} gives $2\underline{\lambda}\|\beta\|_2^2 \leq \E[u^2]  \leq 2\overline{\lambda}\|\beta\|_2^2$ and $\E[u^4]\lesssim C_X^4\|\beta\|_2^4$. We can therefore lower bound 
\begin{align*}
    \E_{\beta,\gamma}[h_{\mathrm{Exp}}(Y_1,Y_2)]
    =\E[f(u)] &\geq \frac{2\bigl(\E[u^2]\bigr)^2}{27\E[u^2]+6\E[u^4]}\\
    &\gtrsim \frac{\|\beta\|_2^2}{1+\|\beta\|_2^2}
    \gtrsim \min\{\|\beta\|_2^2,1\},
\end{align*}
which proves \eqref{eq:exponential-kernel-alternative}.
\end{proof}

\subsection{Proof of the upper bounds in \Cref{thm:main-result}}
\label{app:upper_assembly}
\begin{proof}
    Run the multiscale mean test and sparse covariance test, together with the response-only test for GLMs with more than two response support points, using independent randomisation, privacy budget $\epsilon/3$ and error bound $\eta = 1/12$ for each component. Disregard any component whose validity conditions fail, and reject whenever any remaining component rejects. By \Cref{lemma:composition}, the resulting test is $(\epsilon, 0)$-DP. Its type-I error is at most $1/4$, while its type-II error is at most $1/12$ whenever at least one remaining component's separation condition holds.

    For the mean test, choose
\begin{equation*}
\tilde r =c  \ell_{n,d,\epsilon}^3\left\{\frac{d^{1/4}}{\sqrt n}+\frac{d^{1/4}}{n^{3/4}\epsilon^{1/2}}+\frac1{n\epsilon}\right\}
\end{equation*}
with sufficiently large fixed constant $c>0$, and run it when $\tilde r\leq\min\{c_\mu,\sqrt{\overline\lambda}\}$. Since $\tilde r\gtrsim n^{-1/2}$, the dyadic grid has $J\lesssim\ell_{n,d,\epsilon}$ scales. By \eqref{eq:IT-test_s-radius}, we have  
\begin{equation*}
\rho\left(n,d,\frac{\epsilon}{3J},\frac1{12J}\right) \lesssim\frac{\ell_{n,d,\epsilon}d^{1/4}}{\sqrt n} +\frac{\ell_{n,d,\epsilon}^{3/2}d^{1/4}}{n^{3/4}\epsilon^{1/2}}+\frac{\ell_{n,d,\epsilon}^3}{n\epsilon} \leq\tilde r.
\end{equation*}
By \Cref{lemma:private-multiscale-reduction,lemma:transformation_properties}, this component detects alternatives with $\|\beta\|_2\geq\tilde r/c_\mu$. By \Cref{prop:sparse-upper}, the test for the sparse regime detects alternatives at separation at most $C\{\ell_{n,d,\epsilon}^{1/2}\sqrt{s^*/n}+(n\epsilon)^{-1}\}$.
Taking the minimum of these component bounds gives a sufficient separation of at most $C\ell_{n,d,\epsilon}^3\varrho_{s^*,d}$ for linear regression and GLMs with two response support points, and for larger-support GLMs when $d\leq n$.

For larger-support GLMs, when $d>n$ \Cref{prop:generic-radial-test} gives the separation $C\{n^{-1/4}+(n\epsilon)^{-1/2}\}$ for $\|\beta\|_2\leq r_0$. Combining the response-only and sparse bounds therefore gives a sufficient separation of at most $C\ell_{n,d,\epsilon}^3\varrho_{s^*,n}$. Taking $c_P$ sufficiently small ensures that a component attaining the relevant bound satisfies its validity conditions and that the local separation lies below $r_0$. The combined test thus has total error at most $1/3$ at the claimed separation, and its construction does not depend on $s^*$.
\end{proof}

\section{Technical details on the polynomial-time implementations}
\label{app:polytime_implementation}
\Cref{app:hajek-mechanism} contains all the technical details on the privacy guarantees of the subroutine in \Cref{alg:polytime_test}, and its error guarantees are proven in \Cref{app:polytime-test}. The proof of the polynomial-time sparse test from \Cref{prop:sparse-polytime} is in \Cref{app:sparse-polytime}. 

\subsection{Local H\'ajek mean algorithm with explicit implementation}
\label{app:hajek-mechanism}
\Cref{alg:local-hajek-primitive} is based on Algorithm~1 of \citet{chaudhuri2024DPUstats}, specialised to U-statistics of order two and with the additional wrapper step omitted. We use their local H\'ajek stabilisation construction, while deriving below the explicit sensitivity bound and privacy calibration used in our implementation.

\begin{algorithm}[h]
\caption{LocalH\'ajekMean}
\label{alg:local-hajek-primitive}
\begin{algorithmic}[1]
\Require Data $z_{1:n}$ with $n\ge2$; symmetric kernel $h$ with $\sup h-\inf h\le C_h$, $C_h>0$; tolerance $\xi\ge0$; privacy budget $\epsilon_0>0$.
\State Set $\epsilon'\gets \epsilon_0/6$ and compute
\begin{equation*}
 U_h\gets\binom n2^{-1}\sum_{i<j}h(z_i,z_j),
 \qquad \widehat h_i\gets\frac1{n-1}\sum_{j\ne i}h(z_i,z_j),\quad i\in[n].
\end{equation*}
\State Set
\begin{equation*}
 L\gets\min\left\{\ell\in[n]:\#\left\{i\in[n]:|\widehat h_i-U_h|>\xi+\frac{12C_h\ell}{n}\right\}\le\ell\right\}.
\end{equation*}
\State For $i\in[n]$, set
\begin{equation*}
 w_i\gets\left[1-\frac{n\epsilon'}{12C_h}\left(|\widehat h_i-U_h|-\xi-\frac{12C_hL}{n}\right)_+\right]_+.
\end{equation*}
\State Compute
\begin{equation}
 A\gets U_h+\binom {n}{2}^{-1}\sum_{i<j}\min\{w_i,w_j\}\{h(z_i,z_j)-U_h\}.
 \label{eq:stabilised-Ustat}
\end{equation}
\State Define
\begin{equation}
 G(v):= \frac{4\xi}{n}(1+\epsilon'v)+\frac{57 C_h}{n^2}\left\{v+\left(\epsilon'+\frac{1}{n}\right)v^2+\frac1{\epsilon'}\right\},\qquad v\ge1.
 \label{eq:local-sensitivity-upper-bound}
\end{equation}
\State Set $S\gets\max_{u\in\{0,\ldots,\lceil2/\epsilon'\rceil\}}\exp(-\epsilon'u)G(L+u)$, evaluated as in \Cref{remark:maximisation}
\State Draw $Q$ independently with density $q(t) := \sqrt{2}/\{\pi(1+t^4)\}$.
\State \Return $A+(S/\epsilon')Q$.
\end{algorithmic}
\end{algorithm}

\begin{remark}
\label{remark:maximisation}
The maximisation in Step 6 of \Cref{alg:local-hajek-primitive} does not require checking $O(1/\epsilon')$ values. Writing $G(v)=a_0+a_1v+a_2v^2$, the objective is maximised at either 0, $\lceil2/\epsilon' \rceil$, or at one of the integers immediately below and above the stationary points satisfying
\begin{equation*}
 a_1+2a_2(L+u)-\epsilon'\{a_0+a_1(L+u)+a_2(L+u)^2\}=0.
\end{equation*}
\end{remark}

The following lemma gives the privacy and accuracy guarantees of \Cref{alg:local-hajek-primitive}, following the local H\'ajek analysis of \cite{chaudhuri2024DPUstats}. We specialise their sensitivity argument to order-two U-statistics and track the constants explicitly to obtain the bound in \Cref{lemma:local_sensitivity} and the calibration $\epsilon'=\epsilon_0/6$ used here.

\begin{lemma}\label{lemma:local_hajek}
Algorithm~\ref{alg:local-hajek-primitive} is $(\epsilon_0,0)$-DP for every $\epsilon_0>0$. For i.i.d.\ observations, write $\theta_h:=\E[h(Z_1,Z_2)]$. If $\alpha\in(0,1)$ and
\begin{equation}\label{eq:hajek_good}
\Prob\left\{\max_i|\widehat h_i-U_h|\leq\xi\right\}\geq1-\alpha,
\end{equation}
then its output $\widetilde U_h$ satisfies
\begin{equation}\label{eq:hajek_accuracy}
\Prob\left\{|\widetilde U_h-\theta_h|>C\left[\sqrt{\frac{\var(U_h)}{\alpha}}+\alpha^{-1/3}\left\{\frac{\xi}{n}(1+\epsilon_0^{-1})+\frac{C_h}{n^2}(1+\epsilon_0^{-2})+\frac{C_h}{n^3}(1+\epsilon_0^{-3})\right\}\right]\right\}\leq3\alpha,
\end{equation}
where $C>0$ is an absolute constant. 
\end{lemma}

\begin{proof} 
\medskip \noindent
    \underline{Privacy.} Let
    \begin{equation*}
    \tilde{S}(D):= \max_{u \in \mathbb{Z}_{\geq 0}} \exp(-\epsilon' u) G(L(D)+u)
\end{equation*}
    By \Cref{lemma:local_sensitivity}, we have that $LS_A(D) \leq G(L(D))$, and $|L(D) - L(D') | \leq 1$ for $D, D'$ that differ in one record. Since $G$ is increasing, we have 
    \begin{align*}
        \tilde{S}(D) &\leq\max_{u\ge0}\exp(-\epsilon'u)G(L(D')+u+1)\\
        &=\exp(\epsilon')\max_{v\ge1}\exp(-\epsilon'v)G(L(D')+v)\\
        &\leq\exp(\epsilon')\tilde{S}(D').
    \end{align*}
This implies that $\tilde{S}(D)$ is an $\epsilon'$-smooth upper bound on $\mathrm{LS}_A(D)$ in the sense of \cite{nissim2007smooth}. 

Since $vG'(v)\leq 2G(v)$ for $v>0$ and $L(D)\geq0$, we have
\begin{equation*}
\frac{\dint}{\dint u}\left[e^{-\epsilon'u}G(L(D)+u)\right] \leq e^{-\epsilon'u}G(L(D)+u) \left\{\frac{2}{L(D)+u}-\epsilon'\right\} \leq0,
\qquad u\geq\frac{2}{\epsilon'},
\end{equation*}
so the expression for $S(D)$ in Step 6 of \Cref{alg:local-hajek-primitive} is equal to $\tilde{S}(D)$.

The next part of the proof follows the arguments of Lemma 2.7 of \cite{nissim2007smooth}, with $q$ being the probability density in Step 7 of \Cref{alg:local-hajek-primitive}. Writing $\sigma(D) := S(D)/\epsilon'$, let 
\begin{equation*}
    p_D(x) := \frac{1}{\sigma(D)} q \left( \frac{x-A(D)}{\sigma(D)}\right), \quad \tilde{p}(x) := \frac{1}{\sigma(D')}q \left(\frac{x-A(D)}{\sigma(D')} \right),
\end{equation*}
so that $p_D(x)$ is the output density and $\tilde{p}(x)$ an intermediate density introduced for analysis. We can then bound
\begin{equation}
    \left| \log \frac{p_D(x)}{p_{D'}(x)} \right| \leq \left| \log \frac{p_D(x)}{\tilde{p}(x)} \right| + \left| \log \frac{\tilde{p}(x)}{p_{D'}(x)} \right|. 
    \label{eq:nissim_decomp}
\end{equation}
We bound the first term of \eqref{eq:nissim_decomp}, corresponding to the scale change. Considering the change of variable $t = (x-A(D))/\sigma$, we have that 
\begin{align*}
\left|\frac{\partial}{\partial\log\sigma}\log\left\{\frac1\sigma q\left(\frac{x-A(D)}\sigma\right)\right\} \right|&= |-1 - t(\log q)'(t)|\\
&= \left| \frac{1 - 3t^4}{1 + t^4} \right| \leq 3.
\end{align*}
Thus, by the mean value theorem, we can bound the first term of \eqref{eq:nissim_decomp} by 
\begin{equation*}
\left| \log \frac{p_D(x)}{\tilde{p}(x)}\right| \leq 3 | \log( \sigma(D)) - \log(\sigma(D'))| \leq 3 \epsilon'.
\end{equation*}

Since we have that 
\begin{align*}
     \left|\frac{\partial}{\partial\mu}\log\left\{\frac1{\sigma(D')}q\left(\frac{x-\mu}{\sigma(D')}\right)\right\} \right|
&= \left|-\frac1{\sigma(D')}(\log q)'\left(\frac{x-\mu}{\sigma(D')}\right) \right|, 
\end{align*}
we can similarly use 
\begin{equation*}
    |(\log q)'(t)| \leq  \left| \frac{-4t^3}{1+t^4} \right| \leq 3
\end{equation*}
and the mean value theorem to bound 
\begin{equation*}
\left|\log\frac{\widetilde p(x)}{p_{D'}(x)}\right|
\le\frac{3|A(D)-A(D')|}{\sigma(D')}
\le\frac{3S(D')}{S(D')/\epsilon'}=3\epsilon'.
\end{equation*}

Substituting these bounds into \eqref{eq:nissim_decomp} gives $| \log \{p_D(x)/p_{D'}(x)\}| \leq 6 \epsilon' = \epsilon_0$. 

\medskip\noindent\underline{Accuracy.} For the non-private component, we have by Chebyshev's inequality that 
\begin{equation}
    \Prob \left( |U_h - \theta_h| > \sqrt{\mathrm{Var}(U_h)/\alpha} \right) \leq \alpha.
    \label{eq:hajek_non-priv-Chebyshev}
\end{equation}

On the event in \eqref{eq:hajek_good}, we have that $L(D)=1$ and $w_i=1$ for all $i \in [n]$, so that $A(D) = U_h$.  
Substituting $L(D)=1$ into the definition of $S(D)$ gives 
\begin{equation*}
    S(D)\leq \sup_{u \geq0}\exp(-\epsilon'u)G(1+u).
\end{equation*}
Using
\begin{align*}
\sup_{u\ge0}\exp(-\epsilon'u)(1+u)^k
&\leq2^{k-1}\left\{1+\sup_{u\ge0}u^k\exp(-\epsilon'u)\right\}\\
&=2^{k-1}\left\{1+\left(\frac{k}{e\epsilon'}\right)^k\right\}
\lesssim1+\frac1{(\epsilon')^k},\qquad k=1,2,
\end{align*}
we obtain
\begin{align*}
\frac{S(D)}{\epsilon'}&\lesssim\frac{\xi}{n}\left(1+\frac1{\epsilon'}\right)
+\frac{C_h}{n^2}\left(1+\frac1{\epsilon'}+\frac1{(\epsilon')^2}\right)
+\frac{C_h}{n^3}\left(\frac1{\epsilon'}+\frac1{(\epsilon')^2}+\frac1{(\epsilon')^3}\right)\\
&\lesssim\frac{\xi}{n}(1+\epsilon_0^{-1})+\frac{C_h}{n^2}(1+\epsilon_0^{-2})+\frac{C_h}{n^3}(1+\epsilon_0^{-3}).
\end{align*}
We also have, for $Q \sim q$, that 
\begin{equation}
    \Prob(|Q| > u ) = \frac{2 \sqrt{2}}{\pi} \int_u^\infty  \frac{1}{1 + t^4} \dint t \leq \frac{2\sqrt{2}}{3 \pi u^3}, \quad u >0.
    \label{eq:hajek_noise-tail}
\end{equation}
Combining \eqref{eq:hajek_non-priv-Chebyshev}, the event in \eqref{eq:hajek_good}, and \eqref{eq:hajek_noise-tail} with $u=(2\sqrt2/(3\pi\alpha))^{1/3}$ in a union bound argument proves \eqref{eq:hajek_accuracy}. 
\end{proof}

\begin{lemma}
   For every $\epsilon' >0$, and neighbouring datasets $D, D'$, we have that $|L(D) - L(D')| \leq 1$, and $\mathrm{LS}_A(D) \leq G(L(D))$, with $G$ as in \eqref{eq:local-sensitivity-upper-bound}. 
   \label{lemma:local_sensitivity}
\end{lemma}
\begin{proof}
    Suppose $D, D'$ differ at index $k$; we will use primes to denote any quantity that is computed on $D'$. We can bound 
    \begin{align}
        |U_h - U_h'| = \left|\frac{2}{n(n-1)} \sum_{j\ne k}\{h(z_j,z_k)-h(z_j,z_k')\} \right| \leq \frac{2C_h}{n},
        \label{eq:local_sens_udiff}
    \end{align}
    and for $i \neq k$, we have 
    \begin{align}
        \left|\hat{h}_i - U_h - (\hat{h}_i'- U_h') \right| &= \left| \frac{n-2}{n(n-1)} \{h(z_i,z_k)-h(z_i,z_k')\} -\frac{2}{n(n-1)} \sum_{j\ne i,k}\{h(z_j,z_k)-h(z_j,z_k')\}\right| \notag\\
        &\leq \frac{3C_h(n-2)}{n(n-1)}.
        \label{eq:LS-proj_diff}
    \end{align}
    Therefore, $|\hat{h}_i - U_h| \leq \xi + 12C_hL(D)/n$ implies that 
    \begin{equation*}
        |\hat{h}_i' - U_h'| \leq \xi + \frac{12C_h(L(D)+1)}{n}.
    \end{equation*}
    There are at most $L(D)$ indices of $D$ above the threshold and one changed index. Hence $L(D')\leq L(D)+1$, and interchanging the roles of $D, D'$ proves $|L(D)-L(D')|\leq1$.

    For $T \geq 1$, define the clipping function
    \begin{equation*}
        F_T(x):= x\min\left\{ 1,\frac{n\epsilon'}{12C_h} \left[|x|-\xi-\frac{12C_h T}{n}\right]_+ \right\}
    \end{equation*}
    so that $F_{L(D)}(\hat{h}_i-U_h) = (\hat{h}_i - U_h)(1-w_i)$. We also have the following bounds: 
    \begin{align}
 |x-F_T(x)|&\le\xi+\frac{12C_hT}{n}+\frac{3C_h}{n\epsilon'}(1-\epsilon'T)_+^2,\label{eq:hajek-retained-bound}\\
 |F_T(x)-F_T(y)|&\le\frac{3C_h}{n}\left\{2+\epsilon'\left(\frac{n\xi}{12C_h}+T\right)\right\},\qquad |x-y|\le\frac{3C_h}{n},\label{eq:hajek-fixed-bound}\\
 |F_T(x)-F_{T+1}(y)|&\le\frac{12C_h}{n}\left[\frac{5\epsilon'}4\left(\frac{n\xi}{12C_h}+T-\frac{1}{4}\right)+\frac{3}{2}\right],\qquad |x-y|\le\frac{3C_h}{n}.\label{eq:hajek-shifted-bound}
\end{align}
We defer showing these bounds to the end of this proof. In particular, \eqref{eq:hajek-retained-bound} implies that 
\begin{align}
     \left|(\hat{h}_i - U_h)w_i \right| \leq \xi+\frac{C_h}{n} \left\{ 6L(D)+ \frac{3}{\epsilon'} + 3 \epsilon' L(D)^2 \right\}.
\label{eq:hajek-retained-projection}
\end{align}

    We now use these bounds to control the local sensitivity of $A$. We can rewrite \eqref{eq:stabilised-Ustat} as 
    \begin{align*}
    A &=U_h-\binom{n}{2}^{-1}\sum_{i<j}\{(1-w_i)+(1-w_j)\}\{h(z_i,z_j)-U_h\}+R\\
    &=U_h-\binom{n}{2}^{-1}\sum_i(1-w_i)\sum_{j\ne i}\{h(z_i,z_j)-U_h\}+R\\
    &=U_h-\frac{2}{n}\sum_i(1-w_i)(\hat h_i-U_h)+R
    =U_h-\frac{2}{n}\sum_iF_{L(D)}(\hat h_i-U_h)+R,
    \end{align*}
    where we denote 
    \begin{equation*}
        R:=\binom {n}{2}^{-1}\sum_{i<j}\min\{1-w_i,1-w_j\}\{h(z_i,z_j)-U_h\}.
    \end{equation*}
    
    The first equality follows from
\begin{equation*}
    \min\{w_i,w_j\}=1-(1-w_i)-(1-w_j)+\min\{1-w_i,1-w_j\}
\end{equation*}
and that $\sum_{i<j}\{h(z_i,z_j)-U_h\}=0$. 
The third uses $\sum_{j\ne i}\{h(z_i,z_j)-U_h\}=(n-1)(\hat h_i-U_h)$.

With this form of $A$, and separating out the changed index $k$, we have 
\begin{align}
A-A' &=U_h-U_h'-\frac{2}{n}\{(\hat h_k-U_h)-(\hat h_k'-U_h')\} +\frac{2}{n}\{(\hat h_k-U_h)w_k-(\hat h_k'-U_h')w_k'\} \notag\\ 
&\quad-\frac{2}{n}\sum_{i\ne k}\{F_{L(D)}(\hat h_i-U_h)-F_{L(D')}(\hat h_i'-U_h')\}+R-R' \notag\\ 
&=\frac{2}{n}\bigg[U_h-U_h'+(\hat h_k-U_h)w_k-(\hat h_k'-U_h')w_k' \notag\\
&\quad \hspace{2cm}-\sum_{i\ne k}\{F_{L(D)}(\hat h_i-U_h)-F_{L(D')}(\hat h_i'-U_h')\}\bigg]+R-R'.
\label{eq:hajek-cancellation}
\end{align}
The second equality is from 
noting that 
\begin{align*}
    (\hat{h}_k - \hat{h}_k')-  (U_h - U_h')  &= \frac1{n-1}\sum_{j\ne k}\{h(z_k,z_j)-h(z_k',z_j)\} -  (U_h - U_h')\\
    &=\frac{n-2}{2}(U_h-U_h').
\end{align*}

Since $|L(D) - L(D')| \leq 1$, we will now bound \eqref{eq:hajek-cancellation} in the two different cases. 

\medskip 
\noindent \underline{Case 1: $L(D) = L(D')$.}
 In this case, using \eqref{eq:local_sens_udiff} and \eqref{eq:hajek-retained-projection}, we have
\begin{align*}
    |A - A'| &\leq \frac{2}{n} \left\{ \frac{2C_h}{n} +2 \left[\xi+\frac{C_h}{n} \left\{ 6 L(D)+ \frac{3}{\epsilon'} + 3 \epsilon' L(D)^2 \right\}\right]\right\}\\
    & \quad + \frac{2}{n}\left|\sum_{i\neq k} \{F_{L(D)}(\hat{h}_i - U_h)-F_{L(D')}(\hat{h}_i' - U_h')\}\right|+ |R-R'|.
\end{align*}
In the summation, at most $2L(D)$ terms are non-zero, and using \eqref{eq:hajek-fixed-bound} we can bound 
\begin{equation*}
 |F_{L(D)}(\hat{h}_i - U_h)-F_{L(D)}(\hat{h}_i' - U_h')| \leq \frac{3C_h}{n} \left\{2+\epsilon'\left(\frac{n\xi}{12C_h}+L(D)\right)\right\}.
\end{equation*}

Let $B:=\{i:1-w_i>0\}$ and similarly for $B'$, then $|B|$, $ |B'| \leq L(D)$. For $i$ such that $i \neq k$, we have from \eqref{eq:LS-proj_diff} that 
\begin{equation*}
    |w_i - w_i'| \leq \frac{n\epsilon'}{12C_h} \left| \hat{h}_i - U_h - (\hat{h}_i ' -U_h') \right| \leq \frac{\epsilon'}{4}.
\end{equation*}
 We can therefore write 
\begin{align*}
|R-R'| &\leq \binom{n}{2}^{-1}\Bigg\{ \sum_{\substack{i<j,\ i,j\ne k\\i,j\in B'}} \min\{1-w_i',1-w_j'\}|U_h'-U_h|\\
&\hspace{2cm}+ \sum_{\substack{i<j,\ i,j\ne k\\i,j\in B\cup B'}}
\left|\min\{1-w_i,1-w_j\}-\min\{1-w_i',1-w_j'\}\right| |h(z_i,z_j)-U_h|\\
& +\sum_{\substack{j\ne k\\j\in B\cup B'}}
\left|\min\{1-w_k,1-w_j\}\{h(z_k,z_j)-U_h\}
-\min\{1-w_k',1-w_j'\}\{h(z_k',z_j)-U_h'\}\right| \Bigg\}\\
&\leq \binom{n}{2}^{-1}
\left\{\frac{C_h \epsilon'}{4}\binom{2L(D)}{2}
+\frac{2C_h}{n}\binom{L(D)}{2}+4L(D)C_h\right\}\\
&\leq \frac{C_h}{n^2}\left\{16L(D)+2\epsilon'L(D)^2-\epsilon'L(D)+\frac{4(L(D)^2-L(D))}{n}\right\}.
\end{align*}

Hence, we can bound 
\begin{equation*}
|A-A'|\le\frac{\xi}{n}(4+\epsilon'L(D))+\frac{C_h}{n^2}\left\{64L(D)+4+26\epsilon'L(D)^2-\epsilon'L(D)+\frac{4(L(D)^2-L(D))}{n}+\frac{12}{\epsilon'}\right\}.
\end{equation*}

\medskip 
\noindent \underline{Case 2: $|L(D') -L(D)|=1$.} Without loss of generality, assume that $L(D') = L(D) + 1$. In this case, using \eqref{eq:local_sens_udiff} and \eqref{eq:hajek-retained-projection}, we have
\begin{align*}
|A-A'|&\le\frac{2}{n}\left\{\frac{2C_h}{n}+2\xi+\frac{C_h}{n}\left[6L(D)+6(L(D)+1)+\frac6{\epsilon'}+3\epsilon'\{L(D)^2+(L(D)+1)^2\}\right]\right\}\\
&\quad+\frac{2}{n}\left|\sum_{i\neq k} \left\{F_{L(D)}(\hat h_i-U_h)-F_{L(D)+1}(\hat h_i'-U_h') \right\}\right|+|R-R'|.
\end{align*}
For $i\neq k$, $i \in B'$ requires that $|\hat h_i'-U_h'|>\xi+12C_h(L(D)+1)/n$, which, by \eqref{eq:LS-proj_diff}, implies $|\hat h_i-U_h|>\xi+12C_hL(D)/n$. We therefore have that $B'\setminus\{k\}\subseteq B\setminus\{k\}$, so at most $L(D)$ terms in the sum are non-zero. For each such term, \eqref{eq:hajek-shifted-bound} gives 
\begin{equation*}
|F_{L(D)}(\hat h_i-U_h)-F_{L(D)+1}(\hat h_i'-U_h')|
\le\frac{12C_h}{n}\left\{\frac{5\epsilon'}4\left(\frac{n\xi}{12C_h}+L(D)-\frac{1}{4}\right)+\frac{3}{2}\right\}.
\end{equation*}

Further, we have 
\begin{align*}
|w_i-w_i'|
&\leq\frac{n\epsilon'}{12C_h}\left|\left(|\hat h_i-U_h|-\xi-\frac{12C_hL(D)}n\right)-\left(|\hat h_i'-U_h'|-\xi-\frac{12C_hL(D')}n\right)\right| \\
&\leq\frac{n\epsilon'}{12C_h}\left(\frac{3C_h}{n}+\frac{12C_h}{n}\right)
=\frac{5\epsilon'}4,\qquad i\ne k.
\end{align*}
Using the same grouping argument as in Case 1 leads to  
\begin{align*}
|R-R'| &\leq\frac4{n^2}\left\{\frac{2C_h}{n}\binom{L(D)}2+\frac{5\epsilon'C_h}{4}\binom{L(D)}2+2C_hL(D)\right\}\\
&=\frac{C_h}{n^2}\left\{8L(D)+\frac{5\epsilon'}2(L(D)^2-L(D))+\frac{4(L(D)^2-L(D))}{n}\right\}.
\end{align*}

We therefore have that
\begin{align*}
|A-A'|&\leq\frac{\xi}{n}\left(4+\frac52\epsilon'L(D)\right)\\
&\quad+\frac{C_h}{n^2}\left\{68L(D)+16+\frac{89}{2}\epsilon'L(D)^2+2\epsilon'L(D)+6\epsilon'+\frac{4(L(D)^2-L(D))}{n}+\frac{12}{\epsilon'}\right\}.
\end{align*}
If $L(D')=L(D)-1$, interchanging $D$ and $D'$ gives $|A(D)-A(D')|\leq G(L(D'))\leq G(L(D))$, since $G$ is increasing, so the bound holds for every neighbouring pair $D, D'$.

\medskip
\noindent \underline{Bound for both cases.}
Using $L(D)\ge1$, we have that $16\le16L(D)$ and $2\epsilon'L(D)+6\epsilon'\le8\epsilon'L(D)^2$. Thus, in both cases, we can bound 
\begin{align*}
|A-A'|\leq\frac{4\xi}{n}(1+\epsilon'L(D))+\frac{C_h}{n^2}\left\{84L(D)+\frac{105}{2}\epsilon'L(D)^2+\frac{4L(D)^2}{n}+\frac{12}{\epsilon'}\right\} \\
\leq\frac{4\xi}{n}(1+\epsilon'L(D))+\frac{C_h}{n^2}\left\{57L(D)+57\epsilon'L(D)^2+\frac{4L(D)^2}{n}+\frac{105}{2\epsilon'}\right\}
\end{align*}
where the second inequality is from noting  that $27L(D)\leq(9/2)\epsilon'L(D)^2+81/(2\epsilon')$.

This shows that $|A-A'|\leq G(L(D))$. Taking the supremum over neighbouring $D'$ yields $\mathrm{LS}_A(D)\leq G(L(D))$.

\medskip
\noindent \underline{Showing the bounds in \eqref{eq:hajek-retained-bound}, \eqref{eq:hajek-fixed-bound}, and \eqref{eq:hajek-shifted-bound}.}
To show \eqref{eq:hajek-retained-bound}, we first consider the interval where the clipping factor is strictly between 0 and 1, i.e. 
\begin{equation}
   \xi+\frac{12C_hT}{n} < |x| < \xi+\frac{12C_hT}{n}+\frac{12C_h}{n\epsilon'}
    \label{eq:transition_interval}
\end{equation}
When $x$ satisfies \eqref{eq:transition_interval}, writing $s=|x|-\xi - 12C_hT/n$ we have that 
\begin{align*}
 |x-F_T(x)|
 &=\xi+\frac{12C_hT}{n}+s\left(1-\epsilon'T-\frac{n\epsilon'\xi}{12C_h}-\frac{n\epsilon's}{12C_h}\right), \qquad 0 \leq s \leq \frac{12C_h}{n \epsilon'}\\
 &\le\xi+\frac{12C_hT}{n}+\frac{3C_h}{n\epsilon'}(1-\epsilon'T)_+^2.
\end{align*}
The last inequality uses $\xi\ge0$ and maximises the remaining quadratic in $s\ge0$.
When $x$ does not satisfy \eqref{eq:transition_interval}, we have that $|x-F_T(x)|\le\xi+12C_hT/n$. 

Next, to show \eqref{eq:hajek-fixed-bound}, note that on \eqref{eq:transition_interval}, we have 
\begin{equation}
    F_T'(x) = \frac{n\epsilon'}{12C_h} \left( 2|x| - \xi - \frac{12C_hT}{n}\right) \leq 2 + \epsilon' \left( \frac{n \xi}{12C_h} + T \right).
    \label{eq:hajek_fixed_lipschitz}
\end{equation}
Outside of \eqref{eq:transition_interval}, the derivative of $F_T(x)$ is either 0 or 1. Since $F_T$ is continuous and piecewise differentiable, the right-hand side of \eqref{eq:hajek_fixed_lipschitz} is a global Lipschitz bound on $F_T$, which implies \eqref{eq:hajek-fixed-bound}. 

Finally, we show \eqref{eq:hajek-shifted-bound}. Since $F_T$ and $F_{T+1}$ are odd and non-decreasing, and $0\leq F_{T+1}(r)\le F_T(r)$ for $r\geq0$, we have for $x, y$ of the same sign that 
\begin{equation*}
    |F_T(x) - F_{T+1}(y)| \leq \sup_{r \geq 0} \left\{ F_T \left(r+\frac{3C_h}{n} \right)-F_{T+1}(r)\right\}, \quad |x-y| \leq \frac{3C_h}{n}.
\end{equation*}
Substituting $u := nr/(12C_h)$, the expression in the supremum can be written as 
\begin{equation}
    \frac{12C_h}{n} \left[\left(u+\frac{1}{4} \right)\min\{1,\epsilon'(u+1/4-b)_+\} -u\min\{1,\epsilon'(u-b-1)_+\} \right],
\label{eq:LS-quadratic}
\end{equation}
where $b := n\xi/(12C_h) + T.$

If $\epsilon' \leq 4/5$, then \eqref{eq:LS-quadratic} is maximised at $u=b + (\epsilon')^{-1} - 1/4$, where the term in the square brackets is $3/2 + 5 \epsilon'/4(b-1/4)$. If $\epsilon'\geq 4/5$, then \eqref{eq:LS-quadratic} is maximised at $u=b+1$, which is where the second clipped term starts increasing. At this point, the term in the square brackets equals $b+5/4$. In this case, we can also bound
\begin{equation*}
b+\frac{5}{4}\leq \frac{3}{2}+\frac{5\epsilon'}{4}\left(b-\frac{1}{4}\right)
\end{equation*}
since the difference is $(5\epsilon'/4-1)(b-1/4)\ge0$.

It remains to consider when $x, y$ have different signs. In this case, we can bound 
\begin{equation*}
    |F_T(x) - F_{T+1}(y)| \leq |x| + |y| = |x-y| \leq \frac{3C_h}{n}. 
\end{equation*}

This concludes the proof. 
\end{proof}

\subsection{Proof of \Cref{prop:polytime_test}}
\label{app:polytime-test}
\begin{proof}
The test at each scale uses privacy budget $\epsilon/J$ on the full sample. By \Cref{lemma:local_hajek}, sequential composition and post-processing give $(\epsilon,0)$-DP.  Each scale requires $O(dn^2)$ arithmetic operations, so the total computational cost is $O(Jdn^2)$.

Write  $\theta_{s_j}:=\E[h_{s_j}(Z_1,Z_2)]$ and $\alpha:=\eta/(3J)$. Consider either $\E[Z]=0$, or let $s_j$ be a scale satisfying  $s_j\leq\|\E[Z]\|_2\leq2s_j$. By \Cref{lemma:clipped_kernel}, the event in \eqref{eq:hajek_good} holds with probability at least $1-\alpha$ for $\xi=\xi_{s_j}$. Applying \Cref{lemma:local_hajek} with $\epsilon_0=\epsilon/J$ and $C_h=2R(s_j)$ gives that
\begin{align*}
|\widetilde U_j-\theta_{s_j}|
&\lesssim \sqrt{\frac{s_j^2/n+d/n^2}{\alpha}}\\
&\hspace{1cm}+\alpha^{-1/3}\left\{\frac{l_{s_j}(s_j+\sqrt{d/n})}{n}\left(1+\frac J\epsilon\right)+\frac{R(s_j)}{n^2}\left(1+\frac{J^2}{\epsilon^2}\right)+\frac{R(s_j)}{n^3}\left(1+\frac{J^3}{\epsilon^3}\right)\right\}
\end{align*}
with probability at least $1-3\alpha=1-\eta/J$.

Under the condition in \eqref{eq:polytime-rate} with sufficiently large constant and polylogarithmic factors, we thus have with probability at least $1 - \eta/J$ that 
\begin{align*}
    |\tilde{U}_j - \theta_{s_j}| \lesssim\eta^{-1/2}\Bigg\{\ell_{n,d,\epsilon}^{4/3} \left(\frac{s_j}{\sqrt n}+\frac{\sqrt d}{n}\right) +\ell_{n,d,\epsilon}^{7/3} \left(\frac{s_j}{n\epsilon} +\frac{\sqrt d}{n^{3/2}\epsilon}\right)+\frac{\ell_{n,d,\epsilon}^{10/3}\sqrt d}{n^2\epsilon^2} +\frac{\ell_{n,d,\epsilon}^{13/3}\sqrt d}{n^3\epsilon^3} \Bigg\} \leq \frac{s_j^2}{4}.
\end{align*}

Under the null, these bounds apply at every scale. By \eqref{eq:clipping_bias} and a union bound, with probability at least $1-\eta$, we have 
\begin{equation*}
    \widetilde{U}_j\leq s_j^2/16+s_j^2/4 =5s_j^2/16<s_j^2/2 
\end{equation*}
for all $j$. Under this event, the test does not reject. 

Under an alternative with
$\tilde r\leq\|\E[Z]\|_2\leq\sqrt{\bar\lambda}$, the dyadic grid
contains an index $k$ such that
$s_k\leq\|\E[Z]\|_2\leq2s_k$.
With probability at least $1-\eta/J$,
\begin{equation*}
    \widetilde{U}_k
    \geq\|\E[Z]\|_2^2-s_k^2/16-s_k^2/4
    \geq11s_k^2/16>s_k^2/2.
\end{equation*}
The test therefore rejects with probability at least
$1-\eta/J\geq1-\eta$. Finally, let $c_\mu>0$ be the constant in \Cref{lemma:transformation_properties}. If $\tilde r\leq c_\mu$ and $C_1\geq1/c_\mu$, \eqref{eq:mean-to-regression} ensures that every alternative with $\|\beta\|_2\geq C_1\tilde r$ satisfies $\|\E[Z]\|_2\geq\tilde r$. The upper bound $\|\mathbb{E}[Z]\|_2\leq\sqrt{\overline\lambda}$ follows from $|\tanh|\leq1$ and $\E[XX^\top]\preceq\overline\lambda I$. Hence the mean testing guarantee for $\mathbb{E}[Z]$ transfers to the testing guarantee for $\beta$. 
\end{proof}

The following lemma bounds the clipping bias and the non-private variance of the proposed kernel, and provides the local stability required by the private U-statistic mechanism. 

\begin{lemma}
\label{lemma:clipped_kernel}
Suppose we have i.i.d.~$\{Z_i := X_i \tanh (Y_i)\}_{i=1}^n$ from a model in $\mathcal{P}_{\mathrm{Lin}}$ or  $\mathcal{P}_{\mathrm{GLM}}(\mathcal{F})$, with $\mathcal{F}$ satisfying \Cref{assp:glm-regularity}. Let $\alpha\in(0,1)$ and $0<s\leq\sqrt{\overline{\lambda}}$, and suppose either $\E[Z]=0$ or $s\leq\|\E[Z]\|_2\leq2s$. Define 
\begin{equation*}
    \hat{h}_i^{(s)} := \frac{1}{n-1} \sum_{j \neq i} h_s(Z_i, Z_j), \qquad U_s := \frac{1}{n}\sum_{i=1}^n \hat{h}_i^{(s)},
\end{equation*}

For $C_\mathrm{clip}$ sufficiently large, the kernel $h_s$ in \eqref{eq:clipped_kernel_def} satisfies
\begin{equation}
    \left|\mathbb{E}[h_s(Z_1,Z_2)]-\|\mathbb{E}[Z]\|_2^2\right|\leq \frac{s^2}{16}
    \label{eq:clipping_bias}
\end{equation}
and 
\begin{equation}  
    \operatorname{Var}\left(U_s(Z_{1:n})\right) \lesssim \frac{s^2}{n}+\frac{d}{n^2}.
    \label{eq:clipping_variance}
\end{equation}
Moreover, there exists a constant $C>0$ such that
\begin{equation}
    \Prob \left( \max_{i\in[n]} \left|\hat{h}_i^{(s)}-U_s(Z_{1:n})\right|\leq C l_s \left\{s+\sqrt{d/n}\right\} \right) \geq 1 - \alpha.
    \label{eq:hajek_weight_tol}
\end{equation}
\end{lemma}

\begin{proof}
    For $i \in [n]$, let $R_i := Z_i - \mathbb{E}[Z]$. We can write 
    \begin{equation*}
    Z_i^\top Z_j = \|\mathbb{E}[Z]\|_2^2 + \mathbb{E}[Z]^\top(R_i + R_j) + R_i^\top R_j.
    \end{equation*}
    By sub-Gaussianity of the second term and  \Cref{lemma:sub-Gaussian-inner-product}, we have that 
    \begin{equation}
        \Prob \left( |Z_i^\top Z_j| \geq C_{C_X, \overline{\lambda}} \{\sqrt{dt} + t\} \right) \leq e^{-t}
        \label{eq:inner_tail}
    \end{equation}
    for all $t\geq 1$. Let $u=C_{C_X,\overline\lambda}(\sqrt{dt}+t)$. Since $C_{C_X, \overline{\lambda}} \{\sqrt{d(2^kt)} + 2^kt\} \leq 2^k C_{C_X, \overline{\lambda}}\{\sqrt{dt} + t\} = 2^k u$, a dyadic
decomposition and \eqref{eq:inner_tail} give
\begin{align*}
\E \left[(Z_1^\top Z_2)^2\ind\{|Z_1^\top Z_2|>u\}\right]
&\leq \sum_{k=0}^\infty (2^{k+1}u)^2 \Prob\{|Z_1^\top Z_2|>2^ku\}\\
&\lesssim 24u^2\sum_{k=0}^\infty 4^k \exp(-2^kt)\\
&\lesssim (dt+t^2) \exp(-t).
\end{align*}

Since $l_s \geq 1 + \log \{(d+1)/s\}$, we can take $t=c_1l_s$ with sufficiently large absolute constant $c_1$ and then expanding $C_{\mathrm{clip}}$ to match in \eqref{eq:clipped_kernel_def} ensures that $R(s)\geq u$ and
\begin{equation}
    \E\left[\{h_s(Z_1,Z_2)-Z_1^\top Z_2\}^2\right]\leq
    \E\left[(Z_1^\top Z_2)^2\ind\{|Z_1^\top Z_2|>R(s)\}\right] \leq \frac{s^4}{256}.
    \label{eq:clipped_bias_squared}
\end{equation}

Since $\mathbb{E}[Z_1^\top Z_2] = \|\mathbb{E}[Z]\|_2^2$, we have 
\begin{align*}
    \left|\E[h_s(Z_1, Z_2)] - \|\mathbb{E}[Z]\|_2^2 \right| \leq \left\{\E \left[(h_s(Z_1,Z _2) - Z_1^\top Z_2)^2 \right] \right\}^{1/2} \leq \frac{s^2}{16}.
\end{align*}
This proves the clipping bias in \eqref{eq:clipping_bias}. 

To obtain the bound on the variance in  \eqref{eq:clipping_variance}, let 
\begin{equation*}
v_R(z):=\E[h_s(z,Z_2)],\qquad v(z):=\E[z^\top Z_2]=z^\top\mathbb{E}[Z].
\end{equation*}
Since $Z_2$ is independent of $Z_1$,
\begin{equation*}
v_R(Z_1)-v(Z_1)= \E[h_s(Z_1,Z_2)-Z_1^\top Z_2\mid Z_1],
\end{equation*}
so by a conditional Jensen's inequality and \eqref{eq:clipped_bias_squared}, we have 
\begin{align*}
\E \left[(v_R(Z_1)-v(Z_1))^2\right] &\leq \E\left[ \E\left[(h_s(Z_1,Z_2)-Z_1^\top Z_2)^2\mid Z_1\right]\right]\\
&=\E\left[(h_s(Z_1,Z_2)-Z_1^\top Z_2)^2\right] \leq \frac{s^4}{256}.
\end{align*}

Writing $\Omega:=\E[Z_1Z_1^\top]$, we have $\Omega\preceq\E[XX^\top]$ since
$Z_1=X_1\tanh(Y_1)$ and $|\tanh(Y_1)|\leq1$. It follows that 
\begin{equation*}
\var\{v(Z_1)\} \leq
\E[(Z_1^\top\mathbb{E}[Z])^2]
= \mathbb{E}[Z]^\top\Omega\mathbb{E}[Z]
\leq \overline\lambda\|\mathbb{E}[Z]\|_2^2
\leq 4\overline\lambda s^2.
\end{equation*}
Using $\var(A+B)\leq2\var(A)+2\E[B^2]$ and that
$s\leq\sqrt{\overline\lambda}$, we have
\begin{equation*}
\var\{v_R(Z_1)\} \leq 2\var\{v(Z_1)\} + 2\E[(v_R(Z_1)-v(Z_1))^2] \lesssim s^2.
\end{equation*}

Moreover, independence of $Z_1$ and $Z_2$ gives
\begin{align*}
\E[(Z_1^\top Z_2)^2]
&= \E \left[\tr(Z_1Z_1^\top Z_2Z_2^\top)\right]\\
&=\tr \left(\E[Z_1Z_1^\top]\E[Z_2Z_2^\top]\right)
\leq
\overline\lambda^2 d.
\end{align*}
Since $|h_s(Z_1, Z_2)| \leq |Z_1^\top Z_2|$, we also obtain $\var\{h_s(Z_1,Z_2)\}\lesssim d$. Finally, Hoeffding's variance decomposition for order-two U-statistics \citep[e.g.][Chapter 5]{serfling1980approximation} gives
\begin{equation*}
\var(U_{s})
\leq\frac{4(n-2)}{n(n-1)}\var \left(v_R(Z_1) \right)+
\frac{2}{n(n-1)}\var\left(h_s(Z_1,Z_2)\right)
\lesssim \left(\frac{s^2}{n}+\frac{d}{n^2}\right),
\end{equation*}
which proves \eqref{eq:clipping_variance}.

Finally, we prove \eqref{eq:hajek_weight_tol}. The concentration argument in the proof of the row sum in \Cref{lemma:typical_set_prob}, with $t=\log(6n/\alpha)$ and a union bound over $i$, gives
\begin{equation*}
    \Prob \left( \max_{i\in[n]}\left|\frac{1}{n-1} \sum_{j\ne i}Z_i^\top Z_j-\|\E[Z]\|_2^2\right|\leq C l_s\left(s+\sqrt{\frac dn}\right) \right) \geq 1 - \alpha/2.
\end{equation*}

Taking $t=c l_s$ in \eqref{eq:inner_tail} and choosing $c$ and $C_{\mathrm{clip}}$ sufficiently large gives
\begin{equation*}
    \Prob\left(\max_{i<j}|Z_i^\top Z_j|>R(s)\right)\leq\binom {n}{2} \exp(-c l_s)\leq\alpha/2.
\end{equation*}

On the intersection of these events, no inner product is clipped, so $\widehat h_i^{(s)}=u_i/(n-1)$. Since $U_s$ is the average of the rows,
\begin{equation*}
\max_i|\widehat h_i^{(s)}-U_s|\leq2\max_i|\widehat h_i^{(s)}-\|\E[Z]\|_2^2|\leq C l_s\left(s+\sqrt{d/n}\right).
\end{equation*}
A union bound argument then proves \eqref{eq:hajek_weight_tol}.
\end{proof}

\subsection{Proof of \Cref{prop:sparse-polytime}}
\label{app:sparse-polytime}
\begin{proof}
\noindent \underline{Privacy and computation.} Let $D, D'$ be datasets that are different at only index $k$. We have that 
\begin{align*}
|T(D)-T(D')|
&=\left|\max_{j\in[d]}\left|\frac1n\sum_{i=1}^n[X_{ij}]_{-B}^{B}\tanh(Y_i)\right|
-\max_{j\in[d]}\left|\frac1n\sum_{i=1}^n[X'_{ij}]_{-B}^{B}\tanh(Y'_i)\right|\right|\\
&\leq\max_{j\in[d]}\left|\frac1n\sum_{i=1}^n
\left\{[X_{ij}]_{-B}^{B}\tanh(Y_i)-[X'_{ij}]_{-B}^{B}\tanh(Y'_i)\right\}\right|\\
&=\frac1n\max_{j\in[d]}
\left|[X_{kj}]_{-B}^{B}\tanh(Y_k)-[X'_{kj}]_{-B}^{B}\tanh(Y'_k)\right|\\
&\leq\frac{2B}{n}.
\end{align*}
Since the sensitivity of $T$ is $2B/n$, the privacy guarantee therefore follows from \Cref{lemma:laplace-test} and \Cref{lemma:post-processing}. 

Computing the clipped coordinate means and taking the maximum over $d$ coordinates requires $O(nd)$ operations. 

\medskip 
    \noindent\underline{Clipping bias and concentration.}
We have that 
\begin{align}
    \left| \E \left[[X_j]_{-B}^{B}\tanh(Y) \right] - \E \left[X_j\tanh(Y) \right]\right| &\leq  \E[\{|X_j| - B\}_+] \notag \\
    &= \int _B^\infty \Prob(|X_j| > t ) \dint t \leq 2 \int_B^\infty \exp(-t^2/C_X^2) \dint t \notag \\
    &\leq \frac{C_X^2}{B} \exp(-B^2/C_X^2) = b_B 
    \label{eq:coordinate-bias}
\end{align}
where the tail bound in the second line is due to \Cref{assp:covariate-dist}.

Since 
\begin{equation*}
    \var \left([X_j]_{-B}^{B}\tanh(Y)\right)
\leq\E\left[([X_j]_{-B}^{B}\tanh(Y))^2\right]
\leq\E[X_j^2]\leq\overline\lambda,
\end{equation*}
by Bernstein's inequality and a union bound over the $d$ coordinates, we have \begin{align}
&\Prob\left(\max_{j\in[d]}\left|\frac1n\sum_{i=1}^n[X_{ij}]_{-B}^{B}\tanh(Y_i)
-\E\left[[X_j]_{-B}^{B}\tanh(Y)\right]\right|>q_B\right) \leq2d\exp\left(-\frac{nq_B^2}{2\{\overline\lambda+(2B/3)q_B\}}\right)
\leq\frac{\eta}{2}.
\label{eq:coordinate-concentration}
\end{align}

\medskip 
    \noindent\underline{Type-I error.} Under the null, we have $ \max_{j \in [d]} |\E \left[[X_j]_{-B}^{B}\tanh(Y) \right] |\leq b_B$ by \eqref{eq:coordinate-bias}. By \eqref{eq:coordinate-concentration}, we have that $T\leq b_B + q_B$ with probability at least $1-\eta/2$. Since $\Prob(W>\log(1/\eta))=\eta/2$, a union bound gives type-I error at most $\eta$.

\medskip 
    \noindent\underline{Type-II error.}
Suppose that $\|\beta\|_2 \neq 0$ and let $\beta = \|\beta\|_2 u$ and $V:=u^\top X$. With $m_\nu$ and $s_\vartheta$ as in \Cref{lemma:conditinal-mean}, it implies that 
\begin{align*}
    \max_{j \in [d]}\left| \E \left[[X_j]_{-B}^{B}\tanh(Y) \right] \right| &\geq \|\E[X\tanh(Y)]\|_\infty - b_B\\
    &\geq \frac{\left|u^\top \E[X \tanh(Y)] \right|}{\|u\|_1} - b_B\\
    &\geq \frac{c_\mu}{\sqrt{s^*}} \min\{\|\beta\|_2, 1\} -b_B. 
\end{align*}
where the first inequality is due to \eqref{eq:coordinate-bias} and the third uses \eqref{eq:transformed-mean-directional}. By a union bound with \eqref{eq:coordinate-concentration} and $\Prob(W<-\log(1/\eta))=\eta/2$, we have that 
\begin{equation}
    T + \frac{2B}{n \epsilon} W  \geq\frac{c_\mu}{\sqrt{s^*}}\min\{\|\beta\|_2,1\}-b_B-q_B-\frac{2B\log(1/\eta)}{n\epsilon}
    \label{eq:sparse-polytime-alt}
\end{equation}
with probability at least $1-\eta$. By the choices of $B, b_B$, and $q_B$, we have that 
\begin{equation*}
    b_B+q_B+\frac{2B\log(1/\eta)}{n\epsilon}
\lesssim\sqrt{\frac{\log(4d/\eta)}n}
+\frac{\sqrt{\log(2n)}\log(4d/\eta)}n
+\frac{\sqrt{\log(2n)}\log(1/\eta)}{n\epsilon}.
\end{equation*}
Under \eqref{eq:coordinate-power} with a sufficiently large constant $C>0$, the right-hand side of \eqref{eq:sparse-polytime-alt} exceeds the rejection threshold, and hence the test has type-II error at most $\eta$. 
\end{proof}

\section{Minimax rates for sparse mean testing}
\label{app:sparse-mean}
We can obtain the minimax rates for sparse mean testing by making some modifications to our existing arguments. Fix $C_Z \in (0, \infty)$ and define
\begin{equation*}
\mathcal{P}_{\mathrm{mean}}:=\left\{P:\|\E_P[Z]\|_2\leq1,\ \|Z-\E_P[Z]\|_{\psi_2}\leq C_Z\right\}.
\end{equation*}
For $1\leq s^*\leq d$ and $0<r<1$, let
\begin{equation*}
\mathcal{P}_0:=\{P\in\mathcal{P}_{\mathrm{mean}}:\E_P[Z]=0\},\qquad \mathcal P_1(s^*,r,1):=\left\{P\in\mathcal{P}_{\mathrm{mean}}:\|\E_P[Z]\|_2\geq r,\ \|\E_P[Z]\|_0\leq s^*\right\}.
\end{equation*}
The bound $\|\mathbb E_P[Z]\|_2\leq1$ can be replaced by $\|\mathbb E_P[Z]\|_2\leq M$ for any fixed $M>0$, with the constants in the results below then allowed to depend on $M$.

We define the minimax separation radius as
\begin{align*}
r^*(n,d,s^*,\epsilon,\delta,\mathcal P_{\mathrm{mean}},1):=\inf\Big\{r\in(0,1):\;&\exists\,(\epsilon,\delta)\text{-DP test }\phi\text{ such that}\\
&\sup_{P\in\mathcal P_0}P^{\otimes n}(\phi=1)+\sup_{P\in\mathcal P_1(s^*,r,1)}P^{\otimes n}(\phi=0)\leq \frac{1}{3} \Big\}.
\end{align*}

Recall the definitions of $\varrho_{s^*, d}$ and $\ell_{n, d, \epsilon}$ from \eqref{eq:separation-notation}. 
\begin{proposition}
    \label{prop:sparse-mean-rate}
    There exist constants $c_P$ and $c_\delta$, depending only on $C_Z$, such that, for every $s^*\in[d]$, if $\ell_{n,d,\epsilon}^3\varrho_{s^*,d}\leq c_P$ and $0\leq\delta\leq c_\delta\epsilon$, then
\begin{equation*}
 \varrho_{s^*,d} \lesssim  r^*(n,d,s^*,\epsilon,\delta,\mathcal P_{\mathrm{mean}},1) \lesssim \ell_{n,d,\epsilon}^3\varrho_{s^*,d}.
\end{equation*}
\end{proposition}

The minimax rates for sparse mean testing are therefore the same as those  for signal detection in linear and binary regression in \eqref{eq:lower-linear}.

Our test for the sparse regime is based on testing the sample mean of clipped sparse projections of the observations. Let $\mathcal K:=\{2^j:0\leq j\leq\lfloor\log_2d\rfloor\}\cup\{d\}$. For each $k\in\mathcal K$, let $\mathcal V_k$ be the union of $1/4$-nets of the unit spheres on all supports of size $k$, so that $|\mathcal V_k|\leq(9ed/k)^k$. Let $\eta \in (0, 1/4)$, $B := 1+C_Z\sqrt{2\log(2n^2)}$, and for a privacy budget $\epsilon_0$, consider the test 
\begin{equation}
    \mathbbm{1}\left\{A(D)+\frac{2B}{n\epsilon_0} W> \frac{2B\log(1/\eta)}{n\epsilon_0}\right\}, \quad A(D):=\max_{k\in\mathcal K}\left\{\max_{v\in\mathcal V_k} \left|\frac{1}{n}\sum_{i=1}^n[v^\top Z_i]_{-B}^B \right|-p_B-q_k\right\}, 
    \label{eq:sparse-mean-test}
\end{equation}
where $W \sim \operatorname{Lap}(1)$ and 
\begin{equation*}
     p_B:=\frac{2C_Z^2}{B-1}\exp\left\{-\frac{(B-1)^2}{2C_Z^2}\right\}, \qquad q_k:=B\sqrt{\frac{2\log(4|\mathcal K||\mathcal V_k|/\eta)}n}.
\end{equation*}

\begin{lemma}Suppose \(Z_1,\ldots,Z_n\) are i.i.d. from \(P\in\mathcal P_{\mathrm{mean}}\), with mean \(\mu\), and let \(n\geq2\), \(\epsilon_0>0\) and \(\eta\in(0,1/4)\). The test in \eqref{eq:sparse-mean-test} is \((\epsilon_0,0)\)-DP and has type-I error at most \(\eta\). There exists a universal constant \(C>0\) such that, whenever \(\|\mu\|_0\leq s^*\) and
\begin{equation}
\|\mu\|_2\geq C\left\{p_B+B\sqrt{\frac{s^*\log(ed/s^*)+\log(|\mathcal K|/\eta)}{n}}+\frac{B\log(1/\eta)}{n\epsilon_0}\right\},
\label{eq:sparse-mean-sep}
\end{equation}
its type-II error is at most \(\eta\).
    \label{prop:sparse-mean-upper}
\end{lemma}

\subsection{Proof of the lower bound}
\begin{proof}
    We consider the set of distributions $P_\mu=\mathcal N(\mu, \sigma_0^2I_d)$ with $\|\mu\|_2\leq1$ and $\sigma_0 := C_Z \sqrt{3/8}$, which belong to $\mathcal{P}_{\mathrm{mean}}$. We consider a prior on $\mu$ that chooses, uniformly at random, a support of size $s^*$ and independent signs, with non-zero entries $\pm r/\sqrt{s^*}$, and denote the resulting mixture by $\bar P_{r,s^*}^{(n)}$.

    For $\|\mu\|_2 = \|\mu'\|_2=r$, we have that 
    \begin{align*}
        \E_0 \left[ \frac{\dint P_\mu}{ \dint P_{0}}(Z) \frac{\dint P_{\mu'}}{\dint P_{0}}(Z) \right] &= \E_0 \left[ \exp \left\{\frac{(\mu + \mu')^\top Z - r^2}{\sigma_0^2} \right\} \right]\\
        &= \exp \left\{ \frac{1}{2 \sigma_0^2} \|\mu+\mu'\|_2^2- \frac{r^2}{\sigma_0^2} \right\} = \exp \left\{ \frac{\mu^\top \mu'}{\sigma_0^2} \right\}.
    \end{align*}

    By similar arguments to the proof of \Cref{lemma:gaussian-mixture}, we thus have that 
   \begin{equation}
    \KL(\bar P_{r,s^*}^{(n)}\parallel P_0^{\otimes n}) \leq\frac{(s^*)^2}{d}\left\{\cosh\left(\frac{nr^2}{ \sigma_0^2s^*}\right)-1\right\}.
    \label{eq:sparse-mean-mixture-KL}
    \end{equation}

    Let $r_{\mathrm{np}}$ and $r_{\mathrm{priv}}$ be as defined in \eqref{eq:r_np-def}. Similar calculations give, for sufficiently small $c>0$, that 
    \begin{equation*}
        \KL(\bar P_{r,s^*}^{(n)}\parallel P_0^{\otimes n})\leq 
        \begin{cases}
             c^4,&r=c\min\{ \sigma_0r_{\mathrm{np}},1\},\\
             c^4/(n\epsilon^2),&r=c\min\{ \sigma_0r_{\mathrm{priv}},1\}.
        \end{cases}
\end{equation*}

By considering $P_0$ and $P_{r (1, 0, \dots, 0)}$, which have a KL-divergence of $r^2/(2\sigma_0^2)$, we also have a two-point hypothesis comparison. We can therefore apply \Cref{prop:general_lower_bound} with 
\begin{equation*}
    \alpha\lesssim r, \qquad \tau=0, \qquad \omega \lesssim
\begin{cases}
c^2/\sqrt n,&r= c\min\{ \sigma_0 r_{\mathrm{np}}, 1\},\\
c^2/(n\epsilon),&r=c \min\{\sigma_0 r_{\mathrm{priv}}, 1\}
\end{cases}
\end{equation*}
and taking sufficiently small constants to satisfy the error bounds as in \Cref{app:lower-assembly} gives the lower bound. 
\end{proof}

\subsection{Proof of the upper bound}
\begin{proof}[Proof]
    We can construct a test for the dense regime using \Cref{lemma:single_scale_test} and applying it over a dyadic grid of scales as in \Cref{lemma:private-multiscale-reduction}, since its proof and that of \Cref{lemma:typical_set_prob} only require $Z$ to be sub-Gaussian with a bound on $\|\E[Z]\|_2$ under the alternative. This gives a $(\epsilon_0, 0)$-DP test with a detection boundary of 
    \begin{equation}
         C\ell_{n,d,{\epsilon_0}}^3\left(\frac{d^{1/4}}{\sqrt n} +\frac{d^{1/4}}{n^{3/4}\epsilon_0^{1/2}}+\frac1{n\epsilon_0}\right).
         \label{eq:dense-mean-rate}
    \end{equation}

We can therefore run this test and the test in \eqref{eq:sparse-mean-test} each with privacy budget $\epsilon_0=\epsilon/2$ and $\eta = 1/9$, and reject if either of the tests reject. The type-I error is at most $2\eta$ and the type-II error is at most $\eta$ under either of the separation conditions in \eqref{eq:dense-mean-rate} or \eqref{eq:sparse-mean-sep}. Taking the minimum of the two separation conditions gives a test with rate $C\ell_{n,d,\epsilon}^3\varrho_{s^*,d}$, and a total error of at most $1/3$.
\end{proof}

\begin{proof}[Proof of \Cref{prop:sparse-mean-upper}]
\medskip \noindent\underline{Privacy.} Replacing one observation changes $\bar{Z}_v:=n^{-1} \sum_{i=1}^n [v^\top Z_i]_{-B}^B$ by at most $2B/n$. As in the proof for \eqref{eq:adaptive-covariance-test}, taking absolute values, subtracting fixed thresholds and maximising over $(k,v)$ preserves this sensitivity bound. Thus $|A(D)-A(D')|\leq2B/n$ for neighbouring datasets, and \Cref{lemma:laplace-test,lemma:post-processing} imply that the test in \eqref{eq:sparse-mean-test} is $(\epsilon_0,0)$-DP.

    \medskip\noindent\underline{Clipping bias and concentration.} Since $|v^\top\mu|\leq1$, the clipping argument in \eqref{eq:coordinate-bias} gives
\begin{align}
 \left|\E\left[[v^\top Z]_{-B}^B\right]-v^\top\mu\right|
 &\leq\int_B^\infty\Prob(|v^\top Z|>t) \dint t\notag\\
 &\leq2\int_B^\infty\exp\left\{-\frac{(t-1)^2}{2C_Z^2}\right\} \dint t
 \leq p_B\lesssim n^{-2}.
 \label{eq:sparse-mean-bias}
\end{align}
Hoeffding's inequality and a union bound, give
\begin{equation}
 \Prob\left(\exists k\in\mathcal K,\ v\in\mathcal V_k: \left|\frac{1}{n} \sum_{i=1}^n[v^\top Z_i]_{-B}^B -\E\left[[v^\top Z]_{-B}^B\right]\right|>q_k\right)
 \leq\sum_{k\in\mathcal K}2|\mathcal V_k|\exp\left(-\frac{nq_k^2}{2B^2}\right)=\frac\eta2.
 \label{eq:sparse-mean-concentration}
\end{equation}

\medskip\noindent\underline{Type-I error.} Under the null, \eqref{eq:sparse-mean-bias} and \eqref{eq:sparse-mean-concentration} imply that $A(D)\leq0$ except on an event of probability at most $\eta/2$. Since $\Prob(W>\log(1/\eta))=\eta/2$, the type-I error of \eqref{eq:sparse-mean-test} is at most $\eta$.

\medskip\noindent\underline{Type-II error.} Suppose $\|\mu\|_2>0$ and $\|\mu\|_0\leq s^*$. Choose $k^*\in\mathcal K$ with $s^*\leq k^*\leq2s^*$. There exists $v^*\in\mathcal V_{k^*}$ such that $\|v^*-(\mu/\|\mu\|_2)\|_2\leq1/4$; for this choice of $v^*$, we have that $|(v^*)^\top\mu|\geq 3\|\mu\|_2/4$. On the event in \eqref{eq:sparse-mean-concentration}, we have 
\begin{equation*}
 A(D)\geq \left|\frac{1}{n} \sum_{i=1}^n [(v^*)^\top Z_i]_{-B}^B \right|- p_B -q_{k^*}
       \geq\frac{3}{4}\|\mu\|_2- 2p_B-2q_{k^*}.
\end{equation*}
Since $\Prob(W<-\log(1/\eta))=\eta/2$, the type-II error is at most $\eta$ provided that
\begin{equation}
 \|\mu\|_2\geq4\left(p_B+q_{k^*} +\frac{2B\log(1/\eta)}{n\epsilon_0}\right).
 \label{eq:sparse-mean-power}
\end{equation}
The bound in \eqref{eq:net-card-bounds} implies that $\log|\mathcal V_{k^*}|\lesssim s^*\log(ed/s^*)$, so \eqref{eq:sparse-mean-power} is satisfied under \eqref{eq:sparse-mean-sep}. 
\end{proof}

\section{Additional details on the numerical experiments}
\label{app:numerical}
In this appendix, we give further details on sensitivity to tuning parameters and implementation of the different methods during the numerical experiments. 

\medskip \noindent
\textbf{Sensitivity of \Cref{alg:polytime_test} to tuning parameters.}
The implementation in \Cref{alg:polytime_test} involves three tuning constants: $C_{\mathrm{sep}}$ determines the detection radius, $C_{\mathrm{clip}}$ controls the kernel truncation level, and $C_{\mathrm{tol}}$ controls the concentration tolerance used in the local H\'ajek estimator. Our default choices are $C_{\mathrm{sep}}=5$, $C_{\mathrm{clip}}=0.5$ and $C_{\mathrm{tol}}=0.05$. We assess sensitivity by varying these constants jointly over $C_{\mathrm{sep}}\in\{3,5,7\}$, $C_{\mathrm{clip}}\in\{0.25,0.5,0.75\}$ and $C_{\mathrm{tol}}\in\{0.025,0.05,0.075\}$. We use the Gaussian linear and exponential regression settings described in the main text, fixing $n=4000$. The Gaussian experiment uses $d=500$ and $\epsilon=1$, while the exponential experiment uses $d=50$ and $\epsilon=4$. We report the estimated sum of the type-I and -II errors over 250 repetitions in \Cref{tab:sensitivity-constants}. In these pilot experiments, performance is relatively stable across the values of $C_{\mathrm{tol}}$ considered, but can be sensitive to $C_{\mathrm{sep}}$ and $C_{\mathrm{clip}}$, particularly to larger values of $C_{\mathrm{sep}}$ in exponential regression. 

\begin{table}[ht]
\centering
\begin{tabular}{rrrrrrrr}
  \toprule
 & & \multicolumn{3}{c}{Gaussian linear regression} & \multicolumn{3}{c}{Exponential regression} \\
\cmidrule(lr){3-5}\cmidrule(lr){6-8}
$C_{\mathrm{sep}}$ & $C_{\mathrm{clip}}$ & $C_{\mathrm{tol}}=0.025$ & $C_{\mathrm{tol}}=0.05$ & $C_{\mathrm{tol}}=0.075$ & $C_{\mathrm{tol}}=0.025$ & $C_{\mathrm{tol}}=0.05$ & $C_{\mathrm{tol}}=0.075$ \\ 
  \midrule
3 & 0.25 & 0.184 & 0.240 & 0.200 & 0.008 & 0.016 & 0.024 \\ 
  3 & 0.50 & 0.344 & 0.412 & 0.340 & 0.036 & 0.036 & 0.080 \\ 
  3 & 0.75 & 0.484 & 0.472 & 0.428 & 0.116 & 0.108 & 0.112 \\ 
  5 & 0.25 & 0.032 & 0.032 & 0.056 & 0.068 & 0.036 & 0.048 \\ 
  5 & 0.50 & 0.140 & 0.144 & 0.152 & 0.044 & 0.028 & 0.016 \\ 
  5 & 0.75 & 0.292 & 0.312 & 0.300 & 0.044 & 0.052 & 0.040 \\ 
  7 & 0.25 & 0.144 & 0.092 & 0.088 & 0.936 & 0.960 & 0.968 \\ 
  7 & 0.50 & 0.052 & 0.052 & 0.032 & 0.780 & 0.840 & 0.836 \\ 
  7 & 0.75 & 0.088 & 0.084 & 0.100 & 0.788 & 0.784 & 0.820 \\ 
   \bottomrule
\end{tabular}
\caption{Sensitivity of \Cref{alg:polytime_test} to the tuning parameters $C_{\mathrm{sep}}$, $C_{\mathrm{clip}}$ and $C_{\mathrm{tol}}$, in Gaussian linear regression ($n=4000, d=500, \epsilon=1$) and exponential regression ($n=4000, d=50, \epsilon=4$). Entries give the estimated sum of type-I and type-II error, using 250 null and 250 alternative repetitions per cell} 
\label{tab:sensitivity-constants}
\end{table}

\medskip \noindent \textbf{Sensitivity of the coordinate test.} To quantify the robustness of the coordinate test in \Cref{prop:sparse-polytime} to misspecification in the covariate population quantities $C_X$ and $\overline{\lambda}$ in the rejection threshold, we run the test with $C_X = c_B \sqrt{8/3}$ and $\bar{\lambda} = c_\lambda$, since the true population quantities are $\sqrt{8/3}$ and $1$ respectively. \Cref{tab:coordinate-sensitivity} suggests some robustness to scaling by 0.1 to 5 in both parameters. 

\begin{table}[ht]
\centering
\begin{tabular}{rrrrrr}
\toprule
 & \multicolumn{5}{c}{$c_B$} \\
\cmidrule(lr){2-6}
$c_\lambda$ & $0.10$ & $0.25$ & $1.00$ & $5.00$ & $10.00$ \\
\midrule
0.1 & 0.000 & 0.108 & 0.004 & 0.000 & 0.000 \\
0.25 & 0.000 & 0.000 & 0.004 & 0.000 & 0.000 \\
1 & 0.000 & 0.000 & 0.000 & 0.000 & 0.000 \\
5 & 0.000 & 0.000 & 0.000 & 0.000 & 0.032 \\
10 & 1.000 & 0.000 & 0.000 & 0.000 & 0.432 \\
\bottomrule
\end{tabular}
\caption{Sensitivity of the coordinate test, for Gaussian linear regression with $n=4000$, $d=500$, $\epsilon=1$, $\|\beta\|_0=1$. We vary $C_X=c_B\sqrt{8/3}$ and $\overline\lambda=c_\lambda$. Entries give the estimated sum of type-I and type-II error, using 250 null and 250 alternative repetitions per cell}
\label{tab:coordinate-sensitivity}
\end{table}

\medskip \noindent
\textbf{Implementation of dpHSIC.} For dpHSIC, we use Gaussian kernels with bandwidth $\sigma_X = \sqrt{d}$ for the covariates, and $\sigma_Y=1$ for the response, following the default dimensional scaling in the implementation accompanying \cite{kim2026permutation}. We use 99 permutations. \Cref{tab:dpHSIC_sensitivity} reports sensitivity to bandwidth choice at $n=4000$; our fixed choice gives the third-lowest estimated total error among the configurations considered.

\begin{table}[h!]
\centering
\begin{tabular}{rrrrr}
  \toprule
 & \multicolumn{4}{c}{$\sigma_Y$} \\
\cmidrule(lr){2-5}
$\sigma_X/\sqrt{d}$ & $0.5$ & $1.0$ & $2.0$ & $4.0$ \\ 
  \midrule
0.50 & 0.960 & 0.980 & 1.004 & 0.960 \\ 
  1.00 & 0.948 & 0.832 & 0.724 & 0.856 \\ 
  2.00 & 0.932 & 0.872 & 0.816 & 0.904 \\ 
   \bottomrule
\end{tabular}
\caption{Sensitivity of dpHSIC to kernel bandwidths, for Gaussian linear regression with $n=4000$, $d=500$, and $\epsilon=1$. Entries give the estimated sum of type-I and type-II error, using 250 null and 250 alternative repetitions per cell} 
\label{tab:dpHSIC_sensitivity}
\end{table}

\medskip \noindent
\textbf{Implementation of ToT-F.} We randomly partition the observations into $M=\lfloor n/(d+1)\rfloor$ balanced blocks, the largest number for which each block contains more than $d$ observations, as required for the regression F-test without an intercept. Within each block, we test $\beta=0$ using the F-test at significance level $0.25$, and count the rejections. Following \cite{pmlr-v202-kazan23a}, we add Tulap (truncated-uniform-Laplace) noise to this count and reject when the noisy count exceeds the $0.95$ quantile of the convolution of the $\operatorname{Binomial}(M,0.25)$ null distribution with the noise distribution. The sensitivity analysis in \Cref{tab:aggregate_blocks_sensitivity} supports using the largest feasible number of blocks in this setting.

\begin{table}[h!]
\centering
\begin{tabular}{rrr}
  \toprule
 & \multicolumn{2}{c}{Method} \\
\cmidrule(lr){2-3}
Blocks & ToT-F & PB-LRT \\ 
  \midrule
1 & 0.912 & 0.968 \\ 
  2 & 0.840 & 0.976 \\ 
  3 & 0.732 & 1.008 \\ 
  4 & 0.388 & 0.960 \\ 
  5 & 0.232 & 0.944 \\ 
  6 & 0.152 & 0.908 \\ 
  7 & 0.096 & 0.920 \\ 
   \bottomrule
\end{tabular}
\caption{Sensitivity of ToT-F and PB-LRT to the number of aggregation blocks, for Gaussian linear regression with $n=4000$, $d=500$, and $\epsilon=1$. Entries give the estimated sum of type-I and type-II error, using 250 null and 250 alternative repetitions per cell.} 
\label{tab:aggregate_blocks_sensitivity}
\end{table}

\medskip \noindent
\textbf{Implementation of PB-LRT.} Following Section 4 of \cite{pena2024uncertainty}, we compute the Gaussian regression log-likelihood ratio within each block, for $M$ blocks, clip it to $[0,U]$, and add Laplace noise with scale $U/(M\epsilon)$ to the average. We use the same block partition rule as in ToT-F. We set $U$ to the largest blockwise $0.99$ null quantile and estimate the $0.95$ quantile of the resulting private null statistic using 10,000 simulation draws, following the calibration approach in their Section 4.3.

\section{Auxiliary lemmas}
\begin{lemma}[Inner products of independent sub-Gaussian vectors]
\label{lemma:sub-Gaussian-inner-product}
Let $U,V\in\R^d$ be independent, centred random vectors such that $\|U\|_{\psi_2} \leq a$ and $\|V\|_{\psi_2} \leq b$. There exists an absolute constant $C>0$ such that 
for every $p\geq1$ and $t\geq1$, 
\begin{equation*}
\|U^\top V\|_{L^p} \leq Cab(\sqrt{dp}+p),
\end{equation*}
and consequently
 $$\Prob\left(| U^\top V|>Cab(\sqrt{dt}+t)\right) \leq e^{-t}.$$
\end{lemma}

\begin{proof}
Conditional on $V$, the scalar $U^\top V$ is sub-Gaussian with $\|U^\top V\|_{\psi_2\mid V} \leq a\|V\|_2$. Hence the standard sub-Gaussian moment bound gives
\begin{equation}
 \|U^\top V\|_{L^p(U\mid V)}
 \leq Ca\sqrt p\,\|V\|_2.
 \label{eq:typical-conditional_moment}
\end{equation}
By Proposition~6.2.1 of \cite{vershynin2026},
\begin{equation*}
 \Prob\{\|V\|_2>Cb(\sqrt d+s)\}\leq2e^{-s^2}, \qquad s\geq0.
\end{equation*}
Taking the $L^p$ norm with respect to $V$ in
\eqref{eq:typical-conditional_moment} therefore gives
\begin{align*}
 \|U^\top V\|_{L^p}
 &\leq Ca\sqrt p\,\|\|V\|_2\|_{L^p}\\
 &\leq Cab\sqrt p(\sqrt d+\sqrt p)
 =Cab(\sqrt{dp}+p).
\end{align*}
The tail bound follows from Markov's inequality:
\begin{equation*}
\Prob \left( |U^\top V| > eCab(\sqrt{dt} + t) \right) \leq e^{-t}. 
\end{equation*} 
\end{proof}

\begin{lemma}[{\cite{mcshane1934extension}}]
\label{thm:mcshane-extension}
Let $(\mathcal X,\mathsf d)$ be a metric space, let $\varnothing\neq\mathcal C\subseteq\mathcal X$, and suppose that $f:\mathcal C\to\mathbb R$ is $L$-Lipschitz:
\begin{equation*}
    |f(z)-f(z')|\leq L\mathsf d(z,z') \qquad z,z'\in\mathcal C.
\end{equation*}
Then
\begin{equation*}
    \widehat f(x):= \inf_{z\in\mathcal C} \left\{f(z)+L\mathsf d(x,z)\right\}, \qquad x\in\mathcal X,
\end{equation*}
is an $L$-Lipschitz extension of $f$ to $\mathcal X$. In particular, $\widehat f(x)=f(x)$ for all $x \in \mathcal{C}$ and 
\begin{equation*}
    |\widehat f(x)-\widehat f(x')|
    \leq L\mathsf d(x,x'), \qquad x,x'\in\mathcal X.
\end{equation*}
\end{lemma}

\begin{lemma}[Sequential composition of DP; e.g.~\citealp{dwork2014algorithmic}, Theorem~3.16]
\label{lemma:composition}
Let $\epsilon_1,\epsilon_2>0$ and $\delta_1,\delta_2\ge0$. If $M_1$ and $M_2$ are $(\epsilon_1,\delta_1)$-DP and $(\epsilon_2,\delta_2)$-DP mechanisms using independent randomisation, then releasing the outputs of both mechanisms is $(\epsilon_1+\epsilon_2,\delta_1+\delta_2)$-DP.
\end{lemma}

\end{document}